\documentclass[11pt]{article}
\usepackage{amsmath}
\usepackage{amsfonts}
\usepackage{amssymb}
\usepackage{mathrsfs}
\usepackage{epsfig}
\usepackage{fancybox}
\usepackage{dsfont}
\usepackage{bbm}
\usepackage{bm}
\usepackage{multirow}
\usepackage[dvipsnames]{xcolor}
\usepackage{float}
\usepackage{makecell} 
\usepackage{boldline}
\usepackage{diagbox}
\usepackage{threeparttable}
\usepackage{adjustbox}
\usepackage{booktabs}
\usepackage[noend]{algpseudocode}
\usepackage{caption}
\usepackage[title]{appendix}
\usepackage[colorlinks=true]{hyperref}
\hypersetup{urlcolor=black,citecolor=blue,linkcolor=blue}
\allowdisplaybreaks[4]

\usepackage{tikz}

\renewcommand{\Box}{\framebox{\rule{0.3em}{0.0em}}}

 \def\gph {{\rm gph}}

\newtheorem{theorem}{Theorem}[section]
\newtheorem{theorem*}{Theorem}[subsubsection]
\newtheorem{lemma}{Lemma}[section]
\newtheorem{proposition}{Proposition}[section]
\newtheorem{example}{Example}[section]
\newtheorem{remark}{Remark}[section]
\newtheorem{definition}{Definition}[section]
\newtheorem{assumption}{Assumption}[section]

\newcommand{\setd}{{ d \kern -.15em l}}
\newcommand{\hatsetd}{ d \hat{\kern -.15em l }}
\newcommand{\dd}{\mathsf {d\kern -0.07em l}}

\newcommand{\bgeqn}{\begin{eqnarray}}
\newcommand{\edeqn}{\end{eqnarray}}
\newcommand{\bgeq}{\begin{eqnarray*}}
\newcommand{\edeq}{\end{eqnarray*}}

\newcommand{\R}{{\rm I\!R}}

\newcommand{\inmat}[1]{\mbox{\rm {#1}}}

\newcommand{\Z}{{\cal Z}}
\newcommand{\F}{{\cal F}}

\newcommand{\be}{\begin{equation}}
\newcommand{\ee}{\end{equation}}

\newcommand{\bdz}{{\bm z}}

\def\bbe{{\mathbb{E}}} 

\renewcommand{\Box}{\hfill \rule{2.3mm}{2.3mm}}

\renewcommand{\Box}{\framebox{\rule{0.3em}{0.0em}}}

\newcommand{\N}{{\cal N}}

\def\gph {{\rm gph}}

\def\bbe{{\mathbb{E}}} 

\def\Prob{{\rm{Prob}}}

\title{
Statistical Robustness of Kernel Learning Estimator
with Respect to Data Perturbation
}

\author{
 Sainan Zhang\footnote{Department of Systems Engineering and Engineering Management, The Chinese University of Hong Kong, Shatin, N.T., Hong Kong. Email:  snzhang.m@gmail.com.},~  Huifu Xu\footnote{Department of Systems Engineering and Engineering Management, The Chinese University of Hong Kong, Shatin, N.T., Hong Kong. Email: hfxu@se.cuhk.edu.hk.}~ and~ Hailin Sun\footnote{School of Mathematical Sciences and Institute of Mathematics, Nanjing Normal University, Nanjing 210023, China, Email: hlsun@njnu.edu.cn} 
}

\date{\today}

\begin{document}

\maketitle

\noindent
\begin{abstract}
{Inspired by the recent work \cite{guo2023statistical} on the statistical robustness of empirical risks in reproducing kernel Hilbert space
(RKHS) where the training data are potentially perturbed or even corrupted,
we take a step further in this paper to investigate the 
statistical robustness of the kernel learning estimator (the regularized empirical risk minimizer  or stationary point).
 We begin by deriving qualitative statistical robustness 
of the estimator of the regularized empirical risk minimizer for
a broad class of convex cost functions
when all of the training data are potentially perturbed under 
some topological structures,
and then move on  to consider the quantitative statistical robustness of the stationary solution for a
specific case that the cost function is 
continuously differentiable but not necessarily convex. 
In the latter case,
we derive the first-order optimality condition of the regularized expected risk minimization problem, which is essentially a stochastic variational inequality problem (SVIP) in RKHS, and then use the SVIP as a platform to investigate local and global Lipschitz continuity of the 
stationary solution
against perturbation of the probability distribution under the Fortet-Mourier metric.
A crucial assumption in the analysis is that the perturbed data are independent and identically distributed (iid). In some practical applications, this assumption may not be fulfilled when a small proportion of perceived data is seriously perturbed/contaminated. In this case, we use the influence function to investigate the impact of single data perturbation on 
the expected risk minimizer.
Differing from Steinwart and Christmann~\cite[Chapter 10]{StC08}, we concentrate on constrained expected risk minimization problems. The research is essentially down to the derivation of the implicit function theorem of the SVIP in RKHS. Finally, we illustrate our theoretical analysis with a couple of academic examples.
}
\end{abstract}


\section{Introduction}
We consider a stochastic program in machine learning where the uncertainties are characterized by the random variables in the input-output spaces with probability distributions.
Without knowing the exact probability distribution of the variables,
the stochastic program uses the discretization technique with a set of input-output pairs (training data).
The task in supervised learning is to solve a stochastic program with the resulting empirical distribution over a hypothesise class
and obtain the empirical risk minimizer (ERM).
This 
raises a question as to whether 
ERM converges to its true counterpart as sample size increases, 
which is known as consistency or  
learnability.
The former property can be characterized by the convergence of empirical risks to the expected counterparts uniformly for all functions in hypothesis class,
see \cite{BEHW89},
\cite[Chapter 19]{Ant99} and 
\cite{VaC71}.
Qi et al.~\cite{QCLP22} establish the consistency and convergence rates of stationary solutions and values
of a class of coupled nonconvex and nonsmooth empirical risk minimization problems.
The learnability
relates to the learning rules and 
sensitivity to perturbations in the training set.
For example,
Shalev-Shwartz et al.~\cite{SSS10}  and 
Mukherjee et al.~\cite{MNPR06} state that stability is a necessary and sufficient condition for learnability.
Roger and Wagner~\cite{RoW78} discuss that the sensitivity of a learning algorithm with respect to (w.r.t.) small perturbation in the model controls the variance of the leave-one-out estimate.

A conventional assumption in the stability analysis is that the samples are independent and identically distributed (iid) generated from the true probability measure. 
In practice, however, data are often perturbed for various reasons 
such as data contamination and data/distribution shifts from past/present to future. 
For instance,
in data-driven problems such as biology 
and medical experiments,
data may be corrupted with
outliers (see \cite{LeL20}) which do not follow the true distribution. 
The least squares estimators of Lasso are not stable in corrupted environments (see \cite{Tib96}).
$M$-estimators using loss functions such as Huber or $L_1$ loss are not robust to outliers (see \cite{HuR09}).
Consequently,
 it is necessary to investigate 
 whether empirical risks and 
empirical risk minimizers 
 obtained with 
 perceived 
 contaminated data
 are useful.
 This 
 kind of research 
 stems from Tukey~\cite{Tuk60,Tuk62}
 and  Hampel \cite{Ham71,Ham74},
 and popularized by others particularly the monographs of  Huber \cite{Hub81} and Huber and Ronchetti \cite{HuR09}.
 Steinwart and Christmann \cite{StC08} extend the research
to support vector machine
by studying the sensitivity of
learning estimator using the concept of influence function 
in single data perturbation. 
They also propose to investigate the quality of learning estimators by 
estimating the difference between the distributions of the 
empirical learning estimators obtained with contaminated data and  ideal pure data under the Prokhorov metric 
in the case of 
all data 
being potentially contaminated. 
In a more recent development,
Liu and Pang \cite{LiP23} propose an interval 
CVaR (In-CVaR) approach to reduce the effect of outliers on the kernel learning estimator.
Data shift is another important example of data perturbation. It occurs when
(a) the true probability distribution is known or subjective, but the perceived data are often perturbed due to measurement errors or 
occurrence of unexpected events, (b) the validation
data deviate from the past and present data (\cite{besbes2022beyond}),
(c) information about agent's
subjective probability of the scenarios of the state of nature in decision analysis is poor, which may lead to model uncertainty \cite{maccheroni2006ambiguity}. 
 In either case, it is necessary to estimate the difference between the kernel learning estimator based on perceived data and the one based on real data to ensure the usefulness of the estimator.

Guo et al.~\cite{guo2023statistical} seem to be the first to investigate  qualitative and quantitative statistical robustness of the 
optimal value of the empirical risk minimization problem 
in reproducing kernel Hilbert space (RKHS).
Specifically, they derive sufficient conditions under which the expected risk
 changes stably (continuously) against small perturbation of the probability distributions of
 the underlying random  variables and demonstrate how the cost function and kernel
 affect the stability. Moreover, they
 examine the difference between
    laws of the statistical estimators of the expected
    optimal loss based on pure data and contaminated
    data using the Prokhorov metric and Kantorovich metric,
    and derive
    some  asymptotic qualitative and
    non-asymptotic quantitative statistical robustness results.
    They also identify appropriate metrics under which the statistical estimators
    are uniformly asymptotically consistent.

In practice, it might be more important
to consider the optimal solution 
(regularized empirical risk minimizer/stationary solution) than the optimal value (regularized empirical risk) because the former is 
the focus of the empirical risk minimization problem.
In this paper, we complement the research 
of \cite{guo2023statistical}
to study the former
 in both 
 all data perturbation and single data perturbation cases.
 The former is 
 instrumental 
 as it allows us to assess the relative influence of individual data 
 on the value of the kernel learning estimator.
 The latter is also important in practice because we may not have specific information to identify 
 the structure of the 
 {\color{black} perturbed}
 data
 and consequently we 
 treat all the data 
 as potentially 
perturbed.
In this case, it requires 
us to quantify
 the extent of data perturbation under which the resulting laws of statistical estimators remain stable. 
The main 
contributions of this paper 
can be summarized as follows.

First, we derive 
 asymptotic  qualitative statistical robustness 
of the estimator of the regularized empirical risk minimizer for
a broad class of cost functions (which are not necessarily differentiable)
when all of the training data are potentially perturbed under 
some topological structure.
Specifically, we show that 
the difference between 
the laws of kernel learning estimators  based on 
real data and 
perceived data under the Prokhorov metric
is small when the distance between the 
true 
probability distribution underlying the real data and the true 
probability distribution underlying the perceived data is small and the sample size is sufficiently large (Theorem \ref{t-ML-SR-QL-SR}).
This kind of research is motivated 
by Steinwart and Christmann in \cite[Chapter 10]{StC08} for
support vector machine 
although they do not go into detailed discussions in theory as presented in this paper.

Second, we derive non-asymptotic quantitative statistical robustness 
of the estimator of 
the stationary solution
for
a specific class of cost functions which 
are 
continuously differentiable but not necessarily convex, which gives rise to an explicit relationship between the distance of 
two statistical estimators and the distance between the true probability distributions 
underlying the real data and perceived data 
independent of sample size. To this end, 
 we derive the first-order optimality condition which is in essence 
a stochastic variational inequality problem (SVIP) in RKHS (Theorem~\ref{eq:first-order-optimal}).
Under some moderate conditions, 
we show the local and global Lipschitz continuity of the 
stationary solution mapping
over a relatively compact set of probability measures (Theorem~\ref{thm:lip_solution-lambda}). 
The stability results enable us to derive the desired non-asymptotic quantitative statistical robustness which
guarantees that 
difference between laws of
kernel learning estimators
based on 
pure data and 
perturbed
data under the Kantorovich metric
is linearly bounded by the Kantorovich distance  of the original probability distributions generating
the two sets of data (Theorem~\ref{t:QSR_SGE}, $c$ is twice continuously differentiable but not necessarily convex).
The adoption of the Kantorovich metric 
allows us to derive a linear error bound as opposed 
to the Prokhorov metric envisaged by 
 Steinwart and Christmann in \cite[Chapter 10]{StC08}.
While this kind 
of research 
has been conducted by Jiang and Li \cite{JiL22} for SVIP and by Guo and Xu \cite{GuX23}
for the stochastic generalized equation in finite-dimensional spaces,
it is the first  
time to be presented in the context of 
machine learning.
One of the new challenges 
that we have to tackle is that the decision variable 
depends on random variables.
In the case when the kernel learning estimator
lies in the interior of the set of 
feasible solutions,
we demonstrate quantitative statistical robustness results which 
  does not require 
 second order continuous differentiability of the cost function 
  (Theorem~\ref{thm:Qunt-SR-New}, $c$ is continuously differentiable and convex). 
 We have managed this result by 
 deriving 
 a new implicit function theorem (Theorem~\ref{thm:implicit_thm})
 extended from Dontchev and Rockafellar \cite[Theorem 1H.3]{DoR09}.

Third,  for single data perturbation,
we take a step forward from the existing
work of Steinwart and Christmann \cite[Chapter 10]{StC08} on the unconstrained minimization problem to the
constrained minimization problem.
This is essentially down to derive the implicit function theorem for the SVIP  in RKHS.
We propose to use
proto-derivative to
tackle the differentiation of the normal cone,
and discuss sufficient conditions under which the influence function is bounded.

The rest of the paper is organized as follows.
In Section~\ref{Sec:2}, we  set up
the regularized empirical risk minimization 
models, 
recall 
basic ideas of statistical robustness
and 
$\psi$-weak topology.
In Section~\ref{sec:Quali-SR-KLE}, we
derive qualitative statistical robustness 
of the estimator of the regularized empirical risk minimizer for
a broad class of cost functions
when all of the training data are potentially perturbed.
In Section~\ref{sec:Quant-SR-KLE}, 
we take a step further to
discuss
quantitative statistical robustness 
of the estimator of the regularized empirical risk minimizer/stationary solution
for a more specific class of 
empirical risk minimization problem where 
the cost function is continuously differentiable.
In Section \ref{Sec:Single-data-IF}, we
 derive
an explicit form 
and boundedness 
of the influence function. 
In Section~\ref{sec:numerical tests}, we report some numerical tests on the established theoretical results for both all data perturbation and single data perturbation, and finally
Section~\ref{sec:conclusion} concludes with some remarks.

\section{Problem statement and preliminaries}
\label{Sec:2}

In this section, we give a precise mathematical description of 
the empirical risk minimization
model, 
followed by 
a motivation for statistical robustness
and a recall of some related basic results in applied probability.

\subsection{The expected risk minimization model}
In the framework of learning,
there are two sets of variables:
an input variable ${\bm x}\in X$
and
 an output variable 
 $y\in Y$,
where $X \subset \R^n$ and 
$Y\subset \R$.
Let
$Z:= X\times Y$.
The relation between an input ${\bm x}\in X$ and
an output $y\in Y$ is described by a probability distribution
 $P({\bm x}, y)$.
Let $\F\subset {\cal H}(X)$ be a convex set of measurable functions, where
${\cal H}(X)$  is 
known as hypotheses space.
The problem of learning is down to 
look for the oracle in ${\cal F}$
such that
expected risk 
is
minimized, i.e., 
\begin{equation}
\label{eq:ML-Rf-min}
\min_{{\bm f}\in \F} \;R_{P}({\bm f}) :=\bbe_P[c({\bm z},{\bm f}({\bm x}))],
\end{equation}
where $c(\cdot,\cdot):Z\times \R \rightarrow \R_+$ denotes the loss function,
${\bm z}=({\bm x},y):\Omega \rightarrow Z$ is random vector defined  on an atomless probability space $(\Omega,{\cal B},\mathbb{P})$,
$P= \mathbb{P}\circ {\bm z}^{-1}$ is the probability
measure on $Z$ induced by ${\bm z}$,
 $\bbe_P[c({\bm z},{\bm f}({\bm x}))]$ denotes
 the 
expected loss.

Typically,
the cost  function $c(\cdot,\cdot)$ is
determined by learners.
For instance,
in various 
regression and support vector machine models (see e.g. \cite{SKE19,SSS10}),
     $c({\bm z},{\bm f}({\bm x}))=L(y,{\bm f}({\bm x}))$,
   where $L:\R\times \R \to \R_+$
   takes different forms such as 
    \begin{itemize}
    \item[(i)] 
$L(t,y) = \max\{0, 1-t\cdot y\}$ (the hinge loss
function) which is not differentiable 
(satisfied in Theorem~\ref{t-ML-SR-QL-SR});
    \item[(ii)] 
    $L(t,y)=(t-y)^2/2$ for $|t-y|\leq 1$ and $L(t)=|t-y|-1/2$ for $|t-y|>1$ (the Huber loss
function), which is continuously differentiable but not twice differentiable
(satisfied in Theorem~\ref{thm:Qunt-SR-New});  
  \item[(iii)] 
     $L(t,y)=\frac{1}{2}(t-y)^2$ (the squared loss function)
    and  $L(t,y)=\log (1+e^{-t-y})$ (log-loss function),
    which are $k$-th differentiable for all $k\geq 1$ (satisfied in Theorem~\ref{t:QSR_SGE}).
    \end{itemize}
A key ingredient of model (\ref{eq:ML-Rf-min}) is
the hypothesis space ${\cal H}(X)$.
In this paper, we concentrate on the reproducing kernel Hilbert space, see e.g., \cite{CuZ07, guo2023statistical}.

\begin{definition}
[Reproducing kernel Hilbert space (RKHS)]
\label{defi:RKHS}
Let ${\cal H}(X)$ 
be a Hilbert space with inner product $\langle \cdot,\cdot \rangle$.
A function $k:X\times X \rightarrow \R$ is called
 a kernel of
 ${\cal H}(X)$  
 if there exists a a feature map $\Phi: X\rightarrow {\cal H}$ such that
 $\forall\, {\bm x},{\bm x}'\in X$,
 $k({\bm x},{\bm x}'):= \langle \Phi({\bm x}), \Phi({\bm x}')\rangle $.
 ${\cal H}(X)$  is said to be a reproducing kernel Hilbert space 
 if there is a kernel  $k:
 X \times X\to \R$ such that:
 (a) $
  k_{\bm x}:=
 k_{\bm x}(\cdot)=k(\cdot, {\bm x})\in {\cal H}(X)$ for all ${\bm x}\in X$ and (b)
${\bm f}({\bm x}) =\langle {\bm f}, {\color{black}k_{\bm x}}
\rangle$ 
for all ${\bm f}\in {\cal H}(X)$ and ${\bm x}\in X$.
In this case,
$k$ is said to be a {\em reproducing kernel} of ${\cal H}(X)$ and 
${\cal H}_k$ is written for ${\cal H}(X)$.
The norm in ${\cal H}_k$ is subsequently 
defined by $\|{\bm f}\|_k:=\sqrt{\langle {\bm f},{\bm f} \rangle}$.

\end{definition}

A kernel $k:X\times X \to \R$ is said to be symmetric if $k({\bm x}, {\bm t}) = k({\bm t}, {\bm x})$ for each ${\bm x},{\bm t}\in X$,
positive semidefinite if for any finite set $\{{\bm x}_1,\cdots, {\bm x}_m\}\subset X$,
the $m\times m$ matrix whose
$(i, j)$ entry is 
$k({\bm x}_i, {\bm x}_j)$ is positive semidefinite.
A kernel $k$ is called a {\em Mercer kernel} if it is continuous, symmetric, and positive semidefinite.
Let  $\mathbb{N}$ denote the set of positive integers.
We list some of Mercer kernels as follows,
see e.g.~\cite[Corollary~2.19]{CuZ07}  and
\cite[page 481]{Mur12}.

\begin{itemize}
\item[(i)] \textbf{Polynomial kernel}: 
$k({\bm x}_1,{\bm x}_2)=(\gamma 
{\bm x}_1^T{\bm x}_2
$
where $\gamma>0$ is a constant and
$d \in \mathbb{N}$.

\item[(ii)] \textbf{Gaussian kernel}:  
$k({\bm x}_1,{\bm x}_2)=e^{-\gamma\|{\bm x}_1-{\bm x}_2\|_2^2},
$
where $\gamma >0$ is a constant.

\item[(iii)] \textbf{Inverse multiquadrics kernel:}
$k({\bm x}_1,{\bm x}_2)=(c^2+\|{\bm x}_1-{\bm x}_2\|^2)^{-\alpha}$,
where $\alpha>0$ and $c>0$.
\end{itemize}

In this paper,
we assume that a Mercer kernel $k(\cdot, \cdot)$ is given.
By Moore–Aronszajn theorem \cite{Aro50},
there is a unique RKHS ${\cal H}_k$ for which $k$ is a reproducing kernel.
In particular,
let ${\cal H}_0:={\rm span}\{k_{\bm x}: x \in X\}$,
then ${\cal H}_0$ is dense in
${\cal H}_k$.
We can let 
$$
{\cal F}=\left\{\sum_{i=1}^m \alpha_i k_{{\bm x}_i}: \alpha_i\in \R, {\bm x}_i \in [a_i,b_i]\subset X, i=1,\cdots,m\right\}\subset {\cal H}_0.
$$
Note that kernels allow us to efficiently fit linear models to non-linear data without explicitly transforming them to feature spaces where they are linear.
An interesting case is polynomial kernel $k({\bm x},{\bm x}')= ({\bm x}^T{\bm x}')^2$
with ${\bm x},{\bm x}'\in \R^2$.
There exists a feature mapping $\phi:\R^2\to \R^3$:
$(x_1,x_2)\mapsto (x_1^2,\sqrt{2}x_1x_2, x_2^2)$ such that 
$\langle \phi({\bm x}), \phi({\bm x}')\rangle
=(x_1x_1')^2+2x_1x'_1x_2x'_2+(x_2x'_2)^2=k({\bm x},{\bm x}')$.
In such a case,
the non-linear data in $\R^2$ is almost linear in the feature space and fits kernel-based linear models. 

We make a blanket assumption that  $c({\bm z},{\bm f}({\bm x}))$ is measurable for each ${\bm f}$ and the expected value is well-defined.
A sufficient condition is that
$c({\bm z},{\bm f}({\bm x}))$ is continuous for every fixed ${\bm z}$ 
and is measurable for each fixed ${\bm f}$,
and 
 there is a function $\psi$ such that
$c({\bm z}, {\bm f}({\bm x})) \leq \psi({\bm z})$, $\forall {\bm z}\in Z$, ${\bm f} \in {\cal F}$ with $\bbe_P[\psi({\bm z})]<\infty$.
Moreover,
we make the following 
blanket 
assumption 
to secure the existence
of an optimal solution to problem (\ref{eq:ML-Rf-min}):
there exists some 
$\alpha\in \R$
 such that
the lower-level sets 
\bgeqn
\label{eq:inf_com}
{\rm lev}_{\leq \alpha} R_P
:=\left\{{\bm f} \in \F: R_P({\bm f}) \leq \alpha \right\}, 
\edeqn
are nonempty and relatively compact under the weak 
topology of ${\cal H}_k$.
It is known as inf-compactness condition, which is the weakest condition to guarantee the existence of 
an optimal solution to a minimization problem,
see \cite[Theorem~A.6.9]{StC08}.
By \cite[Theorem 8.2.11]{aubin2009set}, it ensures
that  the optimal values and 
optimal solutions of problem (\ref{eq:ML-Rf-min})
is measurable when $c$ is a Carath\'eodory function.
The assumption 
implies that the set of the optimal solutions 
to problem (\ref{eq:ML-Rf-min})
is bounded, which is
 satisfied when either 
 $\F$ is  
 weakly compact 
 and/or 
 $R_{P}({\bm f})$ is coercive, i.e.,
  $R_{P}({\bm f})\to\infty$ as $\|{\bm f}\|_k\to \infty$.

\subsection{The regularized empirical risk minimization model}

In the case when 
the true probability distribution $P$ is known,
${\cal F}$ is the set of all measurable functions,
and the cost function 
$c$ takes a specific form,
we may obtain a closed form of the optimal solution. For example, when $c$ is a quadratic loss function,
the minimizer is $\bbe_P[y|{\bm x}]$ (\cite{BiN06}).
When $c({\bm z},{\bm f}({\bm x}))=\max\{(1-\nu)({\bm f}({\bm x})-y), \nu(y-{\bm f}({\bm x}))\}$,
for $0<\nu<1$,
the minimizer is a conditional $\nu$-quantile function of $P$.
However, in many practical data-driven problems, 
the true probability distribution $P$ 
is often 
unknown and the problem
is generally ill-posed.
With 
iid samples $\{{\bm z}^i=({\bm x}^i,y^i)\}_{i=1}^N$, 
we consider the regularized empirical risk minimization  problem
\begin{equation}
\label{eq:ML-saa-r}
\min_{{\bm f}\in {{\cal F}}}R_{P_N}^{\lambda_N}({\bm f}) :=\bbe_{P_N}[c({\bm z},{\bm f}({\bm x}))]
+ \lambda_N \|{\bm f}\|_k^2,
\end{equation}
where $\lambda_N>0$ is a regularization parameter,
$
P_N(\cdot) :=\frac{1}{N} \sum_{i=1}^N \delta_{{\bm z}^i}(\cdot)
$
denotes the empirical probability measure and
 $\delta_{{\bm z}^i}(\cdot)$ denotes the Dirac measure at ${\bm z}^i$ for $i=1,\cdots,N$.
 In the case that $\lambda_N=0$,
problem  (\ref{eq:ML-saa-r}) is the sample average approximation (SAA) of problem (\ref{eq:ML-Rf-min}).
In 
general, $\lambda_N$ is driven to $0$ but the choice of 
the value may affect the rate of convergence. Some papers
have been devoted to this, see for instance 
Breheny and Huang \cite{BrH15} for logistic regression models
in a finite-dimensional space, Cucker and Smale 
\cite{CuS02} and Caponnetto and Vito \cite{CaD07} 
for regularized least squares models in 
RKHS.

The coerciveness of the objective function ensures well-definedness of  problem (\ref{eq:ML-saa-r}).
Under the condition that the cost function $c$ is convex w.r.t. the second argument,
the objective function in problem (\ref{eq:ML-saa-r}) is 
strongly convex and has bounded level sets,
and thus (\ref{eq:ML-saa-r}) has a unique solution ${\bm f}_{P_N,\lambda_N}$.
By the representation theorem (see 
e.g.~\cite[Theorem 1]{SHS01} and \cite[Theorem 4.2]{ScS02}), 
a solution of 
problem~(\ref{eq:ML-saa-r}) with ${\cal F}={\cal H}_k$ can be 
written in the following form
\bgeqn
\label{eq:kernel_estimator}
{\bm f}_{P_N,\lambda_N}=\sum_{j=1}^N \alpha_j^* k_{{\bm x}_j},
\edeqn
which means problem (\ref{eq:ML-saa-r}) can be
equivalently written as
\begin{equation}
\label{eq:ML-saa-r-equiv}
\min_{{\bm f}\in {{\cal F}_N(P_N)}}R_{P_N}^{\lambda_N}({\bm f}) :=\bbe_{P_N}[c({\bm z},{\bm f}({\bm x}))]
+ \lambda_N \|{\bm f}\|_k^2,
\end{equation}
where
\bgeq
{\cal F}_N(P_N)
=\left\{{\bm f}=\sum_{j=1}^N \alpha_j k_{{\bm x}_j}\in {\cal F}:
(\alpha_1,\cdots,\alpha_N)\in \R^N
\right\}.
\edeq

We call an optimal solution (\ref{eq:kernel_estimator}) to  problem (\ref{eq:ML-saa-r}) a
{\em kernel learning estimator}
(in this paper, it represents the regularized empirical risk minimizer when $c$ is convex 
and the stationary solution 
when $c$ is not convex)
based on the training data and 
 the corresponding optimal value 
 the regularized empirical risk. 

 \subsection{Data perturbation}
 \label{sec:data_perturbation}
 
In practice,
samples are often obtained from perceived empirical data which are potentially 
perturbed,
which means  that they are not generated by $P$, 
rather they are generated by some 
$Q$ (a perturbation of $P$).
We call ${\bm z}^1,\cdots,{\bm z}^N$ real data or 
pure data which
are not  
perturbed
but are not obtainable in our setup.
With the perceived data $\tilde{{\bm z}}^1,\cdots,\tilde{{\bm z}}^N$,
we define $Q_N:=\frac{1}{N}\sum_{i=1}^{N}\delta_{\tilde{\bm z}^i}$,
and consider the perturbed regularized empirical risk minimization problem
\begin{equation}
\min_{{\bm f}\in {{\cal F}}} \;R_{Q_N}^{\lambda_N}({\bm f}) =\bbe_{Q_N}[c({\bm z},{\bm f}({\bm x}))]
+ \lambda_N \|{\bm f}\|_k^2.
\label{eq:ML-saa-Q_N}
\end{equation}
Let ${\bm f}_{Q_N,\lambda_N}$ be an optimal solution to (\ref{eq:ML-saa-Q_N}). 
In the literature,
there are usually two 
ways to describe convergence of the optimal solution.
One is to look into 
convergence of 
$
{\bm f}_{Q_N,\lambda_N}$
as $N\to \infty$.
By the law of large
numbers, 
 $
 {\bm f}_{Q_N,\lambda_N}\to {\bm f}_{Q}
 $
 as $N\to \infty$
and
$\lambda_N\to 0$
 under some moderate conditions.
This is known as the asymptotic convergence/consistency,
see e.g. \cite{MNPR06,QCLP22,guo2023statistical}.
The other is to investigate convergence of 
$
{\bm f}_Q$ to ${\bm f}_P$ (an optimal solution of problem (\ref{eq:ML-Rf-min})) 
as $Q\to P$, 
which is known as 
stability of the optimal solution ${\bm f}_P$ when $P$ is 
perturbed to $Q$, 
see e.g. \cite{Rom03}.
In this context,
 $P_N$ is not obtainable,
our interest is to 
investigate whether ${\bm f}_{Q_N,\lambda_N}$ 
can be used as 
an estimator of 
${\bm f}_P$ 
when $N$ is sufficiently large and $\lambda_N$ is close to zero. 
This requires not only stability 
of 
${\bm f}_P$, 
but also uniform consistency of ${\bm f}_{Q_N,\lambda_N}$ to ${\bm f}_Q$ for all $Q$ near $P$.
Statistical robustness provides 
a new effective approach for this. 
It examines
the difference between empirical 
distributions of ${\bm f}_{Q_N,\lambda_N}$ 
and ${\bm f}_{P_N,\lambda_N}$. 
This is because  we are guaranteed in theory that 
the latter is a reliable estimator of ${\bm f}_P$ (${\bm f}_{P_N,\lambda_N}\to {\bm f}_{P}$), the former will also be a reasonably 
good estimator when the difference is small.

To quantify the closeness of ${\bm f}_{Q_N,\lambda_N}$ and  ${\bm f}_{P_N,\lambda_N}$,
we need to choose an appropriate metric to quantify the distance between them. 
We observe that solutions ${\bm f}_{P_N,\lambda_N}$  and ${\bm f}_{Q_N,\lambda_N}$ continuously depend on
the whole samples $\{{\bm z}_i\}_{i=1}^N$ and $\{\tilde{\bm z}_i\}_{i=1}^N$ respectively, 
so we can consider them as ${\cal H}_k$-valued random variables, defined on a product of the original probability space, i.e., $(Z^{\otimes N}, {\cal B}(Z)^{\otimes N}, P^{\otimes N})$. 
Consequently, 
we may use the difference between 
the probability distributions 
induced by the  estimators corresponding to ${\bm f}_{Q_N,\lambda_N}$ and ${\bm f}_{P_N,\lambda_N}$ to quantify their distance.
This kind of research is in alignment with
the statistical robustness (\cite{CDS10,GuX20-SR,guo2023statistical,KSZ14}),
which examines the continuity of the difference between 
the induced probability distributions 
w.r.t. the original probability distribution underpinning the random data ($Q$ and $P$). 
We will give formal definitions in Sections~\ref{sec:Quali-SR-KLE} and \ref{sec:Quant-SR-KLE}.

A key assumption in statistical robustness analysis
is that both the real data $\{{\bm z}_i\}_{i=1}^{N}$ and the 
perturbed
data $\{\tilde{\bm z}_i\}_{i=1}^{N}$ are iid.
While this assumption is widely used 
in the literature on statistical robustness
and machine learning \cite{CaD07, KSZ12,KSZ17,Mur12},
it is not 
fulfilled in the case when some 
perturbed
data are corrupted (significantly contaminated) while most of the data are good \cite{HuR09}.
For example, some data (outliers)  
deviate significantly from most of the data set.
Outliers do occur in practice. There are often no or virtually no gross errors in high-quality data, but $1$ percent to $10$ percent of gross errors in routine data seem to be 
more than the exception, 
see Hampel et al.~\cite[pages 25, 27]{Ham86}.
Consequently, the question arises
as to 
what
impact 
a small proportion of data such as outliers
may have on the kernel learning estimator.
This issue has been raised and partly investigated 
by Steinwart and Christmann in \cite{StC08} in the context of 
support vector machines.
Specifically, the authors introduce the
notions of qualitative statistical robustness (\cite[Definition 10.2]{StC08}) in all data perturbation case and the influence function in single data perturbation case and give strong motivations for investigating them albeit they have not presented further discussions on the former. Moreover,
they provide a number of interesting
examples for the influence function of empirical regularized risk minimizer ${\bm f}_{P,\lambda}$ with nice graphical illustrations.
Here we revisit the issue under a broader framework of kernel learning problems in 
hope that we will be able to fill out some gaps for both the all data perturbation case (Sections~\ref{sec:Quali-SR-KLE} and \ref{sec:Quant-SR-KLE}) and the single data perturbation case (Section~\ref{Sec:Single-data-IF}) by giving a comprehensive theoretical treatment of both.

\subsection{\texorpdfstring{$\psi$}{}-weak topology}
The statistical robust analysis requires 
some basic concepts and results about weak topology. 
Here we recall some  
relevant 
materials 
extracted mainly  from 
\cite{Clau16,FoS16} and references therein
  for a more comprehensive discussion.

We write ${\cal C}_b(Z)$ for the set of all bounded and continuous
functions on $Z$
and 
$\mathscr{P}(Z)$
for the set of all probability measures on the measurable space $\left(Z,{\cal B}(Z)\right)$.
The weak topology on $\mathscr{P}(Z)$ is the coarsest topology for which 
the mapping
$P\mapsto \int_{Z} g({\bm z}) P(d{\bm z})$ is continuous
for every
$g\in {\cal C}_b(Z)$.
Let $\psi:Z\rightarrow [0,\infty)$ be 
a gauge function,
that is, 
a continuous function 
with $\psi({\bm z})\geq 1$
outside some compact set.
We denote by ${\cal C}_{\psi}(Z)$ 
the linear space of all continuous functions $g$ on $Z$ for which there exists a constant $\kappa$
such that
$|g({\bm z})|\leq \kappa \psi({\bm z})$ for all ${\bm z}\in Z$.
Let
\bgeqn
\label{eq:M_Zpsi}
{\cal M}_Z^\psi := \left\{P'\in \mathscr{P}(Z): \int_{Z} \psi({\bm z})P'(d{\bm z})< \infty
\right\}.
\edeqn 
${\cal M}_Z^\psi $ defines a subset of probability measures in $\mathscr{P}(Z)$ with $\psi$-finite moment.
In the case when $\psi({\bm z})=\|{\bm z}\|^p$, where $p\geq 1$, we 
write
\bgeqn 
{\cal M}_Z^p := \left\{P'\in \mathscr{P}(Z): \int_{Z} \|{\bm z}\|^pP'(d{\bm z})< \infty \right\}.
\label{eq:M_Z-p}
\edeqn 
 The 
 $\psi$-weak topology on ${\cal M}_Z^{\psi}$ is the coarsest topology 
 for which 
 mapping
 $P\mapsto \int_{Z} g({\bm z}) P(d{\bm z})$
 is continuous for every 
 $g\in {\cal C}_{\psi}(Z)$.
 A sequence  $\{P_N\}_{N=1}^{\infty} 
 \subset {\cal M}_Z^{\psi}$
 is said to 
 $\psi$-weakly 
 converge to $P$, written
 ${P_N} \xrightarrow[]{\tau_\psi} P$,
 if 
 $\int_{Z} g({\bm z}) P_N(d{\bm z})\rightarrow \int_{Z} g({\bm z}) P(d{\bm z})$ 
for every $g\in {\cal C}_{\psi}(Z)$. 
We refer readers to
\cite[Definition A.47]{FoS16}
for more details of the definition.
Moreover, it follows by Corollary 2.62 in \cite{Clau16} that the $\psi$-weak topology on ${\cal M}_Z^\psi$ is generated by the metric
$\dd_\psi:{\cal M}_Z^\psi\times {\cal M}_Z^\psi\to \R$ defined by
\bgeqn
\dd_\psi(P',P''):=\dd_{\inmat{Prok}}(P',P'')+\left|\int_{Z} \psi(z) P'(d{\bm z})-\int_{Z} \psi(z) P''(d{\bm z})\right| \; \inmat{for} \; P', P'' \in {\cal M}_Z^{\psi},
\label{eq:d-psi}
\edeqn
where $\dd_{\inmat{Prok}}: \mathscr{P}(Z)\times \mathscr{P}(Z)\to \R_+$ is the Prokhorov metric defined as follows:
\bgeqn
\dd_{\inmat{Prok}}(P',P''):=\inf\{
\epsilon>0: P'(A) \leq P''(A^\epsilon)+\epsilon\,\;\mbox{for all}\,\;A \in \mathcal B(Z)
\}
\edeqn
where $A^\epsilon:= A + B_\epsilon(0)$ denotes the Minkowski sum of $A$ and the open ball centred at $0$ (w.r.t. the norm in $Z$) and $\mathcal B(Z)$ is the Borel sigma algebra of  $Z$.
We also define the Prokhorov metric in the space of probability measures over ${\cal F}$,
\bgeqn
\dd^{\cal F}_{\inmat{Prok}}(P',P''):=\inf\{
\epsilon>0: P'(A) \leq P''(A^\epsilon)+\epsilon\,\;\mbox{for all}\,\;A \in \mathcal B({\cal F})
\},
\label{eq:defi_Prok_F}
\edeqn
where $\mathcal B({\cal F})$ is Borel sigma algebra of metric space $({\cal F},
\|\cdot\|_k)$.
When $\psi$ takes some specific form,
the corresponding $\psi$-weak topology can be metricized by the following Fortet-Mourier metric.

\begin{definition}[Fortet-Mourier metric]\label{D-Fort-Mou-metric}Let 
\bgeqn
\label{eq:define_L}
\mathcal{F}_{p}(Z):=\left\{h: Z\rightarrow \R: |h({\bm z}')-h({{\bm z}''})|\leq L_{p}({\bm z}',{{\bm z}''})\|{\bm z}'-{{\bm z}''}\|, \forall {\bm z}',{{\bm z}''}\in Z\right\},
\edeqn
{\color{black}
be the set of 
locally Lipschitz continuous functions of growth  order $p$,
}
where $\|\cdot\|$ denotes some norm on $Z$,
$L_{p}({\bm z}',{{\bm z}''}):=\max\{1,\|{\bm z}'\|,\|{{\bm z}''}\|\}^{p-1}$ 
for all ${\bm z}',{{\bm z}''}\in Z$,
and $p\geq 1$ describes the growth of the local Lipschitz constants. 
The $p$-th order Fortet-Mourier metric over $\mathscr{P}(Z)$ is defined by
\bgeqn
\dd_{\rm FM}(P',P''):=\sup_{h\in \mathcal{F}_{p}(Z)}\left|\int_{Z}h({\bm z})P'(d{\bm z})-\int_{Z}h({\bm z})P''(d{\bm z})\right|, \forall P',P''\in \mathscr{P}(Z).
 \label{eq:FM-Kan}
\edeqn
\end{definition} 
In the case when $p=1$, it reduces to 
Kantorovich metric, {\color{black}
in which case we denote the distance by $\dd_K$.} 
We refer readers to 
R\"omisch \cite{Rom03},
and 
Rachev and R\"omisch \cite{RoR02} for a comprehensive overview of the topic.  
It is well-known that the Fortet-Mourier distance metricizes weak convergence on the set of probabilities possessing uniformly 
$p$-th moment 
when $Z$ is 
a finite-dimensional space,
see 
\cite{RoR02}.
Here we include 
a version in general normed space.
\begin{proposition}
\label{prop:metricize}
Assume that $Z$
is a separable Banach space equipped
with norm $\|\cdot\|$. 
Then the Fortet-Mourier metric $\dd_{\rm FM}$ metricizes the $\psi$-weak topology on ${\cal M}_Z^{\psi}$ for 
$\psi({\bm z}):=1+\max\{1,\|{\bm z}\|,\|{\bm z}'\|\}^{p-1}\|{\bm z}-{\bm z}'\|$,
where ${\bm z}'\in Z$ is a fixed point.
\end{proposition}

We defer the proof to the appendix to facilitate reading.
This result is also a step further from \cite[Corollary A.48]{FoS16}
which asserts that $\psi$-weak topology is metrizable but is short of giving a specific metric for the metrization.

\section{Qualitative statistical robustness}
\label{sec:Quali-SR-KLE}
In this section, we investigate from the statistical point of view how ${\bm f}_{Q_N,\lambda_N}$ deviates from ${\bm f}_{P_N,\lambda_N}$
 in terms of the push-forward distributions induced by them. 
To derive the qualitative 
statistical robustness result,
we need the following assumption on 
the cost function $c$ and the feasible set ${\cal F}$. 

\begin{assumption}\label{A:cost-1}
Let $c$ be defined as in (\ref{eq:ML-Rf-min}).
The following hold.
\begin{itemize}
\item[${\rm (C1)}$] For almost every ${\bm z}\in Z$, $c({\bm z},\cdot)$ is convex.
\item[${\rm (C2)}$] There is a
 gauge function $\psi(\cdot)$
such that
\bgeqn
c({\bm z},{\bm f}({\bm x}))\leq \psi({\bm z}), \;\;\forall {\bm z}\in Z \; \inmat{and} \; {\bm f}\in \F.
\label{eq:growth-ml}
\edeqn
\item[${\rm (C3)}$] The function $c: Z\times Y\to \R$ is uniformly continuous
over any compact subset of $Z\times Y$.
\end{itemize}
\end{assumption}

Assumption~\ref{A:cost-1}(C1) is a sufficient condition for the uniqueness of the optimal solution of the regularized problem (\ref{eq:ML-saa-Q_N}) with $\lambda_N>0$.
Assumption~\ref{A:cost-1}(C2)
is trivially satisfied when $Z$ is compact. Our focus in this section is on the case that $Z$ is unbounded.
  Obviously $\psi$ depends on the concrete structure of $c$.
 Consider for example 
$c({\bm z},{\bm f}({\bm x}))
=\max\{0,1-y{\bm f}({\bm x})\}$.
Then
\bgeq
c({\bm z},{\bm f}({\bm x}))
\leq |1-y{\bm f}({\bm x})|
\leq 1+|y|\cdot|\langle {\bm f}, {\color{black}k_{\bm x}}\rangle |
\leq  1+|y|\cdot \|{\bm f}\|_k \|{\color{black}k_{\bm x}}\|_k.
\edeq
 Moreover,
if there exists a positive number $\beta$ such that
$\|{\bm f}\|_k \leq \beta$,
then
we can work out an explicit form of $\psi$ for some specific kernels.
\begin{itemize}
\item If $k$ is a linear kernel, then
$\|{\color{black}k_{\bm x}}\|_k^2
=\|{\bm x}\|^2$ and
$ 
\psi({\bm z}):=1+\beta \|{\bm x}\||y|.
$

\item If $k$ is a Gaussian kernel or Laplacian kernel, then
$\|{\color{black}k_{\bm x}}\|_k^2=0$ and
$
\psi({\bm z}):=1.$

\item If $k$ is a polynomial kernel, then
$\|{\color{black}k_{\bm x}}\|_k^2=(\gamma \|{\bm x}\|^2+1)^{d}$ and
\bgeqn
\psi({\bm z}):=1+\beta (\gamma \|{\bm x}\|^2+1)^{d/2}|y|.
\label{eq:Poly-kel-phi}
\edeqn
\end{itemize}

We also need some conditions on the kernel function.
\begin{assumption}
\label{A:kernel}
 For any compact subset $Z_0\subset Z$, let $X_0$ be its orthogonal projection 
on $X$.
The following hold.
\begin{itemize}
\item[{\rm (K1)}] The set of functions $\{k_{\bm x}: {\bm x}\in X_0\}$
is uniformly continuous over
$X_0$, i.e.,
for any $\epsilon>0$, there exists a constant $\eta>0$ such that
$$
\|k_{{\bm x}'}-{\color{black}k_{{\bm x}''}}\|_k < \epsilon, \forall {\bm x}', {\bm x}'' \in X_0:  \|{\bm x}'-{\bm x}''\| < \eta,
$$
 where $\|\cdot\|$ is
 some norm on $X$.
 \end{itemize}
\end{assumption}
Condition (K1)
is considered in \cite{guo2023statistical}.
It is satisfied by the linear kernel,
Gaussian kernel and polynomial kernel, see details there. 
 Under (K1),
  $\|k_{\bm x}\|_k\leq \beta_{X_0}$ for all ${\bm x}\in X_0$,
 where $\beta_{X_0}$ is a positive constant depending on $X_0$.

Next, we introduce 
a formal definition of  statistical estimator $T_{(\cdot,\lambda)}$ parameterized by $\lambda$,
where $T_{(\cdot,\lambda)}$
maps from a subset of $\mathcal{M} \subset \mathscr{P}(Z)$ to ${\cal H}_k$.
To ease the exposition,
we
write
 $\vec{\bm z}^N$ for $({\bm z}^1,\cdots,{\bm z}^N)$ and
$\widehat{T}(\vec{\bm z}^N,\lambda_N)$
for $T_{P_N,\lambda_N}$ for fixed sample size $N$.

\begin{definition}[Qualitative Statistical robustness \cite{KSZ14}] Let
${\cal M} \subset \mathscr{P}(Z)$ be a set of probability measures
and $\dd_\psi$ be defined as in (\ref{eq:d-psi}) for some gauge function $\psi:Z\to \R$, let $\{\lambda_N\}$ be a sequence of parameters.
 A parameterized  statistical estimator 
$\widehat{T}(\cdot,\lambda_N)$
is said to be robust on
 $\cal{M}$ with respect to $\dd_\psi$ and $\dd_{\inmat{Prok}}^{\cal F}$ if for all
$P \in {\cal M}$ and $\epsilon>0$, there exist $\delta>0$ and $N_0\in \mathbb{N}$ such that
\bgeq
Q\in {\cal M}, \dd_\psi(P,Q) \leq \delta \Longrightarrow
\dd_{\inmat{Prok}}^{{\cal F}}\left(P^{\otimes N}\circ 
\widehat{T}(\cdot,\lambda_N)^{-1},
Q^{\otimes N}\circ 
\widehat{T}(\cdot,\lambda_N)^{-1}
\right) \leq \epsilon\,\,\inmat{for}\,N\geq N_0,
\edeq
where $\dd_{\inmat{Prok}}^{\cal F}$ is defined as in (\ref{eq:defi_Prok_F}).
\end{definition}
As explained in Section~\ref{sec:data_perturbation}, 
we care about 
the distance between the laws induced by ${\bm f}_{P_N,\lambda_N}$ and ${\bm f}_{Q_N,\lambda_N}$. 
It not only requires ${\bm f}_Q\to {\bm f}_P$ (stability), but also requires
${\bm f}_{P_N,\lambda_N}\to {\bm f}_P$ (consistency)
and 
\bgeqn 
{\bm f}_{Q_N,\lambda_N}\to {\bm f}_Q
\label{eq:uniform-covg-Q_N}
\edeqn 
uniformly for all $Q$ close to $P$ (uniform consistency).
Here ``uniform'' is needed 
because unlike $P$ (which is fixed albeit we do not know where it is in the space $\mathscr{P}(Z)$), 
$Q$ could be anywhere around $P$.
The uniform convergence (\ref{eq:uniform-covg-Q_N})
requires $Q_N$ to converge to $Q$
uniformly over a set containing 
$P$ and this is known as uniform 
Glivenko-Cantelli (GC)
 property. 
\begin{definition}[Uniform Glivenko-Cantelli property]\label{d-UGC}
 Let $\psi$ be a gauge function and $\dd_\psi$ be defined as in (\ref{eq:d-psi}).
Let ${\cal M}$ be a subset of $\mathcal{M}_Z^\psi$,
where $\mathcal{M}_Z^\psi$ is defined as in (\ref{eq:M_Zpsi}).
The metric space
$(\cal M,\dd_\psi)$
is said to have
Uniform Glivenko-Cantelli (UGC) property
if for every $\epsilon>0$ and $\delta>0$, there exists $N_0\in \mathbb{N}$
such that
\bgeq
P^{\otimes N}\left(\vec{\bm z}^N:
\dd_\psi(Q,Q_N) \geq \delta\right) \leq \epsilon\,\,\,
\inmat{for all}\,\,\, Q\in {\cal M}, N\geq N_0.
\edeq
\end{definition}
 This uniform 
 GC property in turn requires
 the perceived data to satisfy some topological structure. In other words,
 one can secure the GC property or 
 (\ref{eq:uniform-covg-Q_N}) only when 
 the perceived data is of reasonable structure such as light-tail, see \cite[Corollary 3.5]{KSZ12}.
 We are now ready to state the main result of this section.

\begin{theorem}[Qualitative statistical robustness of kernel learning estimator]
\label{t-ML-SR-QL-SR}
Let
\bgeq
{\cal M}_{Z,M_1}^{\psi^\gamma} := \left\{ P'\in \mathscr{P}(Z): \int_{Z} \psi({\bm z})^{\gamma}  P'(dz) \leq M_1\right\},
\edeq
where $M_1> 0$ and $\gamma>1$ are some positive constants.
Let 
$\widehat{\bm f}(\vec{\tilde{\bm z}}^N,\lambda_N):={\bm f}_{Q_N,\lambda_N}$.
Assume: 
(a) 
{\rm (C1)-(C3)} and {\rm (K1)}
hold;
(b)  there exists a positive number $\beta$ such that
$\|{\bm f}\|_k\leq \beta$ for all ${\bm f}\in \F$;
(c) 
the true probability distribution $P \in {\cal M}_{Z,M_1}^{\psi^\gamma}$.
Then the following assertions hold.
\begin{itemize}
 \item[(i)] If 
 $\lambda_N\downarrow \tau$ as $N\to \infty$ 
where $\tau$ is any small positive number, 
then for any $\epsilon>0$, 
there exist positive numbers
 $\delta>0$ and $N_0\in \mathbb{N}$ 
 such that
\bgeqn
Q \in {\cal M}_{Z,M_1}^{\psi^\gamma}, \; \dd_{\psi}(P,Q) \leq \delta
\Longrightarrow
\dd^{\cal F}_{\inmat{Prok}}\left(P^{\otimes N}\circ 
\widehat{\bm f}(\cdot,\lambda_N)^{-1},
Q^{\otimes N}\circ 
\widehat{\bm f}(\cdot,\lambda_N)^{-1}
\right) \leq \epsilon
\label{eq:stat-robust-vt_n-SP}
\edeqn
for all $N\geq N_0$,
 where $\dd_{\psi}$ is defined as in (\ref{eq:d-psi}) and $\dd^{\cal F}_{\inmat{Prok}}$ is defined as in (\ref{eq:defi_Prok_F}).

\item[(ii)] If, in addition, 
(d) problem (\ref{eq:ML-Rf-min}) has a unique solution
 ${\bm f}_{P}$,
 then 
(\ref{eq:stat-robust-vt_n-SP}) holds
for all $N\geq N_0$
 and $\lambda_N\downarrow 0$. 
\end{itemize}
\end{theorem}

Part (i) of the theorem says that
the kernel learning estimator ${\bm f}_{Q_N,\lambda_N}$
is qualitatively statistically robust for all $\lambda_N$ so long as it is lower bounded by a positive number $\tau$. In that case, 
${\bm f}_{Q_N,\lambda_N}$ converges to ${\bm f}_{Q,\tau}$ as $N\to \infty$. The
underlying reason that we want to have a positive number $\tau$ for lower bound of $\lambda_N$ is that
${\bm f}_{Q,\tau}$ is single-valued and hence continuous in $Q$ 
under the $\psi$-weak topology
whereas
${\bm f}_{Q,0}$ is not necessarily 
so.
Continuity is a key ingredient for deriving the desired qualitative statistical robustness. In Part (ii) of the theorem, we have removed 
such a  condition but require the original problem (\ref{eq:ML-Rf-min}) to have a unique solution.
Under such a circumstance, 
${\bm f}_{Q,0}$ is continuous in $Q$ and we assert that 
${\bm f}_{Q_N,\lambda_N}$
is qualitatively statistically robust for all $\lambda_N\geq 0$. 
In other words, if $\lambda_N= 0$ for all $N$, i.e., we solve
problem (\ref{eq:ML-Rf-min}) with perceived data, then the resulting kernel learning estimator is still statistically robust so long the other conditions of the theorem are fulfilled. 
Next, we make 
some remarks on the conditions of the theorem.
Condition (a) is justified earlier. 
Condition (b) is 
used in the literature, see e.g.~\cite{guo2023statistical, NoK09} and references therein. 
Condition (d) is more restrictive, it is satisfied when 
$\bbe_{P}[c({\bm z},{\bm f}({\bm x}))]$ is strictly quasi-convex.
A sufficient condition 
is that  $c({\bm z},{\bm f}({\bm x}))$
is strongly convex in ${\bm f}$ for almost every 
${\bm z}$, i.e., there exists a positive integrable function 
$\sigma({\bm z})$ such that
$c({\bm z},{\bm f}({\bm x})) - \sigma({\bm z})\|{\bm f}\|^2$ 
is convex for almost every ${\bm z}$, where $\bbe_P[\sigma({\bm z})]>0$.

\vspace{0.2cm}

\textbf{Proof of Theorem~\ref{t-ML-SR-QL-SR}.}
We only prove Part (i) as Part (ii) follows
from Part (i) directly.
Observe that 
\bgeq
\sup_{{\bm f}\in {\cal F}}\left|R_{Q'}^{\lambda'}({\bm f}) 
-R_{Q}^{\tau}({\bm f})\right|
\leq\sup_{{\bm f}\in {\cal F}} \left|\bbe_{Q'}[c({\bm z},{\bm f}({\bm x}))]-\bbe_{Q}[c({\bm z},{\bm f}({\bm x}))]\right|
+ \sup_{{\bm f}\in {\cal F}}(\lambda'-\tau) \|{\bm f}\|_k^2.
\edeq
Under the condition that $\lambda'\downarrow \tau$
and condition (c),
the second term at the 
right-hand-side (rhs)
of the inequality above goes to $0$.
On the other hand, following a similar proof to \cite[Theorem 2]{guo2023statistical}, we can show under condition (b) that the first term at the rhs of the inequality converges to $0$ 
as
$Q' \xrightarrow[]{\tau_\psi} Q$
for any $Q\in {\cal M}^{\psi^\gamma}_{Z,M_1}$.
Let $Q'= Q_N$. By \cite[Corollary 3.5]{KSZ12},
$Q_N \xrightarrow[]{\tau_\psi} Q$ as $N\to\infty$.
Likewise, 
$P_N \xrightarrow[]{\tau_\psi} P$ as $N\to\infty$.
Under condition (a),
${\bm f}_{Q_N,\lambda_N}$ is the unique solution of problem (\ref{eq:ML-saa-Q_N}). Thus,
we can use Berge's maximum theorem (see Appendix \ref{Sec:Berge-max})
to assert that ${\bm f}_{Q_N,\lambda_N}\to
{\bm f}_{Q,\tau}$
and
${\bm f}_{P_N,\lambda_N}\to
{\bm f}_{P,\tau}$ as
$N\to \infty$
and ${\bm f}_{Q,\tau}\to
{\bm f}_{P,\tau}$ as  $Q\overset{\tau_{\psi}}{\to }P$.
Next, note that
\bgeq
\dd^{\cal F}_{\rm Prok}(P^{\otimes N}\circ \widehat{\bm f}(\cdot,\lambda_N)^{-1}, Q^{\otimes N} \circ \widehat{\bm f}(\cdot,\lambda_N)^{-1})
&\leq & \dd^{\cal F}_{\rm Prok}(P^{\otimes N}\circ \widehat{\bm f}(\cdot,\lambda_N)^{-1}, \delta_{{\bm f}_{P,\tau}})\\
&& +\dd^{\cal F}_{\rm Prok}(\delta_{{\bm f}_{P,\tau}},\delta_{{\bm f}_{Q,\tau}}) \\
&& + \dd^{\cal F}_{\rm Prok}(\delta_{{\bm f}_{Q,\tau}}, Q^{\otimes N}\circ \widehat{\bm f}(\cdot,\lambda_N)^{-1}),
\edeq
where $\delta_a$ denotes the Dirac measure at $a\in {\cal F}$.
Let  $\epsilon$ be any small positive number. We consider
\begin{subequations}
\label{eq:dF_epsilon}
\begin{align}
&\dd^{\cal F}_{\rm Prok}(P^{\otimes N}\circ \widehat{\bm f}(\cdot,\lambda_N)^{-1}, \delta_{{\bm f}_{P,\tau}}) \leq \frac{\epsilon}{3}, 
\\
&\dd^{\cal F}_{\rm Prok}(\delta_{{\bm f}_{P,\tau}},\delta_{{\bm f}_{Q,\tau}}) \leq \frac{\epsilon}{3},
\\
&\dd^{\cal F}_{\rm Prok}(\delta_{{\bm f}_{Q,\tau}}, Q^{\otimes N}\circ \widehat{\bm f}(\cdot,\lambda_N)^{-1})\leq \frac{\epsilon}{3}.
\end{align}
\end{subequations}
Since $({\cal F},\|\cdot\|_k)$ is a compact metric space, 
then 
it is separable and hence
a Polish space.
By definition of Prokhorov metric $\dd^{\cal F}_{\inmat{Prok}}$ and Strassen's theorem \cite{Hub81},
inequalities (\ref{eq:dF_epsilon}) are implied respectively by
\begin{subequations}
\label{eq:df_prob}
\begin{align}
&\Prob\left(\|\widehat{\bm f}(\vec{\bm z}^N,\lambda_N)- {{\bm f}_{P,\tau}}\|_k \geq \frac{\epsilon}{3}\right)\leq  \frac{\epsilon}{3},
\label{eq:df_prob-a}
\\
& \|{{\bm f}_{P,\tau}}-{{\bm f}_{Q,\tau}}\|_k 
\leq \frac{\epsilon}{3},
\label{eq:df_prob-b}\\
&\Prob\left(\|\widehat{\bm f}(\vec{\tilde{\bm z}}^N,\lambda_N)- {{\bm f}_{Q,\tau}}\|_k \geq \frac{\epsilon}{3}\right)\leq \frac{\epsilon}{3}.\label{eq:df_prob-c}
\end{align}
\end{subequations}
By Berge's maximum theorem 
(see Theorem \ref{T-berge-max}),
the optimal solution ${\bm f}_{Q,\lambda}$ is continuous in $(Q,\lambda)$,
there exists $\delta>0$ such that when $\dd_{\psi}(P',P)<\delta$ and 
$|\lambda_N-\tau|<\delta$,
we have
$\|{\bm f}_{P',\lambda'}-{\bm f}_{P,\tau}\|_k<\frac{\epsilon}{6}$ ($\leq \frac{\epsilon}{3}$).
Thus we are left to show 
(\ref{eq:df_prob-a}) and (\ref{eq:df_prob-c}).
We only prove the latter as the former is only a special case of the latter with $Q=P$. 
Note that for fixed positive number $\tau$, 
problem
\bgeq
\min_{{\bm f}\in {{\cal F}}} R_{Q}^{\tau}({\bm f}) =\bbe_{Q}[c({\bm z},{\bm f}({\bm x}))]
+ \tau \|{\bm f}\|_k^2
\edeq
has a unique optimal solution ${\bm f}_{Q,\tau}$ for any $Q\in \mathscr{P}(Z)$. Moreover, 
for any $Q\in \mathscr{P}(Z)$,
the objective function satisfies the growth condition 
\bgeq
R_{Q}^{\tau}({\bm f}) \geq
R_{Q}^{\tau}({\bm f}_{Q,\tau}) + \tau\|{\bm f}-
{\bm f}_{Q,\tau}\|_k^2, ~~\forall {\bm f}\in {\cal F}. 
\edeq
By \cite[Lemma A.1]{guo2021existence},
\bgeq 
\|{\bm f}_{Q',\lambda}-
{\bm f}_{Q,\tau})\|_k\leq \frac{3}{\tau} \left(\sup_{{\bm f}\in {\cal F}} |R_{Q'}^{\lambda}({\bm f})-R_{Q}^{\tau}({\bm f})|\right)^{\frac{1}{2}},~~\forall Q\in \mathscr{P}(Z).
\edeq
Let $Q'=Q_N, \lambda=\lambda_N$. Then
\bgeqn 
\label{eq:solu_epsilon}
\|{\bm f}_{Q_N,\lambda_N}-
{\bm f}_{Q,\tau}\|_k &\leq& \frac{3}{\tau} \left(\sup_{{\bm f}\in {\cal F}} |R_{Q_N}^{\lambda_N}({\bm f})-R_{Q}^{\tau}({\bm f})|\right)^{\frac{1}{2}}\nonumber\\
&\leq &
\frac{3}{\tau} \left( \sup_{{\bm f}\in {\cal F}}|\bbe_{Q_N}[c({\bm z},{\bm f}({\bm x}))]-\bbe_{Q}[c({\bm z},{\bm f}({\bm x}))]|+\sup_{{\bm f}\in {\cal F}}(\lambda_N-\tau)\|{\bm f}\|_k^2 \right)^{\frac{1}{2}}\nonumber\\
&\leq &
\frac{3}{\tau} \left( \sup_{{\bm f}\in {\cal F}}|\bbe_{Q_N}[c({\bm z},{\bm f}({\bm x}))]-\bbe_{Q}[c({\bm z},{\bm f}({\bm x}))]|+(\lambda_N-\tau)\beta^2 \right)^{\frac{1}{2}}.
\edeqn
For $Q_N\overset{\tau_{\psi}}{\to }Q$,
$\bbe_{Q_N}[c({\bm z},{\bm f}({\bm x}))]-\bbe_{Q}[c({\bm z},{\bm f}({\bm x}))] \to 0$ 
for each fixed 
${\bm f}\in {\cal F}$. Moreover, 
by Assumption~\ref{A:cost-1}${\rm (C3)}$, the function
$c$ is uniformly continuous over any compact subset of $Z\times Y$,
which ensures
 $c({\bm z},{\bm f}({\bm x}))$ is uniformly 
 continuous in ${\bm f}\in {\cal F}$.
Since ${\cal F}$ is compact, by the finite covering theorem, we can show 
that 
$
\bbe_{Q_N}[c({\bm z},{\bm f}({\bm x}))]\to \bbe_{Q}[c({\bm z},{\bm f}({\bm x}))]
$
as $N\to\infty$ uniformly for all
${\bm f}\in {\cal F}$.

Let $\delta$  be such that 
when 
$\dd_{\psi}(Q_N,Q)\leq \delta$ and $|\lambda_N-\tau| \leq \delta$, the rhs of 
(\ref{eq:solu_epsilon}) is less than $\frac{\epsilon}{6}$.
On the other hand, it follows by Corollary 3.5 in \cite{KSZ12} that
$({\cal M}^{\psi^\gamma}_{Z,M_1},\dd_{\psi})$ has the UGC property which implies that there exists $N_0$ such that
\bgeqn
Q^{\otimes N}\left( \tilde{\vec{\bm z}}^N: \dd_{\psi}(Q_N,Q)\geq \frac{\delta}{2}\right) \leq \frac{\epsilon}{3}
\label{UGC-prok-P-N-x=-to-Px}
\edeqn
for all $N\geq N_0$ and $Q\in {\cal M}^{\psi^\gamma}_{Z,M_1}$ including $Q=P$. 
Let $\dd_{\psi}(Q,P)\leq \frac{\delta}{2}$ such that 
\bgeqn
\label{eq:df_prob-b2}
\|{\bm f}_{P,\tau}-{\bm f}_{Q,\tau}\|_k\leq \frac{\epsilon}{6}~~ (\leq \frac{\epsilon}{3}).
\edeqn
By (\ref{UGC-prok-P-N-x=-to-Px})
\bgeqn
\label{eq:Uniform-consist}
\frac{\epsilon}{3}
&\geq& Q^{\otimes N}
\left(\tilde{\vec{\bm z}}^N: \dd_{\psi}(Q_N,Q)
\geq \frac{\delta}{2} \right) \nonumber \\
&\geq&
Q^{\otimes N}
\left(
\tilde{\vec{\bm z}}^N: 
\dd_{\psi}(Q_N,P) \geq \frac{\delta}{2}+\dd_{\psi}(Q,P)
\right)\nonumber\\
&{\geq}&
Q^{\otimes N}
\left(\tilde{\vec{\bm z}}^N: 
\dd_{\psi}(Q_N,P)
\geq \delta
\right)~~(\inmat{because}~\dd_{\psi}(Q,P)\leq \frac{\delta}{2})\nonumber\\
&\geq&
Q^{\otimes N}
\left(\tilde{\vec{\bm z}}^N
:
\left|
{\bm f}_{Q_N,\lambda_N}-{\bm f}_{P,\tau}
\right|
\geq \frac{\epsilon}{6}
\right)\nonumber\\
&\geq &
Q^{\otimes N}
\left(
\tilde{\vec{\bm z}}^N:
\left|
{\bm f}_{Q_N,\lambda_N}-{\bm f}_{Q,\tau}
\right|
\geq|{\bm f}_{P,\tau}-{\bm f}_{Q,\tau}|
+\frac{\epsilon}{6}
\right)\nonumber\\
&\overset{(\ref{eq:df_prob-b2})}{\geq}& Q^{\otimes N}
\left(
\tilde{\vec{\bm z}}^N:
\left|
{\bm f}_{Q_N,\lambda_N}-{\bm f}_{Q,\tau}
\right|
\geq \frac{\epsilon}{3}
\right)    \nonumber\\
&=&
Q^{\otimes N}
\left(
\tilde{\vec{\bm z}}^N:
\left|
\widehat{\bm f}(\tilde{\vec{\bm z}}^N,\lambda_N)-{\bm f}_{Q,\tau}
\right|
\geq \frac{\epsilon}{3}
\right),
\forall Q\in {\cal M}^{\psi^\gamma}_{Z,M_1}.\nonumber
\edeqn
The proof is complete.
\hfill $\Box$

\begin{remark}
From the proof of Theorem~\ref{t-ML-SR-QL-SR},
we can draw the following additional conclusions which are not explicitly stated in the theorem.
\begin{itemize}

\item[(i)] \underline{Stability of ${\bm f}_{P,\lambda}$}. For fixed $\lambda$, 
${\bm f}_{Q,\lambda}\to {\bm f}_{P,\lambda}$ as  $Q\overset{\tau_{\psi}}{\to }P$ which is known as stability of ${\bm f}_{P,\lambda}$ w.r.t. perturbation of probability distribution from $P$ to $Q$.

\item[(ii)] \underline{Uniform consistency of ${\bm f}_{Q_N,\lambda}$}.
For fixed $\lambda$, 
${\bm f}_{Q_N,\lambda}$ converges to 
${\bm f}_{Q,\lambda}$ as $N\to \infty$ uniformly for all 
$Q\in {\cal M}_{Z,M_1}^{\psi^\gamma}$.
This result is known  as uniform consistency which ensures
the kernel learning estimator  converges to its true counterpart when the sample size goes to infinity
and such convergence is uniform so long as $Q$ 
lies in ${\cal M}_{Z,M_1}^{\psi^\gamma}$. The latter in turn requires a specific structure of perceived data.
\end{itemize}
   
\end{remark}
The next example shows how the qualitative statistical robustness results may be 
established in support vector machine (SVM).

\begin{example}
Consider a specific 
constrained SVM 
(see \cite[Section 1]{StC08})
\bgeq
 \inf_{{\bm f}\in {\cal F}}\; \bbe_{P}[\max\{0, 1-y {\bm f}({\bm x})\}], 
 \edeq
 where 
 ${\cal F}=\{{\bm f}\in {\cal H}_k: \|{\bm f}\|_k\leq \beta\}$,
$\beta>0$ is a constant,
${\cal H}_k$ is a RKHS 
with polynomial kernel
$k({\bm x},{\bm x}'):= (\gamma \langle{\bm x},{\bm x}'\rangle+1)^d$,
for $x,x'\in \R^n$, $d\in \mathbb{N}$,
$\gamma>0$ is a constant,
$P$ is the probability measure of random vector ${\bm z}$ with  
support $Z$.
 The regularized formulation of the problem can be written as
\bgeqn 
{\bm f}_{Q_N,\lambda_N} := \arg\min_{{\bm f}\in {\cal F}}\; \bbe_{Q_N}[\max\{0, 1-y {\bm f}({\bm x})\}]
+\lambda_N\|{\bm f}\|_k^2,
\label{eq:A2-exa}
\edeqn
where
$Q_N=\frac{1}{N}\sum_{i=1}^N \delta_{\tilde{z}^i}$ is the empirical probability measure constructed by perturbed data $\{\tilde{\bm z}_i\}_{i=1}^N$,
$\lambda_N \downarrow
\tau >0$ as $N\to \infty$.
We can set 
$\psi({\bm z})$ as in 
(\ref{eq:Poly-kel-phi}).
Then all of the conditions in  Theorem~\ref{t-ML-SR-QL-SR} are satisfied, which means 
that the 
statistical estimator 
${\bm f}_{Q_N,\lambda_N}$ defined in (\ref{eq:A2-exa})
is qualitatively statistical robust.
\end{example}

\section{Quantitative statistical robustness}
\label{sec:Quant-SR-KLE}
The qualitative 
statistical robustness results guarantee that  
the laws of ${\bm f}_{Q_N,\lambda_N}$ and
 ${\bm f}_{P_N,\lambda_N}$ 
are close when $Q$ is sufficiently close to $P$. However, it is short of giving an explicit relation between $\epsilon$ and $\delta$
in (\ref{eq:stat-robust-vt_n-SP}), which means that for a given 
$\epsilon$, it could be the case that $\delta$ should be very small which is undesirable.  In this section, we propose to address the issue by deriving quantitative statistical robustness of 
${\bm f}_{Q_N,\lambda_N}$ 
where the relation between $\epsilon$ and $\delta$ is explicitly described. To this end, 
we will confine our discussions
to a specific class of 
differentiable functions.

\subsection{First order optimality condition}
\label{Sec:3}

A key step to establish
the quantitative statistical robustness of the regularized empirical risk minimizer/stationary solution ${\bm f}_{P_N,\lambda_N}$ is to show the continuity 
of ${\bm f}_{P,\lambda}$
w.r.t a small perturbation of $P$ and $\lambda$.
In this section, 
we investigate the effect of perturbation of probability $P$
on the optimal solutions to problems
(\ref{eq:ML-Rf-min})
and (\ref{eq:ML-saa-r}).
This kind of research is well-known in the literature of stochastic programming \cite{Rom03} 
but
much in machine learning. 

We begin by deriving the first-order optimality 
conditions of the two problems.
To this end, 
we investigate the differentiability of 
functional
$R_{P}({\bm f}) = \bbe_{P}[c({\bm z},{\bm f}({\bm x}))]$. 
Our aim is to establish 
\bgeqn\label{eq:df}
D_{\bm f}(\bbe_{P}[c({\bm z},{\bm f}({\bm x}))]) 
=\bbe_{P}[D_{\bm f}c({\bm z},{\bm f}({\bm x}))],
\mbox{ for } {\bm f}\in {\cal F},
\edeqn
see Appendix \ref{Def:Gateaux-diffentibility} for the definition of $D_{\bm f}$.
We need to 
come up with some 
additional
conditions on function $c$ to ensure \eqref{eq:df}.

\begin{assumption}
\label{ass:assuption-exchange}
Let $c$ be defined as in (\ref{eq:ML-Rf-min})
and $ 
\Psi:Z\rightarrow \R_+$ is an integrable function.
The following hold.
\begin{itemize}
   \item[${\rm (C3)'}$] 
There exists $\Psi$ such that, for almost every ${\bm z}\in Z$, 
$$
|c({\bm z},w_1)-c({\bm z},w_2)|\leq 
\Psi({\bm z}) |w_1-w_2|,\quad \forall w_1,w_2\in Y.
$$
\item[{\rm (C4)}] For almost every ${\bm z}\in Z$, 
$c({\bm z},w)$ is continuously differentiable in $w$.
\end{itemize}
\end{assumption}

Let 
\bgeqn
\label{eq:set_Q}
\hat{\mathcal{P}}:=\{Q \in \mathscr{P}(Z):\bbe_Q[\psi({\bm z})] < \infty,\; \bbe_Q[\Psi({\bm z})]  < \infty\},
\edeqn
where $\psi$ is defined as in ${\rm (C2)}$.
${\rm (C3)'}$
requires 
Lipschitz continuity of $c({\bm z},w)$ in $w$.
${\rm (C3)'}$ will be used to derive the Lipschitz continuity
of the integrand 
$c({\bm z},{\bm f}({\bm x}))$  in ${\bm f}$ under (K1).
The assumption is 
standard in stability analysis and algorithm design 
in
stochastic programming,
see e.g.~\cite[Theorem~7.44]{SDR09}.
This condition is trivially satisfied when $Z$ is compact.
In \cite{NoK09}, Norkin and Keyzer
commented on page 1208 that 
``compactness of $Z$ or $Y$ is commonly accepted in the statistical learning
literature, where it allows us to apply exponential concentration measure inequalities
for bounded random variables as developed by Bernstein, McDiarmid, and Hoeffding'',
see for example Cucker and Smale 
\cite{CuS02,CuS02-2},
Bousquet and Elisseeff \cite{BoE02},  
Sch\"olkopf and Smola \cite{ScS02}, 
Poggio and Smale \cite{PoS03}, 
De Vito et al.~\cite{DCR05}.
(C4) is required when we derive the first-order optimality condition of problems (\ref{eq:ML-Rf-min})
 and  (\ref{eq:ML-saa-r}).

\begin{proposition}\label{pro:derivative-objective} 
Assume:
(a) 
{\rm (C1), (C2)}, ${\rm (C3)'}$, and {\rm(C4)}
hold;
(b) $k$ 
satisfies that 
$\bbe_{Q}[\Psi({\bm z})\|{\color{black}k_{\bm x}}\|_k]<\infty$ for all
$Q\in  \hat{\mathcal{P}}$;
(c) for each fixed ${\bm f}\in {\cal F}$,
$\bbe_Q[|c'_2
({\bm z},{\bm f}({\bm x}))|\|{\color{black}k_{\bm x}}\|_k]<\infty$, 
for all $Q\in \hat{\mathcal{P}}$.
Then
$$
D_{\bm f}\left(\bbe_{Q}[c({\bm z},{\bm f}({\bm x}))] \right)
=\bbe_Q[c'_2
({\bm z},{\bm f}({\bm x})){\color{black}k_{\bm x}}], \inmat{ for each } Q\in \hat{\mathcal{P}},
$$
where
$D_{\bm f}(R_{Q}({\bm f}))$ is the  G\^{a}teaux derivative of functional $R_{Q}$ w.r.t.~${\bm f}$ in an open neighborhood 
containing 
${\cal F}$
and 
$c'_2$ 
denotes the derivative of $c(\cdot,\cdot)$ w.r.t. the second
argument,
$c'_2
({\bm z},{\bm f}({\bm x})){\color{black}k_{\bm x}}$
is ${\cal H}_k$-valued random element.\footnote{For a basic probability space $(\Sigma, S, {\bm P})$ and let $(F, S')$ be some measurable space. A random element in $F$ (more precisely, an $F$-valued random element) is any mapping 
$
X:\Sigma\to F,
$
measurable with respect to the $\sigma$-algebras $S$ and $S'$, see  \cite[page 86]{BuK13}.
}
\end{proposition}

The
proposition states that 
under some appropriate conditions,
the 
interchange 
of 
the differentiation in ${\bm f}$ and the integration can be made
in the setting of infinite dimensional space and hence
the 
objective function
in 
(\ref{eq:ML-Rf-min}) is continuously differentiable.
The main difference between Proposition~\ref{pro:derivative-objective} and \cite[Theorem~7.44]{SDR09} is that the former considers the derivative of the nonlinear functional defined over function space ${\cal H}_k$,
where the chain rule, the definition of the derivative in infinite dimensional space, and the properties of RKHS are considered,
whereas the latter focuses on the functions defined over Euclidean space $\R^m$ (i.e., ${\bm f}$ is independent of ${\bm x}$ in our context), see also \cite[Theorem 2.7.2]{Cla83}. Note that we can 
represent $c({\bm z},{\bm f}({\bm x}))$ artificially 
as a function $\hat{c}({\bm z},{\bm f})$ and then invoke \cite[Theorem 7.44]{BoS00} or \cite[Theorem 2.7.2]{Cla83} but this may require additional conditions, we leave interested readers to explore. 
For functions defined over Banach space,
if the integrand $c({\bm z},{\bm f}({\bm x}))$ is continuous jointly in ${\bm z}$ and ${\bm f}$, and $c({\bm z},{\bm f}({\bm x}))$ is convex in ${\bm f}$ for all ${\bm z}\in Z$,
then the differentiation and the integration can be 
interchanged,
see \cite[Proposition~2.175]{BoS00},
where convexity plays an important role and the monotone convergence theorem is applied,
while in Proposition~\ref{pro:derivative-objective},
we use uniform integrability and the Lebesgue Dominated Convergence Theorem.

\noindent{\bf Proof of Proposition~\ref{pro:derivative-objective}.}
Under Assumption~\ref{ass:assuption-exchange}${\rm (C3)'}$,
$\bbe_{Q}[c({\bm z},{\bm f}({\bm x}))]$ is well-defined. Moreover,
it follows by Assumption~\ref{ass:assuption-exchange}${\rm (C4)}$
that for any ${\bm f}_1,{\bm f}_2\in {\cal F}$,
\bgeqn
\label{eq:Lipschitz-f}
|R_{Q}({\bm f}_1)-R_{\color{black}Q}({\bm f}_2)|
&\leq& \bbe_{Q}[\Psi({\bm z})|{\bm f}_1({\bm x})-{\bm f}_2({\bm x})| ]\nonumber\\
\qquad & = & \bbe_Q[\Psi({\bm z})|\langle {\bm f}_1, {k_{\bm x}}\rangle -\langle {\bm f}_2, {\color{black}k_{\bm x}}\rangle |] \nonumber \\ 
\qquad &\leq & \bbe_Q[\Psi({\bm z})\|{\color{black}k_{\bm x}}\|_k\|{\bm f}_1- {\bm f}_2\|_k]
=  L \|{\bm f}_1-{\bm f}_2\|_k,
\edeqn
where $L:=\bbe_{Q}[\Psi({\bm z})\|{\color{black}k_{\bm x}}\|_k]$.
This shows that $R_Q({\bm f})$
is Lipschitz continuous over  ${\cal F}$.
 For any 
 fixed
 ${\bm f}\in {\cal F}$,
${\bm h}\in {\cal H}_k$ with $\|{\bm h}\|_k=1$
and 
$t>0$
 consider the ratio
$$
 \mathfrak{R}_t({\bm z};{\bm h}):=t^{-1} [c({\bm z},({\bm f}+t{\bm h})({\bm x}))-c({\bm z},{\bm f}({\bm x}))].
$$
Inequality  (\ref{eq:Lipschitz-f})  implies 
$| \mathfrak{R}_t({\bm z};{\bm h})|\leq \Psi({\bm z})\|{\color{black}k_{\bm x}}\|_k \|{\bm h}\|_k$.
By Assumption~\ref{ass:assuption-exchange}${\rm (C4)}$,
$c({\bm z},\cdot)$ is continuously differentiable, 
then
$$
\lim_{t\downarrow 0}  \mathfrak{R}_t({\bm z};{\bm h})
=c'_2({\bm z},{\bm f}({\bm x})){\bm h}({\bm x})
=c'_2({\bm z},{\bm f}({\bm x}))\langle {\bm h},{\color{black}k_{\bm x}}\rangle.
$$
Consequently, by the Lebesgue Dominated Convergence Theorem
\bgeq
\lim_{t\downarrow 0}  \frac{R_{Q}({\bm f}+t{\bm h})-R_{Q}({\bm f})}{t}
&=& \lim_{t\downarrow 0} \bbe_{Q}[\mathfrak{R}_t({\bm z};{\bm h}]
=\bbe_Q[\lim_{t\downarrow 0}  \mathfrak{R}_t({\bm z};{\bm h})]\nonumber \\
&=& 
\bbe_Q[ 
c'_2({\bm z},{\bm f}({\bm x}))\langle {\bm h},{\color{black}k_{\bm x}}\rangle],
\edeq
which shows that  
$R_{\color{black}Q}({\bm f})$ is directionally differentiable. Thus
$$
R'_{\color{black}Q}({\bm f};{\bm h}) = \bbe_Q[c'_2 ({\bm z},{\bm f}({\bm x})) \langle {\bm h},
 {\color{black}k_{\bm x}}\rangle].
 $$
Since $R'_{\color{black}Q}({\bm f};{\bm h})$ is linear and continuous in ${\bm h}$, then
$R_{\color{black}Q}$ is G\^{a}teaux differentiable.
To show continuous differentiability of $R_{\color{black}Q}({\bm f})$, we note that for fixed ${\bm z}$, $c({\bm z},w)$ is continuously differentiable in $w$, and 
for fixed ${\bm x}$, with the continuity 
of $k(\cdot,\cdot)$,
${\bm f}({\bm x})=\langle {\bm f}, {\color{black}k_{\bm x}}\rangle$ is continuously differentiable in ${\bm f}$. This guarantees that for fixed ${\bm z}$, $D_{\bm f}(c({\bm z},{\bm f}({\bm x})))$ 
is continuous in ${\bm f}$.
This ensures that
$
D_{\bm f}(R_{\color{black}Q}({\bm f}))(\cdot)=\bbe_Q[c'_2 ({\bm z},{\bm f}({\bm x}))\langle \cdot,{\color{black}k_{\bm x}}\rangle ]
$
is continuous in ${\bm f}$. 
\hfill $\Box$

With Proposition~\ref{pro:derivative-objective},
we are ready to derive the first-order optimality conditions of (\ref{eq:ML-Rf-min})
 and  (\ref{eq:ML-saa-r}).

\begin{theorem}[First-order optimality conditions of (\ref{eq:ML-Rf-min})
 and  (\ref{eq:ML-saa-r})
 ]
\label{eq:first-order-optimal}
Assume the setting and conditions of  Proposition~\ref{pro:derivative-objective}.
Then the first order optimality conditions of (\ref{eq:ML-Rf-min})
 and 
 (\ref{eq:ML-saa-r})
 can be written 
 respectively as
\bgeqn
0\in 
\bbe_{P}[ c'_2({\bm z},{\bm f}({\bm x})){\color{black}k_{\bm x}}]+\N_\F({\bm f}),
\label{eq:ML-Rf-min-opti-KKT}
\edeqn
and
\bgeqn
0\in 
\bbe_{P_N}[c'_2 ({\bm z},{\bm f}({\bm x})) {\color{black}k_{\bm x}}]+ 2\lambda_N {\bm f}+ 
\N_\F({\bm f}),
\label{eq:ML-Rf-min-opti-saa-KKT}
\edeqn
where
$\N_\F({\bm f}):=\left\{{\bm d}\in {\cal H}_k: \langle {\bm d},{\bm g}-{\bm f} \rangle \leq 0,\forall {\bm g}\in {\cal F}\right\}$
denotes the normal cone in ${\cal H}_k$ for a convex set ${\cal F}\subset {\cal H}_k$ at ${\bm f}$ 
(see \cite[page 48]{BoS00}).
\end{theorem}

\noindent{\bf Proof}.
Note that $D_{\bm f}(\|{\bm f}\|_k^2)=D_{\bm f}(\langle {\bm f},{\bm f}\rangle)=2{\bm f}$.
We can apply Proposition~\ref{pro:derivative-objective} to problems (\ref{eq:ML-Rf-min})
 and (\ref{eq:ML-saa-Q_N}) 
 to 
 obtain
 (\ref{eq:ML-Rf-min-opti-KKT}) and (\ref{eq:ML-Rf-min-opti-saa-KKT}).
\hfill
\Box

Systems (\ref{eq:ML-Rf-min-opti-KKT}) and 
(\ref{eq:ML-Rf-min-opti-saa-KKT})
are 
SVIP in 
RKHS ${\cal H}_k$. 
Differing from Guo and Xu \cite{GuX23} where the decision variables are deterministic, here the decision variable is infinite dimensional and is a function of ${\bm x}$ which is random. 
In the case when
$c({\bm z},\cdot)$ is convex, any solutions to (\ref{eq:ML-Rf-min-opti-KKT}) and 
(\ref{eq:ML-Rf-min-opti-saa-KKT})
correspond to 
the
the solutions
of (\ref{eq:ML-Rf-min}) and 
(\ref{eq:ML-saa-r}).
Let ${\bm f}_{P}$ (resp.~${\bm f}_{P_N,\lambda_N}$) be
a solution to (\ref{eq:ML-Rf-min-opti-KKT}) (resp.~to  
(\ref{eq:ML-Rf-min-opti-saa-KKT})). 
In order to 
secure 
statistical robustness
of the solution
against perturbation of $P$
(resp.~$P_N$), 
we need to ensure
the perturbation is under the appropriate topology of probability measures/distributions which is more restrictive than the usual topology of weak convergence. The next subsection prepares us for this and we will come back to re-explain in Remark~\ref{rem:calM}.

\subsection{Lipschitz continuity of the stationary solution 
}
We begin by examining the continuity of ${\bm f}_{P,\lambda}$ as $P,\lambda$ varies.
We propose to do so by exploiting a generic 
stability result for an abstract generalized 
equation in Banach space in \cite{BoS00} because
SVIP (\ref{eq:ML-Rf-min-opti-saa-KKT})
may be regarded as a special
generalized equation. To this end, 
we need 
to make some 
additional 
assumptions
on 
the cost function and the kernel function and derive some intermediate technical results in the next proposition.

\begin{assumption}
\label{assu-solution-P}
Let $c$ be defined as in (\ref{eq:ML-Rf-min}) and $L_{p}({\bm z}',{\bm z}'')$ be given in (\ref{eq:define_L}).
The following hold.
\begin{itemize}
\item[${\rm (C2)'}$] 
 There exist an exponent 
{\color{black}$p_0>1$}
and a constant $C_0>0$ such that 
$$
|c({\bm z},w)| \leq C_0 (\|{\bm z}\|^{p_0-1} + |w|^{p_0-1}+1),\;
\forall ({\bm z},w)\in Z\times Y.
$$
\item[${\rm (C3)''}$] Let $\Psi$ be defined as in ${\rm (C3)'}$.
In addition to ${\rm (C3)'}$, there exists $\bar{p}_0 (= p_0-1)>1$ and $\bar{C}_0>0$ such that 
$$
\Psi({\bm z})\leq \bar{C}_0\max\{1,\|{\bm z}\|\}^{\bar{p}_0-1}.
$$
\item[${\rm (C4)'}$] 
There exist a constant $C_1>0$ and $p_1\geq 1$
such that for any $w_1,w_2 \in \R,\; {\bm z}_1,{\bm z}_2\in Z$
$$
| c'_2({\bm z}_1,w_1)- c'_2 ({\bm z}_2,w_2)|\leq 
C_1L_{p_1}({\bm z}_1,{\bm z}_2)(\|{\bm z}_1-{\bm z}_2\|+|w_1-w_2|).
$$
\item[{\rm (C5)}]
There exist a constant $C_2>0$ and $p_2\geq 1$ such that for any $w_1,w_2 \in Y,\; {\bm z}_1,{\bm z}_2\in Z$,
$$
| c''_2 ({\bm z}_1,w_1)- c''_2({\bm z}_2,w_2)|\leq 
C_2 L_{p_2}({\bm z}_1,{\bm z}_2)(\|{\bm z}_1-{\bm z}_2\|+|w_1-w_2|).
$$ 
\end{itemize}
\end{assumption}

Condition ${\rm (C2)'}$
characterizes the growth condition of 
function 
$c$ whereas 
${\rm (C4)'}$ and (C5)
stipulate locally Lipschitz 
continuity of 
$ c'_2 ({\bm z},w)$ 
and 
and global Lipschitz continuity of $c''_2 ({\bm z},w)$ w.r.t. the second argument $w$.
These assumptions are more specific and/or restrictive than those of 
${\rm (C3)'}$ and {\rm (C4)}
in terms of 
growth in ${\bm z}$ and Lipschitz continuity w.r.t.~$w$.
The condition is trivially satisfied by the quadratic loss function $c({\bm z},w)=\frac{1}{2}|y-w|^2$.
In the case when $Z$ is bounded, 
${\rm (C2)'}$ is automatically satisfied whereas ${\rm (C4)'}$ and ${\rm (C5)}$ 
reduce to global 
Lipschitz continuity in both ${\bm z}$ and $w$.

\begin{assumption}
\label{ass:assumption-Z-k}
Assume:
\begin{itemize}
\item[${\rm (K1)'}$]  Kernel function $k$ is Lipschitz continuous, i.e., there exists constant $C_3>0$ such that
$$
\|k_{{\bm x}_1}-k_{{\bm x}_2}\|_k\leq C_3
\|{\bm x}_1-{\bm x}_2\|,\quad \forall{\bm x}_1,{\bm x}_2\in X.
$$
\end{itemize}
\end{assumption}

Condition ${\rm (K1)'}$
 is stronger than condition (K1). However, it can be easily satisfied when $X$ is a compact set in that most kernel functions in machine learning are locally Lipschitz continuous. The assumption
implies that 
\bgeqn
|{\bm f}({\bm x}_1)-{\bm f}({\bm x}_2)|\leq \|{\bm f}\|_k \|k_{{\bm x}_1}-k_{{\bm x}_2}\|_k\leq 
 C_3 \|{\bm f}\|_k 
\|{\bm x}_1-{\bm x}_2\|,
\label{eq:f-LLip}
\edeqn 
for any ${\bm x}_1,{\bm x}_2\in X$.
In the case when $\|{\bm f}\|_k $ is bounded, the condition ensures local Lipschitz continuity of ${\bm f}$ over $X$.

\begin{remark}
\label{rem:well-defined}
It is possible to unify 
conditions ${\rm (C2)}'$, ${\rm (C3)}''$, ${\rm (C4)}'$ and ${\rm (K1)'}$ such that
$\bbe_Q[\psi({\bm z})] < \infty$,
$\bbe_Q[\Psi({\bm z})]  < \infty$,
$
\bbe_Q[|c'_2({\bm z},{\bm f}({\bm x}))|\|k_{\bm x}\|_k]
<\infty,
$
and
$
\bbe_Q[\Psi({\bm z})\|k_{\bm x}\|_k]<\infty
$
when $Q$ is restricted to a specific set of probability distributions. We streamline the idea as follows.
 
 \begin{itemize}
 \item[(i)] 
Under ${\rm (C2)'}$ and ${\rm (C3)''}$, we are guaranteed that
\bgeq
\hat{\mathcal{P}}=\{Q \in \mathscr{P}(Z):\bbe_Q[\psi({\bm z})] < \infty,\; \bbe_Q[\Psi({\bm z})]  < \infty\}={\cal M}_Z^{\max\{p_0-1,\bar{p}_0-1\}}.
\edeq
 \item[(ii)] Let ${\bm z}_0=({\bm x}_0,y_0)\in Z$ be fixed.
Under ${\rm (K1)'}$,
\bgeqn
\label{eq:kx_norm}
\|{k}_{\bm x}\|_k\leq \|k_{{\bm x}_0}\|_k+C_3\|{\bm x}-{\bm x}_0\|\leq (\|k_{{\bm x}_0}\|_k+2C_3\max\{1,\|{\bm x}_0\|\})\max\{1,\|{\bm z}\|\}=:\widehat{M}_1\max\{1,\|{\bm z}\|\},
\edeqn
which ensures
$\int_{X}\|k_{\bm x}\|_k Q_{X}(d{\bm x})<\infty$ for all $Q\in {\cal M}_Z^1$,
where $Q_{X}$ denotes the marginal distribution w.r.t. ${\bm x}$
and ${\cal M}_{Z}^{p}$ is defined as in (\ref{eq:M_Z-p}).
\item[(iii)]
Let ${\bm f}$ and ${\bm z}_0$ be fixed. Then 
\bgeqn
\label{eq:Lp}
L_{p_1}({\bm z},{\bm z}_0)\leq \max\{1,\|{\bm z}_0\|\}^{p_1-1}\max\{1,\|{\bm z}\|\}^{p_1-1}
\edeqn
and 
$\|{\bm z}-{\bm z}_0\|\leq \|{\bm z}\|+\|{\bm z}_0\|\leq 2\max\{1,\|{\bm z}_0\|\}\max\{1,\|{\bm z}\|\}$.
Combining with (\ref{eq:f-LLip}), we have
\bgeqn
\label{eq:z_f-difference}
\|{\bm z}-{\bm z}_0\|+|{\bm f}({\bm x})-{\bm f}({\bm x}_0)|
&\leq& (1+ C_3 \|{\bm f}\|_k) 
\|{\bm z}-{\bm z}_0\| \nonumber \\
&\leq&  2(1+C_3\|{\bm f}\|_k)\max\{1,\|{\bm z}_0\|\}\max\{1,\|{\bm z}\|\}
\edeqn
and subsequently
\bgeq
| c'_2 ({\bm z},{\bm f}({\bm x}))|
&\overset{\rm (C4)'}{\leq} & 
 |c'_2({\bm z}_0,{\bm f}({\bm x}_0))|+C_1L_{p_1}({\bm z},{\bm z}_0)(\|{\bm z}-{\bm z}_0\|+|{\bm f}({\bm x})-{\bm f}({\bm x}_0)|)  \\
&\overset{(\ref{eq:z_f-difference})}{\leq }& |c'_2({\bm z}_0,{\bm f}({\bm x}_0))|
+C_1L_{p_1}({\bm z},{\bm z}_0)\Big(2(1+C_3\|{\bm f}\|_k)\max\{1,\|{\bm z}_0\|\}\max\{1,\|{\bm z}\|\}\Big)\\
&\overset{(\ref{eq:Lp})}{\leq }& |c'_2({\bm z}_0,{\bm f}({\bm x}_0))|\\
&&
+C_1\Big(\max\{1,\|{\bm z}_0\|\}^{p_1-1}\max\{1,\|{\bm z}\|\}^{p_1-1}\Big)
\Big(2(1+C_3\|{\bm f}\|_k)\max\{1,\|{\bm z}_0\|\}\max\{1,\|{\bm z}\|\}\Big)\\
&\leq& \Big(|c'_2({\bm z}_0,{\bm f}({\bm x}_0))|
+2C_1(1+C_3\|{\bm f}\|_k)\max\{1,\|{\bm z}_0\|\}^{p_1} \Big)\max\{1,\|{\bm z}\|\}^{p_1}\\
&=:&\widehat{M}_2\max\{1,\|{\bm z}\|\}^{p_1},
\forall {\bm z}\in Z.
\edeq
Thus,  for each fixed ${\bm f}\in {\cal F}$,
$\bbe_Q[|c'_2 ({\bm z},{\bm f}({\bm x}))|]<\infty$ for all $Q\in {\cal M}_Z^{p_1}$.
\item[(iv)]
By the results in Parts (ii) and (iii) in this remark,
we obtain
\bgeq
 | c'_2 ({\bm z},{\bm f}({\bm x}))|\|k_{\bm x}\|_k 
\leq \widehat{M}_1\widehat{M}_2
\max\{1,\|{\bm z}\|\}^{p_1+1},
\forall {\bm z}\in Z.
\edeq
Thus, for each fixed ${\bm f}\in {\cal F}$,
$
\bbe_Q[|c'_2({\bm z},{\bm f}({\bm x}))|\|k_{\bm x}\|_k]
<\infty,
$
for all $Q\in  {\cal M}_{Z}^{p_1+1}$.
\item[(v)]
By Part (ii) in the remark,
under ${\rm (C3)''}$ and ${\rm (K1)'}$,
\bgeq
\Psi({\bm z})\|k_{\bm x}\|_k
\leq   \bar{C}_0 \max\{1,\|{\bm z}\|\}^{\bar{p}_0-1}\widehat{M}_1\max\{1,\|{\bm z}\|\}
=  \bar{C}_0 \widehat{M}_1 \max\{1,\|{\bm z}\|\}^{\bar{p}_0}.
\edeq
Thus $
\bbe_Q[\Psi({\bm z})\|k_{\bm x}\|_k]<\infty$,
for all $Q\in {\cal M}_Z^{\bar{p}_0}$.
\end{itemize}
Part (i) means that 
$\hat{\mathcal{P}}={\cal M}_Z^{\max\{p_0-1,\bar{p}_0-1\}}$
and
Parts (iv) and (v) mean that the conditions (b) and (c) in Proposition~\ref{pro:derivative-objective} hold for $Q\in {\cal M}_Z^{\max\{p_1+1,\bar{p}_0\}}$.
Therefore, the optimality condition in proposition~\ref{pro:derivative-objective} holds for $Q\in {\cal M}_Z^{\max\{p_0-1,\bar{p}_0-1\}}\cap {\cal M}_Z^{\max\{p_1+1,\bar{p}_0\}}={\cal M}_Z^{\max\{p_0-1,\bar{p}_0,p_1+1\}}$.
\end{remark}

Next, we  derive the second order derivative of $\bbe_Q[c({\bm z},{\bm f}({\bm x}))]$ w.r.t. ${\bm f}$,
that is, the derivative of $\bbe_Q[c'_2({\bm z},{\bm f}({\bm x}))k_{\bm x}]$, when $Q$ is restricted to a specific set of probability distributions.
\begin{proposition}
\label{prop:secon-order-d}
Let
${\rm span}\{k_{\bm x}:{\bm x}\in X\}\subset {\cal H}_k$ be the space spanned by $k_{{\bm x}}$ and 
\bgeqn
\label{eq:tilde_P}
\widetilde{\mathcal{P}}_{\bm h}:=\left\{Q\in \mathscr{P}(Z):\bbe_{Q}[|c'_2({\bm z},{\bm f}({\bm x}))|\|k_{\bm x}\|_k]<\infty, \bbe_{Q}[|c''_2 ({\bm z},{\bm f}({\bm x}))|\|T_{\bm x}{\bm h}\|_k]<\infty  \right\}
\edeqn
for ${\bm h}\in {\cal H}_k$,
where
$T_{\bm x}:{\cal H}_k\to {\rm span}\{k_{\bm x}: {\bm x}\in X\}$
is
a projection mapping on ${\rm span}\{k_{\bm x}:{\bm x}\in X\}$
 and 
 \bgeqn
T_{\bm x}{\bm h} := k_{\bm x} \langle k_{\bm x}, {\bm h}\rangle.
\label{eq:tx1}
\edeqn
Assume that ${\rm (C2)'}$,  ${\rm (C4)}'$, ${\rm (C5)}$ and
${\rm (K1)'}$ hold.
Then $\widetilde{\mathcal{P}}_{\bm h}\subset {\cal M}_Z^{\max\{p_1+1,p_2+2\}}$ 
for all ${\bm h}\in {\cal H}_k$ 
and
$$
D_{{\bm f}}(\bbe_{Q}[c'_2({\bm z},{\bm f}({\bm x}))k_{\bm x}])=\bbe_Q[c''_2({\bm z},{\bm f}({\bm x}))T_{\bm x}], \;\;\forall Q\in {\cal M}^{\max\{p_1+1,p_2+2\}}.
$$
\end{proposition}

The proof is included in Appendix~\ref{sec:proof_2d}.
Note that
\bgeqn
\langle T_{\bm x}{\bm h_1}, {\bm h_2} \rangle = \langle {\color{black}k_{\bm x}}, {\bm h_1}\rangle\langle {\color{black}k_{\bm x}}, {\bm h_2}\rangle
\label{eq:tx2}
\edeqn
for any  ${\bm h}, {\bm h}_1, {\bm h}_2\in {\cal H}_k$.
Observe that $T_{\bm x}$ is a linear operator, {\color{black}and from the definition in (\ref{eq:defi-Lnorm})}, we have
\bgeq
\|T_{\bm x}\|_{\cal L}
&=&\sup_{\|{\bm h}\|_k\leq 1} \|T_{\bm x}{\bm h}\|_k
=\sup_{\|{\bm h}\|_k\leq 1}\|{\color{black}k_{\bm x}}\langle {\color{black}k_{\bm x}}, {\bm h}\rangle\|_k
\leq \sup_{\|{\bm h}\|_k\leq 1}\|{\color{black}k_{\bm x}}\|_k |\langle {\color{black}k_{\bm x}}, {\bm h}\rangle|\\
&\leq&  \sup_{\|{\bm h}\|_k\leq 1}\|{\color{black}k_{\bm x}}\|_k \|{\color{black}k_{\bm x}}\|_k \|{\bm h}\|_k
\leq \|{\color{black}k_{\bm x}}\|_k^2
{\color{black}=k({\bm x},{\bm x})},
\edeq 
which implies that for every ${\bm x}\in X$,
$T_{\bm x}$ is a bounded linear operator,
and thus 
$T_{\bm x}\in {\cal L}({\cal H}_k)$,
which means $T_{\bm x}$ is a ${\cal L}({\cal H}_k)$-valued random element.
The next proposition states that under some moderate conditions,
both $ c'_2 ({\bm z},{\bm f}({\bm x}){\color{black}k_{\bm x}}
$ 
and 
{\color{black}
$c''_2({\bm z}, {\bm f}(x))T_{\bm x}$}
are 
locally Lipschitz continuous in ${\bm z}$
uniformly for all ${\bm f}$ in a neighborhood
${\bm f}_0$.

\begin{proposition}
\label{prop:Lip-c'-c''}
Let ${\rm (C2)'}$, ${\rm (C3)}''$, ${\rm (C4)}'$ and
${\rm (K1)'}$ hold
and ${\bm f}_0\in {\cal F}$ is fixed.
Let
${\cal V}_{{\bm f}_0}:=\{{\bm f}\in {\cal F}: \|{\bm f}-{\bm f}_0\|_k\leq \epsilon_{\cal V}\}$
be a neighborhood of ${\bm f}_0$.
Then the following assertions hold.
\begin{itemize}
\item[(i)]
There exists a positive constant $C_{{\bm f}_0}>0$ such that 
for any ${\bm d}\in {\cal H}_k$ with $\|{\bm d}\|_k\leq 1$,
$$
 |\langle (c'_2 ({\bm z}_1,{\bm f}({\bm x}_1))k_{{\bm x}_1}-
c'_2 ({\bm z}_2,{\bm f}({\bm x}_2))k_{{\bm x}_2}, {\bm d}\rangle |
\leq  C_{{\bm f}_0} L_{p_1+1}({\bm z}_1,{\bm z}_2)\|{\bm z}_1-{\bm z}_2\|,
\forall
{\bm z}_1, {\bm z}_2 \in Z,\; {\bm f}\in {\cal V}_{{\bm f}_0},
$$
where $L_p({\bm z}_1,{\bm z}_2)$
is defined as in (\ref{eq:define_L})
for all ${\bm z}_1,{\bm z}_2\in Z$, 
and $p_1\geq 1$ is defined as in ${\rm (C4)}'$.
\item[(ii)] 
In addition,
under ${\rm (C5)}$,
there exists 
a positive constant $\widehat{C}_{{\bm f}_0}>0$ such that 
for any ${\bm d}_1,{\bm d}_2\in {\cal H}_k$ with $\|{\bm d}_1\|_k\leq 1$, $\|{\bm d}_2\|_k\leq 1$,
\bgeq
\left|\left\langle \big( c''_2 ({\bm z}_1,{\bm f}({\bm x}_1))T_{{\bm x}_1}
- c''_2 ({\bm z}_2,{\bm f}({\bm x}_2))T_{{\bm x}_2}\big){\bm d}_1, {\bm d}_2\right\rangle\right|
\leq
\widehat{C}_{{\bm f}_0} L_{p_2+2}({\bm z}_1,{\bm z}_2) \|{\bm z}_1-{\bm z}_2\|,
\forall
{\bm z}_1,{\bm z}_2\in Z, {\bm f}\in {\cal V}_{{\bm f}_0}.
\edeq
where $p_2\geq 1$ is defined as in ${\rm (C5)}$.
\end{itemize}
\end{proposition}

{\color{black}
The proof is standard, we include an outline of it in Appendix~\ref{sec:proof_Lip2d}.
By Remark~\ref{rem:well-defined} and Propositions~\ref{prop:secon-order-d} and \ref{prop:Lip-c'-c''},
we know that the derivative result and the Lipschitz continuity of the derivative hold for all probability distributions in ${\cal M}_Z^{\max\{p_0-1,\bar{p}_0,p_1+1\}}\cap {\cal M}_Z^{\max\{p_1+1,p_2+2\}}={\cal M}_Z^{\max\{p_0-1,\bar{p}_0,p_1+1,p_2+2\}}$.

To ease exposition in the forthcoming discussions, 
we set
\bgeqn 
p:=\max\{p_0-1,\bar{p}_0,p_1+1,p_2+2\}.
\label{eq:defi-p}
\edeqn

With the proposition, we move on to investigate the stability of}
the following
generalized equation
\bgeqn
0\in 
\bbe_{Q}[ c'_2 ({\bm z},{\bm f}({\bm x})) k_{\bm x}]+ 2\lambda{\bm f}+ \N_{{\cal F}}({\bm f}),
\label{eq:KKT-reg}
\edeqn
which is the first-order optimality condition of  
the regularized optimization problem
\bgeqn
\label{eq:regularized-bbQ}
\min_{{\bm f}\in {\cal F}} \;\bbe_{Q}[c({\bm z},{\bm f}({\bm x}))]+\lambda\|{\bm f}\|_k^2,
\edeqn
where $Q\in \mathscr{P}(Z)$,
$\lambda \geq 0$.
We are interested in 
existence of a unique solution ${\bm f}$
to (\ref{eq:KKT-reg}) and local Lipschitz continuity of ${\bm f}$ w.r.t.~variation of $(Q,\lambda)$
near $(P,\lambda_0)$ for some $\lambda_0>0$, which is, in essence, to derive an implicit function theorem for (\ref{eq:KKT-reg}).
Observe first that under {\color{black}
the inf-compactness condition (\ref{eq:inf_com}),
}
the set of 
solutions to (\ref{eq:regularized-bbQ}) is nonempty and bounded, which means that
the set of 
solutions to (\ref{eq:KKT-reg})  is nonempty and bounded.
Let ${\bm f}_{Q,\lambda}$ be a 
stationary point
to
(\ref{eq:regularized-bbQ}).
Our ultimate interest is global Lipschitz 
continuity of ${\bm f}_{Q,\lambda}$ in $(Q,\lambda)$
which is 
{\color{black}is an important step 
towards 
establishing} quantitative statistical robustness of 
the solution in the forthcoming discussions.
To this end, we consider the case
that for 
$Q=P$ and $\lambda=\lambda_0$,
(\ref{eq:KKT-reg}) has a
solution ${\bm f}_{P,\lambda_0}$,  
and demonstrate that 
${\bm f}_{(\cdot,\cdot)}$ is locally Lipschitz continuous 
in a neighborhood of $(P,\lambda_0)$.
We then take a step further to derive
sufficient conditions under which 
${\bm f}_{(\cdot,\cdot)}$ is globally Lipschitz continuous.

{\color{black}The proof will be based on an existing 
stability result about}
an abstract generalized equation established by Bonnans and Shapiro \cite{BoS00}.
Let $S,W$ be Banach spaces,
$\phi:S\rightarrow W$ be a continuously differentiable mapping, and ${\cal N}:S\rightrightarrows W$ be a set-valued mapping.
Consider the following abstract generalized equation:
find $s\in S$ such that 
\bgeqn
\label{eq:abstract-GE}
0\in \phi(s)+{\cal N}(s).
\edeqn
$s_0$ is called 
a {\em strong regular solution} of the abstract generalized equation
(\ref{eq:abstract-GE})
if there exist neighborhoods ${\cal V}_S$ and ${\cal V}_W$ 
of $s_0\in S$ and $0\in W$
respectively
such that for every $\delta\in {\cal V}_W$,
the linearized abstract generalized equation
$
\delta\in \phi(s_0)+D\phi(s_0)(s-s_0)+{\cal N}(s),
$
which is parameterized by $\delta$, 
has a unique solution in ${\cal V}_S$, 
denoted by $\zeta(\delta)$,
and the mapping $\zeta:{\cal V}_W\rightarrow {\cal V}_S$ is Lipschitz continuous with constant $\beta$,
that is 
\bgeqn
\label{eq:linear_GE}
\|\zeta(\delta)-\zeta(\tilde{\delta})\|_Z\leq \beta \|\delta-\tilde{\delta}\|_W,\;\forall \; \delta,\tilde{\delta} \in {\cal V}_W,
\edeqn
where $\|\cdot\|_S$ and $\|\cdot\|_W$ are the respective norms in Banach spaces $S$ and $W$. 
A combination of the existence of $\zeta(\cdot)$ and the Lipschitz conditions (\ref{eq:linear_GE}) is
known as strong regularity condition in the generalized equations, see~\cite[Definition 5.12]{BoS00}.

\begin{lemma}(\cite[page 415]{BoS00})
\label{lem:Lipschitz-perturbed-phiz}
Let $V$ be an open neighborhood of $s_0$
and consider the Banach space $C^1(V,W)$ of continuously differentiable mappings $\phi:V\rightarrow W$
equipped with 
norm 
{\color{black}
$
\|\phi\|_{1,V}
:=
\sup_{s\in V} \|\phi(s)\|+\sup_{s\in V}\|D\phi(s)\|.
$
}
If $s_0$ is a strongly regular solution of the generalized equation~(\ref{eq:abstract-GE}),
then for all $\tilde{\phi}$ in a neighborhood of $\phi$ with respect to the norm $\|\cdot\|_{1,V}$,
the generalized equation 
$
0\in \tilde{\phi}(s)+{\cal N}(s)
$
has a Lipschitz continuous (and hence unique) solution 
$\bar{s}(\tilde{\phi})$
in a neighborhood of $s_0$.
\end{lemma}

{\color{black}
We are now ready to present 
our stability results about 
the generalized equation (\ref{eq:KKT-reg}).
}
Let ${\cal V}_{{\bm f}_0}\subset {\cal F}$ be an open neighborhood of some function ${\bm f}_0\in {\cal F}$,
where ${\cal F}$ is the feasible set of problem (\ref{eq:ML-Rf-min}).
Let $C^1({\cal V}_{{\bm f}_0},{\cal H}_k)$
be a space of continuously differentiable mappings $\psi:{\cal V}_{{\bm f}_0} \to {\cal H}_k$
 equipped with the norm:
 \bgeqn
 \label{eq:definition-norm-1V}
 \|\psi\|_{1,{\cal V}_{{\bm f}_0}}:=\sup_{{\bm f}\in {\cal V}_{{\bm f}_0}} \|\psi({\bm f})\|_k
 +\sup_{{\bm f}\in {\cal V}_{{\bm f}_0}} \|D\psi({\bm f})\|_{\cal L},
 \edeqn
 where  $\|\cdot\|_{\cal L}$ is defined as in (\ref{eq:defi-Lnorm}),
 ${\cal V}_{{\bm f}_{0}}$ is given in Proposition~\ref{prop:Lip-c'-c''}.
Note that $\left(C^1({\cal V}_{{\bm f}_0},{\cal H}_k),\|\cdot\|_{1,{\cal V}_{{\bm f}_0}}\right)$ is a Banach space, see \cite[page 415]{BoS00}.
 The norm defined as in (\ref{eq:definition-norm-1V}) is standard in Euclidean space $\R^m$. Specifically,
 let $K\subset \R^m$ be a domain (an open connected subset),
and $C^1(K,\R)$ denotes the set of continuously differentiable real-valued functions defined over $K$,
endowed with bounded norm $\|{\bm f}\|_{C^1(K,\R)}:=\sup_{{\bm x}\in K}|{\bm f}({\bm x})|+\sup_{{\bm x}\in K}\|D{\bm f}({\bm x})\|$,
where $D{\bm f}({\bm x})$ is the derivative of ${\bm f}$ at ${\bm x}\in K$.
$\left(C^1(K,\R),\|\cdot\|_{C^1(K,\R)}\right)$ is a Banach space,
see e.g.~\cite[page 18]{CuZ07}.

\begin{theorem}
[Lipschitz continuity of 
the solution of (\ref{eq:ML-Rf-min-opti-KKT}) and (\ref{eq:KKT-reg})]
\label{thm:lip_solution-lambda}
Let 
$p$ be defined as in (\ref{eq:defi-p}) and
${\cal M}^p_Z$ be defined as in (\ref{eq:M_Z-p}), let
${\bm f}_{Q,\lambda}$ be a solution to 
(\ref{eq:KKT-reg})
and 
${\cal V}_{{\bm f}_{P,\lambda_0}}=\{{\bm f}\in {\cal F}: \|{\bm f}-{\bm f}_{P,\lambda_0}\|_k\leq \epsilon_{\cal V}\}$ be a neighborhood
of ${\bm f}_{P,\lambda_0}$ 
under the norm 
$\|\cdot\|_k$.
Let $\psi:{\cal H}_k\to {\cal H}_k$ be defined as
$$
\psi({\bm f}):=\bbe_P[ c'_2 ({\bm z},{\bm f}({\bm x})){\color{black}k_{\bm x}}]+2\lambda_0{\bm f}.
$$
Under
${\rm (C2)'}$, ${\rm (C3)''}$, ${\rm (C4)'}$, and ${\rm (C5)}$,

the following assertions hold.
\begin{itemize}

\item[(i)] 
Let
 $Q=P\in {\cal M}^{p}_Z$ and 
$\lambda=\lambda_0\in [\tau,\bar{\lambda}]$
for some $\tau>0$ and $\bar{\lambda}>0$.
Let ${\cal V}_{\psi}$ be a neighborhood of $\psi$.
If ${\bm f}_{P,\lambda_0}$ is
a strongly regular solution
of (\ref{eq:KKT-reg}), then  (\ref{eq:KKT-reg}) has a unique 
solution $\bar{\bm f}_{\tilde{\psi}}$ for 
$\tilde{\psi}\in {\cal V}_{\psi}$
 such that
\bgeqn
\label{eq:Lip-Solution-1}
\|\bar{\bm f}_{\tilde{\psi}_1}-\bar{\bm f}_{\tilde{\psi}_2}\|_k
\leq 
\kappa_{P,\lambda_0}
\|\tilde{\psi}_1-\tilde{\psi}_2\|_{1,{\cal V}_{{\bm f}_{P,\lambda_0}}}, \forall \tilde{\psi}_1, \;\tilde{\psi}_2
\in {\cal V}_{\psi}.
\edeqn
\item[(ii)] 
Then there exist  constants $\delta>0$  and
$\bar{C}_{{\bm f}_{P,\lambda_0}}>0$ such that
\bgeqn
\label{eq:Lip-Solution}
\|{\bm f}_{Q_1,\lambda_1}-{\bm f}_{Q_2, \lambda_2}\|_k
\leq 
2 \bar{C}_{{P,\lambda_0}}\kappa_{P,\lambda_0} (\dd_{\rm FM}(Q_1,Q_2) + |\lambda_1-\lambda_2|)
\edeqn
for all $Q_1, Q_2 \in {\cal M}_Z^p$ 
satisfying $\dd_{\rm FM}(Q_1,P) \leq \delta$ and
$\dd_{\rm FM}(Q_2,P)\leq \delta$, and $\lambda_1, \lambda_2\in[\tau, \bar{\lambda}]$ satisfying $|\lambda_1-\lambda_0|\leq \delta$ and $|\lambda_2-\lambda_0|\leq \delta$.

\item[(iii)] For some constants $\gamma>1$,
$p\geq 1$ and $M_2$, 
let 
\bgeqn
 \label{eq:M-phi-bounded}
{\cal M}_{Z,M_2}^{p\gamma} := \left\{Q\in \mathscr{P}(Z): \int_{Z} \|{\bm z}\|^{p\gamma}Q(d{\bm z})\leq M_2 \right\}.
\edeqn
Let $\tau_{\|\cdot\|^p}\times \tau_{\R}$ denote the 
product topology of $\|\cdot\|^p$-weak topology and the standard topology on $\R$ (\cite[page 81]{Mun00}).
If, in addition,
(a) there is a continuous 
 ${\bm f}_{Q,\lambda}$ to
(\ref{eq:regularized-bbQ}) such that
${\bm f}_{(\cdot,\cdot)}: {\cal M}_{Z,M_2}^{p\gamma}\times [\tau,\bar{\lambda}]  \to {\cal H}_k$ 
is continuous under 
$\tau_{\|\cdot\|^p}\times \tau_{\R}$,
(b) the strong regularity condition for
(\ref{eq:KKT-reg})
holds at ${\bm f}_{P', \lambda'}$ for 
any $P'\in {\cal M}_{Z,M_2}^{p\gamma} $ and $\lambda'\in[\tau, \bar{\lambda}]$, 
then 
there exists a constant 
$C_\kappa>0$ such that for any $(Q', \lambda'), (Q'', \lambda'') \in {\cal M}_{Z,M_2}^{p\gamma}\times [\tau, \bar{\lambda}]$,
\bgeq
\|{\bm f}_{Q', \lambda'}-{\bm f}_{Q'', \lambda''}\|_k
\leq 
2 C_\kappa (\dd_{\rm FM}(Q',Q'')+|\lambda'-\lambda''|).
\edeq
\item[(iv)]
Let $\mathscr{P}_N:=\{\frac{1}{N}\sum_{i=1}^N \delta_{\tilde{\bm z}^i}: \tilde{\bm z}^i\in Z\}$.
If, in addition, $\mathscr{P}_N\subset {\cal M}_{Z,M_2}^{p\gamma}$, then 
\bgeqn
\label{eq:Lip-Solution-PN}
~~\|{\bm f}_{Q_N^1,\lambda_N}-{\bm f}_{Q_N^2,\lambda_N}\|_k
\leq \frac{2C_\kappa}{N}
 \sum_{l=1}^N \max\{1,\|{\bm z}_1^l\|,\|{\bm z}_2^l\|\}^{p-1} \|{\bm z}_1^l-{\bm z}_2^l\|,
\edeqn
where $Q_N^1=\frac{1}{N}\sum_{l=1}^N\delta_{{\bm z}_1^l}$,
$Q_N^2=\frac{1}{N}\sum_{l=1}^N\delta_{{\bm z}_2^l}$,
and $\lambda_N\in [\tau,\bar{\lambda}]$ for all $N\in \mathbb{N}$.
 \end{itemize}
\end{theorem}


{\color{black} Before presenting a proof}, we make a few comments about the conditions and results of this theorem.
Part~(i) of the theorem states {\color{black}
that}
the local Lipschitz continuity of the solution mapping of the system (\ref{eq:KKT-reg})
in terms of variations of
$
\int_Z  c'_2 ({\bm z},{\bm f}({\bm x})){\color{black}k_{\bm x}} 
Q(d{\bm z})
+2\lambda {\bm f}$
as $Q$ 
and $\lambda$ vary. 
Part~(ii) quantifies the continuity 
in terms of $Q$ and $\lambda$
under the
{\color{black}product topology
$\tau_{\|\cdot\|^p}\times \tau_{\R}$,
where topology
$\tau_{\|\cdot\|^p}$ can be metricized by the Fortet-Mourier metric (see Proposition~\ref{prop:metricize})
and $\tau_{\R}$ can be metricized by
the standard metric in $\R$ (see \cite[Example 2, page 120]{Mun00}).}
Part (iii) of the theorem says that if there is a continuous solution trajectory to the system 
(\ref{eq:KKT-reg}) over 
${\cal M}_{Z,M_2}^{p\gamma}\times [\tau,\bar{\lambda}]$, and the strong regularity condition holds at every point of the trajectory,
then the solution mapping is globally Lipschitz continuous in $Q$ and $\lambda$ over 
${\cal M}_{Z,M_2}^{p\gamma}\times [\tau,\bar{\lambda}]$
 under the product topology of the $\|\cdot\|^p$-weak topology and the standard topology in $\R$.
The continuity holds when 
(\ref{eq:KKT-reg})
has a unique solution
for every $P\in {\cal M}_{Z,M_2}^{p\gamma}$, 
$\lambda_0\in [\tau,\bar{\lambda}]$,
and strong regularity condition holds (see Lemma~\ref{lem:Lipschitz-perturbed-phiz}) although 
our interest is not restricted to this case.
Note that ${\cal M}_{Z,M_2}^{p\gamma}\subset {\cal M}_{Z}^p$.
We need the boundedness of
the 
$p\gamma$-moment because
it ensures 
relative compactness
of ${\cal M}_{Z,M_2}^{p\gamma}$ under topology $\tau_{\|\cdot\|^{p}}$ (see \cite[Lemma~2.69]{Clau16})
required in the proof of 
Part~(iii).
{\color{black}
Part (iv) 
is a specific version of Part (iii) when $Q'$ and $Q''$
are empirical probability distributions 
and $\mathscr{P}_N\subset {\cal M}_{Z,M_2}^{p\gamma}$.
The result 
prepares us for the 
statistical robustness of the regularized 
stationary point in Theorem~\ref{t:QSR_SGE}.
Finally, we note that  
all of the results cover 
(\ref{eq:ML-Rf-min-opti-KKT}) as a special case with $\lambda$ being fixed as a constant $0$.

A key condition required in this theorem is 
strong regularity of the solution ${\bm f}_{P,\lambda_0}$ in Part~(i) and 
 ${\bm f}_{Q,\lambda}$ in Part~(iii) for all $Q\in {\cal M}_{Z,M_2}^{p\gamma}$ and $\lambda\in [\tau,\bar{\lambda}]$.
To see how these conditions may be possibly satisfied, we consider the case where the 
cost function $c({\bm z},w)$ is convex in $w$.
In this case, $R_P({\bm f})=\bbe_P[c({\bm z},{\bm f}({\bm x}))]$ is convex in ${\bm f}$ and subsequently 
$
R_P({\bm f})+\lambda \|{\bm f}\|_k^2$ is strongly convex for 
$\lambda>0$. The strong convexity 
of $R_P({\bm f})+\lambda \|{\bm f}\|_k^2$
ensures that the regularized problem
(\ref{eq:regularized-bbQ}) has a unique optimal solution and the second order growth condition 
holds at the solution.
The discussions above show that the strong regularity conditions may be satisfied when
$c({\bm z},w)$ is convex in $w$ and $\lambda>0$.
 It is important to note that 
strong regularity condition does not necessarily imply 
convexity of $c({\bm z},\cdot)$.

Another important condition required in this theorem is  continuous differentiability 
of the cost function in $w$.
In the literature of machine learning, 
some cost functions are not 
continuously differentiable. 
A potential way to tackle 
this is smoothing (\cite{ralph2005implicit}) so that the smoothed cost function is twice continuously differentiable in $w$.
We can then perform the analysis 
with the smoothed problem
and drive the smoothing parameter to zero, again we leave this for future exploration.

\vspace{0.3cm}
\noindent{\it Proof of Theorem~\ref{thm:lip_solution-lambda}.}
{\color{black}
We use Lemma~\ref{lem:Lipschitz-perturbed-phiz} to prove the results. 
We begin by identifying 
$s_0$, $S$, $V$, $W$, $\phi$ 
and $\|\cdot\|_{1,V}$
of the lemma
in the context of generalized equation (\ref{eq:KKT-reg}), i.e.,
$s_0={\bm f}_0$,
$S={\cal H}_k$,
$V={\cal V}_{{\bm f}_0}$
$W=\R$,
$\phi(s)=\psi({\bm f})$,
and $\|\cdot\|_{1,V}$ 
corresponds to $\|\cdot\|_{1,{\cal V}_{{\bm f}_0}}$.
}
Part~(i) follows directly from Lemma~\ref{lem:Lipschitz-perturbed-phiz}.

Part (ii). {\color{black}In Part (i), the implicit function is defined over the space of parameter $\psi$. Since $\psi$ is determined by $(Q,\lambda)$,
then we can describe the implicit function in terms of the latter.
}
Let $\tilde{\psi}_1,\tilde{\psi}_2\in {\cal V}_{\psi}$ and 
 ${\bm f}_{Q_i, \lambda_i}=\bar{\bm f}_{\tilde{\psi}_i}$ for $i=1,2$. 
By (\ref{eq:Lip-Solution-1}),
\bgeqn
\label{eq:Lip-Solution-f}
\|{\bm f}_{Q_1, \lambda_1}-{\bm f}_{Q_2, \lambda_2}\|_k
=\|\bar{{\bm f}}_{\tilde{\psi}_1}-\bar{{\bm f}}_{\tilde{\psi}_2}\|_k
\leq \kappa_{{P,\lambda_0}} \|\tilde{\psi}_1-\tilde{\psi}_2\|_{1,{\cal  V}_{{\bm f}_{P,\lambda_0}}}.
\edeqn
To prove (\ref{eq:Lip-Solution}), it suffices to show that
$$
\|\tilde{\psi}_1-\tilde{\psi}_2\|_{1,{\cal  V}_{{\bm f}_{P,\lambda_0}}}
\leq 2\bar{C}_{P,\lambda_0} \kappa_{P,\lambda_0} (\dd_{\rm FM}(Q_1,Q_2) + |\lambda_1-\lambda_2|).
$$
Since ${\cal H}_k$ is a Hilbert space,
 the dual space ${\cal H}^*_k$ is canonically identified with ${\cal H}_k$.
 Under Assumption~\ref{assu-solution-P}${\rm (C3)''}$, ${\rm (C4)'}$ and ${\rm (C5)}$,
we have by  Proposition~\ref{prop:Lip-c'-c''}
that for any ${\bm d}, {\bm d}_i\in {\cal H}_k$ with $\|{\bm d}\|_k\leq 1$, $\|{\bm d}_i\|_k\leq 1$, $i=1,2$,
$$
\frac{1}{\max\{C_{{\bm f}_{P,\lambda_0}},\hat{C}_{{\bm f}_{P,\lambda_0}}\}}\langle {\bm d}(\cdot),  c'_2({\bm z},{\bm f}({\bm x})){\color{black}k_{\bm x}}\rangle \in {\cal F}_p(Z)$$
and 
$\frac{1}{\max\{C_{{\bm f}_{P,\lambda_0}},\hat{C}_{{\bm f}_{P,\lambda_0}}\}}\langle c''_2({\bm z},{\bm f}({\bm x})) T_{{\bm x}}{\bm d}_1,{\bm d}_2\rangle
\in {\cal F}_p(Z).
$
Note that 
$$
\|(\tilde{\psi}_1-\tilde{\psi}_2)({\bm f})\|_k=\sup_{{\bm d}\in {\cal H}_k,\|{\bm d}\|_k\leq 1}\langle {\bm d},(\tilde{\psi}_1-\tilde{\psi}_2)({\bm f})\rangle.
$$
It follows by the definition of $\|\cdot\|_{1,{\cal V}_{{\bm f}_{P,\lambda_0}}}$,
\bgeqn
 \label{eq:norm_W}
&&\|\tilde{\psi}_1-\tilde{\psi}_2\|_{1,{\cal  V}_{{\bm f}_{P,\lambda_0}}} \nonumber \\
&=& \sup_{{\bm f}\in {\cal V}_{{\bm f}_{P,\lambda_0}}} \|(\tilde{\psi}_1-\tilde{\psi}_2)({\bm f})\|_k+ \sup_{{\bm f}\in {\cal V}_{{\bm f}_{Q,\lambda_0}}}\|D_{\bm f}((\tilde{\psi}_1-\tilde{\psi}_2)({\bm f}))\|_{\cal L}\nonumber \\
&\leq &\sup_{\substack{{\bm f}\in {\cal V}_{{\bm f}_{P,\lambda_0}} \nonumber \\
{\bm d}\in {\cal H}_k, \|{\bm d}\|_k\leq 1}}
\langle {\bm d},\bbe_{Q_1}[c'_2({\bm z},{\bm f}({\bm x})){\color{black}k_{\bm x}}]- \bbe_{Q_2}[ c'_2 ({\bm z},{\bm f}({\bm x})){\color{black}k_{\bm x}}]\rangle+2|\lambda_1-\lambda_2|\|{\bm f}\|_k\nonumber \\
&&+\sup_{\substack{{\bm f}\in {\cal V}_{{\bm f}_{P,\lambda_0}}, {\bm d}_1,{\bm d}_2\in   {\cal H}_k\\
\|{\bm d}_1\|_k\leq 1,\|{\bm d}_2\|_k\leq 1}} 
|\langle \big(\bbe_{Q_1}[c''_2({\bm z},{\bm f}({\bm x}))
T_{{\bm x}}
]-\bbe_{Q_2}[ c''_2 ({\bm z},{\bm f}({\bm x})){\color{black}T_{{\bm x}}
}] \big){\bm d}_1,{\bm d}_2\rangle| +2|\lambda_1-\lambda_2| \nonumber \\
&\leq& 2\bar{C}_{{P,\lambda_0}} (\dd_{\rm FM}(Q_1,Q_2)+|\lambda_1-\lambda_2|) 
\edeqn
for $Q_1,Q_2\in {\cal M}_Z^p$,
where $\bar{C}_{{P,\lambda_0}}=\max\{C_{{\bm f}_{P,\lambda_0}},\hat{C}_{{\bm f}_{P,\lambda_0}}, 
\|{\bm f}_{P,\lambda_0}\|+\epsilon_{\cal V}\}$,
with $\epsilon_{\cal V}$ being given in Proposition~\ref{prop:Lip-c'-c''}.
{\color{black}
Let $\tilde{\psi}_1=\tilde{\psi}_i$ and $\tilde{\psi}_2={\psi}$,
then there exists $\delta>0$ such that 
$\tilde{\psi}_i\in {\cal V}_{\psi}$ for all $(Q_i,\lambda_i)\in {\cal M}_Z^p\times \R_+$ satisfying $\dd_{\rm FM}(Q_i,P)\leq \delta$ and $|\lambda_i-\lambda_0|\leq \delta$, $i=1,2$.
Combining (\ref{eq:norm_W}) and (\ref{eq:Lip-Solution-f}),
(\ref{eq:Lip-Solution}) holds.
}

Part (iii).
We use the finite covering theorem to prove the result.
Observe first that $[\tau, \bar{\lambda}]$ is compact 
under the standard topology,
\cite[Example 1, page 14]{Mun00} and 
\cite[Theorem~27.1]{Mun00},
and ${\cal M}_{Z,M_2}^{p\gamma}$ is 
relatively compact
under $\|\cdot\|^{p}$-weak topology
when $p\geq 1$ and $\gamma>1$ (see \cite[Lemma~2.69]{Clau16})).
Moreover, 
it follows from  
and Proposition~\ref{prop:metricize}
that
$\dd_{\rm FM}$ metricizes the $\|\cdot\|^p$-weak topology.
 Thus, by virtue of Prokhorov's theorem \cite{Prok56}
 and the compactness of $[\tau,\bar{\lambda}]$, 
we can construct a $\delta$-net ${\cal Q}^\lambda_J :=\{(Q_1,\lambda_1), \cdots, (Q_J, \lambda_J)\}$  
in 
${\rm cl}({\cal M}_{Z,M_2}^{p\gamma})\times [\tau, \bar{\lambda}]$ under the metric $\dd_{\rm FM}(\cdot, \cdot) \times d_{\R}(\cdot,\cdot)$ 
such that
 ${\cal M}_{Z,M_2}^{p\gamma}\times [\tau, \bar{\lambda}]  \subset \cup_{j=1}^{J} B((Q_j, \lambda_j),\delta)$, 
 where 
$d_{\R}(\lambda_1,\lambda_2):=|\lambda_1-\lambda_2|$,
 ${\rm cl}(S)$ denotes the topological closure of the set $S\subset \mathscr{P}(Z)$,
 $B(Q_j,\delta)$ denotes a closed ball 
 centered at $Q_j$ with radius $\delta$ under the metric $\dd_{\rm FM}(\cdot, \cdot) 
 \times d_{\R}(\cdot,\cdot)$.
 Note that the argument can also be looked at from a different perspective. By \cite[Theorem~2.15]{HuR09},
$\mathscr{P}(Z)$ is a Polish space. Thus we can apply the finite covering theorem to
${\rm cl}({\cal M}_{Z,M_2}^{p\gamma})\times [\tau, \bar{\lambda}]$ which is weakly compact
under topology $\tau_{\|\cdot\|^p}\times \tau_{\R}$.
Since ${\bm f}_{(\cdot, \cdot)}$  is assumed to be defined continuously over
${\cal M}_{Z,M_2}^{p\gamma}\times [\tau, \bar{\lambda}]$ and the strong regularity condition holds 
at every point ${\bm f}_{Q,\lambda}$ for $(Q,\lambda)\in {\cal M}_{Z,M_2}^{p\gamma}\times [{\color{black}\tau}, \bar{\lambda}]$, then
we may set $\delta$ sufficiently small such that ${\bm f}_{Q,\lambda}$ is the unique 
solution to (\ref{eq:regularized-bbQ}) for $(Q,\lambda)\in B(Q_j,\delta)$, $j=1,\cdots,J$. Following a similar argument to Part (ii), we can show that
 \bgeqn
 \label{eq:f_Lipschitz_Qj}
 \|{\bm f}_{Q'_j, \lambda'_j}-{\bm f}_{Q''_j, \lambda''_j}\|_k\leq 2 C_{{Q_j, \lambda_j}} (\dd_{\rm FM}(Q'_j,Q''_j) + |\lambda'_j-\lambda''_j|)
 \edeqn
for any $(Q'_j, \lambda'_j), (Q''_j, \lambda''_j)\in 
 B((Q_j,\lambda_j), \delta)
 $ satisfying 
 $\dd_{\rm FM}(Q'_j,Q_j)+|\lambda'_j-\lambda_j|\leq \delta$ and 
 $\dd_{\rm FM}(Q''_j,Q_j)+|\lambda''_j-\lambda_j|\leq \delta$
 for $j=1,\cdots,J$, where
$C_{Q_j,\lambda_j}:={\bar C}_{Q_j,\lambda_j}\kappa_{Q_j,\lambda_j}$
is a positive constant
depending on $(Q_j,\lambda_j)$.

For any $(Q', \lambda'), (Q'', \lambda'') \in {\cal M}_{Z,M_2}^{p\gamma}\times[{\color{black}\tau}, \bar{\lambda}]$
and $\theta\in [0,1]$, define
$(Q(\theta), \lambda(\theta)):= ((1-\theta)Q', (1-\theta)\lambda')+ (\theta Q'', \theta\lambda'')$.
Let $(\hat{Q}_1, \hat{\lambda}_1)\in {\cal Q}^\lambda_J$ 
be such that $(Q', \lambda')\in B((\hat{Q}_1, \hat{\lambda}_1),\delta)$ 
and $\theta_1$ be the smallest value in $(0,1)$ such that $(Q(\theta_1), \lambda(\theta_1))$ 
lies at the boundary of $B((\hat{Q}_1, \hat{\lambda}_1),\delta)$ and
in the next ball $B((\hat{Q}_2, \hat{\lambda}_2),\delta)$, 
where $(\hat{Q}_2, \hat{\lambda}_2)\in {\cal Q}^\lambda_J$. 
Next, we let $\theta_2$ be the smallest value in 
$[\theta_1,1)$ such that $(Q(\theta_2), \lambda(\theta_2))$
lies at the boundary of $B((\hat{Q}_2, \hat{\lambda}_2),\delta)$ and in the next ball labelled 
 $B((\hat{Q}_3, \hat{\lambda}_3),\delta)$,
where $(\hat{Q}_3, \hat{\lambda}_3)\in {\cal Q}^\lambda_J$. 
Continuing the process, 
we let $\theta_{J-1}$ be the smallest value in
$[\theta_{J-2},1)$ such that $(Q(\theta_{J-1}), \lambda(\theta_{J-1}))$ lies at the boundary of $B((\hat{Q}_{J-1}, \hat{\lambda}_{J-1}),\delta)$
and in the next ball $B((\hat{Q}_{J}, \hat{\lambda}_{J}),\delta)$,
where 
$(\hat{Q}_{J-1}, \hat{\lambda}_{J-1})\in  {\cal Q}^\lambda_J$,
$(\hat{Q}_{J}, \hat{\lambda}_{J})\in  {\cal Q}^\lambda_J$,
and $(Q'', \lambda'')\in B((\hat{Q}_{J}, \hat{\lambda}_{J}),\delta)$.

Under condition~(a), we can set $\delta$ to be sufficiently small such that
inequality (\ref{eq:Lip-Solution}) holds in each of the ball.
Consequently, by inequality (\ref{eq:f_Lipschitz_Qj}) for all $j=1,\cdots,J$, 
we have
\bgeq
\|{\bm f}_{Q', \lambda'}-{\bm f}_{Q'', \lambda''}\|_k
&\leq& \|{\bm f}_{Q', \lambda'}-{\bm f}_{Q(\theta_1), \lambda(\theta_1)}\|_k+\sum_{j=1}^{J-2}\|{\bm f}_{Q(\theta_j), \lambda(\theta_j)} -{\bm f}_{Q(\theta_{j+1}), \lambda(\theta_{j+1})}\|_k \nonumber  \\
&+&\|{\bm f}_{Q(\theta_{J-1}),\lambda(\theta_{J-1})} -{\bm f}_{Q'',\lambda''}\|_k  \nonumber  \\
&\leq& 2C_{\hat{Q}_1,\hat{\lambda}_1} (\dd_{\rm FM}(Q', Q(\theta_1))+| \lambda'-\lambda(\theta_1) |) \nonumber  \\
&+&\sum_{j=1}^{J-2}2C_{\hat{Q}_{j}, \hat{\lambda}_j}(\dd_{\rm FM}(Q(\theta_j),Q(\theta_{j+1})) + | \lambda(\theta_j)-\lambda(\theta_{j+1}) | ) \nonumber \\
&& +2C_{\hat{Q}_{J},\hat{\lambda}_J}\dd_{\rm FM}(Q(\theta_{J-1}), Q'') +| \lambda(\theta_{J-1})-\lambda'' | \nonumber \\
&\leq& 2\max_{j=1,\cdots,J}C_{\hat{Q}_j,\hat{\lambda}_j} \big[(1-\theta_1)(\dd_{\rm FM}(Q',Q'') + |\lambda' -\lambda''|)\nonumber \\
&+&\sum_{j=1}^{J-2}(\theta_{j+1}-\theta_{j})(\dd_{\rm FM}(Q',Q'')+ |\lambda' -\lambda''|)\nonumber \\
&&+ 
 \theta_{J-1}(\dd_{\rm FM}(Q',Q'') + |\lambda' -\lambda''|)
\big]\nonumber\\
&=& 2C_\kappa(\dd_{\rm FM}(Q',Q'') + |\lambda' -\lambda''|),\nonumber
\edeq
where $C_\kappa:=\max_{j=1,\cdots,J}\left\{C_{\hat{Q}_j,\hat{\lambda}_j}\right\}$.

Part (iv). Note that 
\bgeq
\dd_{\rm FM}(Q_N^1,Q_N^2)
&=& \sup_{h\in {\cal F}_p(Z)} \left|\int_{Z}h({\bm z}) Q_N^1(d{\bm z})-\int_{Z}h({\bm z}) Q_N^2(d{\bm z})\right|\\
&=&  \sup_{h\in {\cal F}_p(Z)} \left|\frac{1}{N}\sum_{l=1}^N h({\bm z}^l_1) 
-\frac{1}{N}\sum_{l=1}^N h({\bm z}^l_2) \right|\\
&\leq & \frac{1}{N}\sum_{l=1}^N |h({\bm z}^l_1)-h({\bm z}^l_2)|
\leq \frac{1}{N}\sum_{l=1}^N \max\{1,\|{\bm z}_1^l\|,\|{\bm z}_2^l\|\}^{p-1} \|{\bm z}_1^l-{\bm z}_2^l\|. 
\edeq
Then by Part (iii), 
\bgeq
\|{\bm f}_{Q_N^1,\lambda_N}-{\bm f}_{Q_N^2,\lambda_N}\|_k  \leq 2C_\kappa \dd_{\rm FM}(Q_N^1,Q_N^2)
\leq \frac{2C_\kappa}{N} 
 \sum_{l=1}^N \max\{1,\|{\bm z}_1^l\|,\|{\bm z}_2^l\|\}^{p-1} \|{\bm z}_1^l-{\bm z}_2^l\|.
\edeq
The proof is complete.
\hfill $\Box$
}

\subsection{Quantitative statistical robustness}
\label{sec:all_data}
 We are now ready to present our desired quantitative 
statistical robustness results.
We present them according to 
the conditions required on the cost function because this determines 
the scope of the applicability of the results.

\subsubsection{Cost function is twice continuously differentiable}
In Theorem \ref{thm:lip_solution-lambda} (iv), 
we have derived an error bound for ${\bm f}_{Q_N,\lambda_N}$
based on two different sets of samples, which allows
us to calculate the deterministic quantity 
$\|{\bm f}_{Q_N^1,\lambda_N} - {\bm f}_{Q_N^2,\lambda_N}\|_k$ for any two given samples.
In this section, we derive an error bound for the difference of the probability distributions of ${\bm f}_{Q_N^1, \lambda_N}$ and ${\bm f}_{Q_N^2,\lambda_N}$
when the samples vary randomly: one is constructed with perceived data and the other is constructed with real data. Differing from (\ref{eq:Lip-Solution-PN}),
the new error bound to be established will allow us to 
estimate the difference between the
 cumulative distribution functions
of the 
two estimators.  The next theorem states this.

\begin{theorem}[Quantitative statistical robustness of 
stationary point]
\label{t:QSR_SGE}
Assume the setting and
conditions of 
Theorem~\ref{thm:lip_solution-lambda}.
Let 
$Z$ be a compact set and
$\widehat{\bm f}(\tilde{\bm z}^1,\cdots,\tilde{\bm z}^N,\lambda_N):={\bm f}_{Q_N,\lambda_N}$
be a statistical estimator of ${\bm f}_{Q_N,\lambda_N}$.
Then
\bgeqn
\label{eq:Quant-SR}
\dd_K(P^{\otimes N}\circ \widehat{{\bm f}}(\cdot,\lambda_N)^{\;-1},
Q^{\otimes N}\circ \widehat{\bm f}(\cdot,\lambda_N)^{\;-1})\leq
 \widetilde{C}_{\kappa} \dd_{K}(P,Q),
 {\color{black} 
\;\forall P, Q\in {\cal M}_{Z,M_2}^{p\gamma}}
\edeqn
for all $\lambda_N\in [{\color{black}\tau},\bar{\lambda}]$ with $N\in \mathbb{N}$,
{\color{black}
where 
{\color{black} 
${\cal M}_{Z,M_2}^{p\gamma}$ is defined as in (\ref{eq:M-phi-bounded}) } and 
$\dd_K$ is the special case of $\dd_{FM}$ defined as in (\ref{eq:FM-Kan}) with $p=1$. 

}
\end{theorem}

\vspace{0.3cm}
\noindent 
\textbf{Proof.}
 Since $Z$ is compact,
there exists some $M>0$ such that
$\mathscr{P}_N\subset {\cal M}_{Z,M_2}^{p\gamma}$,
{\color{black} where $\mathscr{P}_N$ is defined as in Theorem~\ref{thm:lip_solution-lambda}.}
By Part (iv) of Theorem~\ref{thm:lip_solution-lambda},
we have
\bgeq
\|\widehat{{\bm f}}(\tilde{{\bm z}}^1,\cdots, \tilde{{\bm z}}^N,{\lambda_N})-
\widehat{{\bm f}}({\bm z}^1,\cdots,{\bm z}^N,{\lambda_N})\|_k 
&=&
\|{\bm f}_{Q_N,\lambda_N}-{\bm f}_{P_N,\lambda_N}\|_k \nonumber \\
&\leq&
\frac{
2C_{\kappa}
}{N}
 \sum_{l=1}^N \max\{1,\|\tilde{\bm z}^l\|,\|{\bm z}^l\|\}^{p-1} \|\tilde{\bm z}^l-{\bm z}^l\|\leq \frac{\widetilde{C}_{\kappa}}{N}\sum_{l=1}^N\|\tilde{\bm z}^l-{\bm z}^l\| 
\edeq
for any 
$Q_N=\frac{1}{N}\sum_{l=1}^N\delta_{\tilde{\bm z}^l}$,
$P_N=\frac{1}{N}\sum_{l=1}^N\delta_{{\bm z}^l}$,
and $\lambda_N\in [{\color{black}\tau},\bar{\lambda}]$ with $N\in \mathbb{N}$,
where $\widetilde{C}_{\kappa}:=2C_{\kappa}\sup_{{\bm z}\in Z}\max\{1,\|{\bm z}\|\}^{p-1}$. 
Let ${\cal G}$ 
be a set of Lipshitz continuous functional $g:{\cal H}_k\rightarrow \R$ with its modulus is $1$.
Then
\bgeq
\left|g\left(
\widehat{{\bm f}}(\tilde{{\bm z}}^1,\cdots, \tilde{{\bm z}}^N,{\lambda_N})
\right) -
g\left(\widehat{{\bm f}}({{\bm z}}^1,\cdots, {{\bm z}}^N,{\lambda_N})\right)\right|
&\leq& 
\left\|
\widehat{{\bm f}}(\tilde{{\bm z}}^1,\cdots, \tilde{{\bm z}}^N,{\lambda_N})
-
\widehat{{\bm f}}({{\bm z}}^1,\cdots, {{\bm z}}^N,{\lambda_N})
\right\|_k\\
&\leq& \frac{\widetilde{C}_{\kappa} }{N} \sum_{i=1}^{N} \|\tilde{{\bm z}}^i-{\bm z}^i\|.
\edeq

By \cite[Lemma~1]{GuX20-SR},
\bgeq
&&\dd_K\left(P^{\otimes N}\circ 
\widehat{{\bm f}}(
\cdot,{\lambda_N})
^{\;-1},
Q^{\otimes N}\circ 
\widehat{\bm f}(
\cdot,{\lambda_N})
^{\;-1}
\right)\\
&=&\sup_{g \in {\cal G}}\left|
\int_{{\cal H}_k} g(t) P^{\otimes N}\circ 
\widehat{\bm f}(
\cdot,{\lambda_N})
^{\;-1}
(dt) -
\int_{{\cal H}_k} g(t) Q^{\otimes N}\circ 
\widehat{\bm f}(
\cdot,{\lambda_N})
^{\;-1}
(dt)
\right|\\
&=&
\sup_{g \in {\cal G}}\left|
\int_{Z^{\otimes N}}g\left(
\widehat{\bm f}(
\vec{{\bm z}}^N,{\lambda_N})
\right) P^{\otimes N}
\left(d\vec{{\bm z}}^N\right) -
\int_{Z^{\otimes N}} g\left(
\widehat{\bm f}(
\vec{{\bm z}}^N,{\lambda_N})
\right) Q^{\otimes N}
\left(d\vec{{\bm z}}^N\right)
\right|\\
&\leq &
 \sum_{i=1}^N 
 \frac{\widetilde{C}_{\kappa}}{N} \dd_{K} (P,Q)
= \widetilde{C}_{\kappa} \dd_{K} (P,Q)
\edeq
for all 
$\lambda_N\in [{\color{black}\tau},\bar{\lambda}]$ with
$N\in \mathbb{N}$.
The proof is complete.
\hfill $\Box$

{\color{black} 
The 
theoretical result is useful in at least two cases: (a) $P$ and $Q$ are known but in actual calculations, only empirical data of $Q$ are used. This is either because
 errors 
occurring in the process of data generation and processing, or the distribution of validation data (for the future)  
is shifted from the distribution of training data (in the past); (b) the difference between $Q$ and $P$
is known in the sense that
the shift of the distribution   
is within a controllable range.  
In the case only when perceived data (the sample data of $Q$) is known, we will not be able to say much about the quality of the learning estimator.
Since the result is built on Theorem \ref{thm:lip_solution-lambda}, 
it might be desirable 
to relax some of the conditions imposed on 
Theorem~\ref{thm:lip_solution-lambda} such as strong regularity to extend the applicability of 
Theorem~\ref{t:QSR_SGE}.
}

\begin{remark}
\label{rem:calM}
Before concluding this subsection,
we note in both the qualitative and quantitative statistical robustness,
the perturbation is under some topology of probability measures in ${\cal M}^{\psi^{\gamma}}_{Z,M_1}$ and ${\cal M}^{p\gamma}_{Z,M_2}$ respectively, 
which
 posts some restriction on the tail distribution and thus
are more restrictive than the usual weak topology.

\begin{itemize}

\item[(i)] 
The proof of quantitative statistical robustness %
in Theorem~\ref{t:QSR_SGE}
requires the statistical estimator $\widehat{\bm f}(\vec{\bm z}^N,\lambda_N)$ 
to be globally Lipschitz continuous
in $\vec{\bm z}^N$.
 The global Lipschitz continuity
is derived by virtue of 
the finite covering theorem to be applied over a relatively closed line segment connecting two probability measures in the set ${\cal M}^{p\gamma}_{Z,M_2}$. 
This explains why we require 
${\cal M}^{p\gamma}_{Z,M_2}$ to be 
relatively compact,
see \cite[Lemma 2.69]{Clau16}.

\item[(ii)] The qualitative statistical robustness presented in Theorem~\ref{t-ML-SR-QL-SR} for $\widehat{\bm f}({\vec{\bm z}^N,\lambda_N})$,
where the probability is restricted to
set 
$
{\cal M}_{Z,M_1}^{\psi^\gamma},
$
because we need to use  
Uniform Glivenko–Cantelli property, see \cite[Corollary 3.5]{KSZ12}.
\end{itemize}
\end{remark}

\subsubsection{Cost function is merely  continuously differentiable}

One of the key conditions in Theorem~\ref{t:QSR_SGE} is 
 the second order continuous differentiability of the cost function $c$
 which might be undesirable in some 
 practical applications.
However, we may get rid of this 
under the circumstances when
the optimal solution
${\bm f}_{Q_N,\lambda_N}$ 
of problem (\ref{eq:ML-saa-r-equiv})
lies in the interior of ${\cal F}$ 
and $c$ is convex w.r.t. the second argument.
The next theorem addresses this.

\begin{theorem}
\label{thm:Qunt-SR-New}
Let $\psi$ be defined as in (\ref{eq:growth-ml}) and  $M_3>0$ be a positive constant. Define the set of probability distributions 
\bgeq
{\cal M}_{Z,M_3}:=\left\{P'\in \mathscr{P}(Z): \int_{Z} \psi({\bm z}) P'(d{\bm z})\leq M_3\right\}
\edeq
Assume: (a) {\color{black} ${\rm (C1)}$, ${\rm (C2)}$, ${\rm (C3)'}$ and ${\rm (C4)'}$ hold 
with $\sup_{{\bm z}\in Z} L_p({\bm z},{\bm z})$ being
bounded;
(b) ${\rm (K1)'}$ hold;
}
(c) $\beta>\sqrt{\frac{M}{\tau}}$,
where $\beta$ is defined as in Theorem~\ref{t-ML-SR-QL-SR}, (d) For $N$ sufficiently large,
$P_N, Q_N\in {\cal M}_{Z,M_3}$.
Then
\bgeqn
\label{eq:Quant-SR-A}
\dd_K(P^{\otimes N}\circ \widehat{{\bm f}}(\cdot,\lambda_N)^{\;-1},
Q^{\otimes N}\circ \widehat{\bm f}(\cdot,\lambda_N)^{\;-1})\leq
 \widetilde{C}_{\kappa} \dd_{K}(P,Q), \forall P,Q  \in {\cal M}_{Z,M_3}
\edeqn
for all $\lambda_N\geq \tau$, 
with $N\in \mathbb{N}$ sufficiently large, where $\tau$ is any fixed positive
number.
\end{theorem}

\noindent
\textbf{Proof.} Let $N$ be sufficiently large such that condition (d) is satisfied,
and $\vec{\bm z}^N$ and $\lambda_N\geq \tau$ be 
fixed. 
We proceed the proof in three steps.

\underline{Step 1}. We show that the optimal solution ${\bm f}_{Q_N,\lambda_N}$ lies in the interior of 
${\cal F}$.
To see this, we note that
 $c({\bm z},{\bm f}({\bm x}))\geq 0$ for any ${\bm z}\in Z$ and ${\bm f}\in {\cal H}_k$. Thus
\bgeqn
\label{eq:ineq_R}
R_{Q_N}^{\lambda_N}({\bm f}) =\bbe_{Q_N}[c({\bm z},{\bm f}({\bm x}))]
 + \lambda_N \|{\bm f}\|_k^2
\geq \lambda_N \|{\bm f}\|_k^2
\geq \tau \|{\bm f}\|_k^2. 
\edeqn
Moreover, it follows from 
{\color{black}
Assumption~\ref{A:cost-1}${\rm (C2)}$} 
that $c({\bm z},0)\leq \psi({\bm z})$ for all ${\bm z}\in Z$.
Consequently,
\bgeq 
R_{Q_N}^{\lambda_N}(0) =\bbe_{Q_N}[c({\bm z},0)]
\leq \bbe_{Q_N}[\psi({\bm z})]
\leq M_3.
\edeq
Note also that
$$
R_{Q_N}^{\lambda_N}({\bm f}_{Q_N,\lambda_N})\leq R_{Q_N}^{\lambda_N}(0) \leq M_3
$$
because ${\bm f}_{Q_N,\lambda_N}$ is optimal.
Combining with (\ref{eq:ineq_R}),
we obtain 
$
\tau \|{\bm f}_{Q_N,\lambda_N}\|_k^2\leq M_3
$
and thus
$
\|{\bm f}_{Q_N,\lambda_N}\|_k\leq \sqrt{\frac{M}{\tau}}.
$
Under condition (c), 
this implies that 
$
\|{\bm f}_{Q_N,\lambda_N}\|_k <\beta,
$
which means the optimal solution ${\bm f}_{Q_N,\lambda_N}$ lies in the interior of 
${\cal F}$.

\underline{Step 2}. By Theorem~\ref{eq:first-order-optimal},
the optimality condition of problem (\ref{eq:ML-saa-r-equiv}) is 
$$
0\in \psi(\vec{\bm z}^N,{\bm f})
 +{\cal N}_{\cal F}({\bm f}), {\color{black}\quad\forall Q_N\in \hat{\cal P}},
$$
{\color{black}
where $\hat{\cal P}$ is defined as in (\ref{eq:set_Q}) and}  $\psi(\vec{\bm z}^N,{\bm f}):=
\frac{1}{N} \sum_{l=1}^N 
 c'_2 ({\bm z}^l,{\bm f}({\bm x}^l)) k_{{\bm x}^l}+ 2\lambda_N{\bm f}$.

Then $\psi(\vec{\bm z}^N,{\bm f}_{Q_N,\lambda_N})
=0$,
which means that the 
optimality condition reduces to
\bgeqn 
\psi(\vec{\bm z}^N,{\bm f}) =0
\label{eq:Opi-cond-no-cone}
\edeqn 
in this case.
Observe that $Q_N$ is determined by $\vec{\bm z}^N$. So we 
can write $\widehat{\bm f}(\vec{\bm z}^N,\lambda_N)$ for
${\bm f}_{Q_N,\lambda_N}$. 
Moreover, since $\lambda_N$ is fixed, we may write  
$\widehat{\bm f}(\vec{\bm z}^N)$ 
for $\widehat{\bm f}(\vec{\bm z}^N,\lambda_N)$ to ease the exposition. 

Let ${\cal B}_{\vec{\bm z}^N}\subset Z^{\otimes N}$ be 
a neighborhood of $\vec{\bm z}^N$. We use the
implicit function theorem, 
Theorem~\ref{thm:implicit_thm}
in Appendix~\ref{app:Implicit},
to show that equation (\ref{eq:Opi-cond-no-cone})
defines an implicit function $\widehat{\bm f}(\cdot)$ 
mapping from ${\cal B}_{\vec{\bm z}^N}$ to 
${\cal B}_{\widehat{\bm f}}$, where
${\cal B}_{\widehat{\bm f}}$ denotes a neighborhood of $\widehat{\bm f}(\vec{\bm z}^N)$  and 
the mapping is locally Lipschitz continuous 
w.r.t. variation of $\vec{\bm z}^N$ in ${\cal B}_{\vec{\bm z}^N}$.
It suffices to verify the conditions of 
the theorem.
Observe first that 
$\psi(\vec{\bm z}^N,
{\bm f}(\vec{\bm z}^N))=0$
and $\psi(\cdot,\cdot)$ is continuous 
in a neighborhood
of point $(\vec{\bm z}^N,
\widehat{\bm f}(\vec{\bm z}^N))$, written ${\cal B}_{\vec{\bm z}^N}\times {\cal B}_{\widehat{\bm f}}$.
Second, since $c({\bm z},\cdot)$ is convex for every ${\bm z}$, then
$\frac{1}{N} \sum_{l=1}^N 
 c'_2 ({\bm z}^l,\langle \cdot, k_{{\bm x}^l}\rangle ) k_{{\bm x}^l}$
is monotone. Moreover, since $\lambda_N>0$, then
$\psi(\vec{\bm z}^N,\cdot)$ is strongly monotone,
that is,
\bgeq
 \langle \psi(\vec{\bm z}^N, {\bm f}_1)-\psi(\vec{\bm z}^N,{\bm f}_2), {\bm f}_1-{\bm f}_2\rangle 
 \geq \lambda_N \|{\bm f}_1-{\bm f}_2\|_k^2.
 \edeq
Third, 
we can show that $\psi(\cdot,{\bm f})$ is Lipschitz continuous on
${\cal B}_{\bm z}$ uniformly for ${\bm f} \in {\cal B}_{\widehat{\bm f}}$.
To see this, for any ${\bm f} \in {\cal B}_{\widehat{\bm f}}$ and $\vec{\bm z}^N_{1},\vec{\bm z}^N_{2}\in 
Z$, we have 
\bgeq
&&\|\psi(\vec{\bm z}_1^N,{\bm f})-\psi(\vec{\bm z}_2^N,{\bm f})\|_k\\
&=&   \Big\|\frac{1}{N} \sum_{l=1}^N 
 c'_2 ({\bm z}^{1,l},{\bm f}({\bm x}^{1,l})) k_{{\bm x}^{1,l}}+ 2\lambda_N{\bm f}
 -\frac{1}{N} \sum_{l=1}^N 
 c'_2 ({\bm z}^{2,l},{\bm f}({\bm x}^{2,l})) k_{{\bm x}^{2,l}}- 2\lambda_N{\bm f}
 \Big \|_k\\
 &\leq & \frac{1}{N}\sum_{l=1}^N \|c'_2 ({\bm z}_{1,l},{\bm f}({\bm x}_{1,l})) k_{{\bm x}_{1,l}}-c'_2 ({\bm z}_{2,l},{\bm f}({\bm x}_{2,l})) k_{{\bm x}_{2,l}}\|_k\\
 &\leq &\frac{1}{N}\sum_{l=1}^N L_{c'}\|{\bm z}_{1,l}-{\bm z}_{2,l}\|+\|{\bm f}({\bm x}_{2,l})-{\bm f}({\bm x}_{1,l})\|_k\\
 &\leq &\frac{1}{N}\sum_{l=1}^N L_{c'}(\|{\bm z}_{1,l}-{\bm z}_{2,l}\|+\|{\bm f}\|_k\|k_{{\bm x}_{2,l}}-k_{{\bm x}_{1,l}}\|_k)\\
 &\leq &\frac{1}{N}\sum_{l=1}^N L_{c'}(\|{\bm z}_{1,l}-{\bm z}_{2,l}\|+\beta\|k_{{\bm x}_{2,l}}-k_{{\bm x}_{1,l}}\|_k)\\
  &\leq &\frac{1}{N}\sum_{l=1}^N L_{c'}(1+\beta L_k)\|{z}_{1,l}-{z}_{2,l}\|\\
&=&  \frac{1}{N} L_{c'}(1+\beta L_k)\|\vec{\bm z}_{1}-\vec{\bm z}_{2}\|,
\edeq
where $L_{c'}:=C_1 \sup_{{\bm z}\in Z} L_p({\bm z},{\bm z})$ is bounded.
Thus, by Theorem~\ref{thm:implicit_thm}
in Appendix~\ref{app:Implicit}, we claim that
there exists a unique implicit function $\widehat{\bm f}:
{\cal B}_{\bm z}\to 
{\cal B}_{{\bm f}}$ 
such that 
$$
\psi(\vec{\bm z}^N,{\bm f}(\vec{\bm z}^N)) =0,\;\; \forall \vec{\bm z}^N \in {\cal B}_{\bm z}
$$
and 
\bgeqn
\|\widehat{\bm f}(\vec{\bm z}^N_1)-\widehat{\bm f}(\vec{\bm z}^N_2)\|_k\leq \frac{L_{c'}(1+\beta L_k)}{N \lambda_N}\|\vec{\bm z}^N_1-\vec{\bm z}^N_2\|, \forall 
\vec{\bm z}^N_1, \vec{\bm z}^N_2\in {\cal B}_{\bm z}.
\label{eq:f_N-local-Lip}
\edeqn

\underline{Step 3}. We show that 
the implicit function can be extended over 
$Z^{\otimes N}$ and 
inequality (\ref{eq:f_N-local-Lip}) holds for any $\vec{\bm z}^N_1, 
\vec{\bm z}^N_2\in Z^{\otimes N}$.
Observe that the optimal solution $\widehat{\bm f}(\vec{\bm z}^N)$ exists for any  $\vec{\bm z}^N\in Z^{\otimes N}$.
To ease the exposition,
we will write $\vec{\bm z}$ for $\vec{\bm z}^N$.

We will use the finite covering theorem to
show the global Lipschitz continuity (with the same Lipschitz modulus).
Let 
$\vec{\bm z}', \vec{\bm z}''\in Z^{\otimes N}$ be any
two fixed vectors. For $\theta\in [0,1]$, let $\vec{\bm z}(\theta):=(1-\theta)\vec{\bm z}'+\theta \vec{\bm z}''$ and 
$$
Z_{[\vec{\bm z}', \vec{\bm z}'']}
=\left\{\vec{\bm z}(\theta): \theta\in [0,1]\right\}
$$
be the line segment connecting $\vec{\bm z}', \vec{\bm z}''$. For any fixed point, $\vec{\bm z}(\theta) \in Z_{[\vec{\bm z}', \vec{\bm z}'']}$, we 
can use the implicit function theorem detailed 
in Step 2 to show that the equation
\bgeq
\psi(\vec{\bm z}(\theta),{\bm f}) =0
\edeq
defines a unique implicit function in a $\delta$-neighborhood
of $\vec{\bm z}(\theta)$, denoted by 
${B}(\vec{\bm z}(\theta),\delta)$ 
and 
\bgeq
\|\widehat{\bm f}(\vec{\bm z}^{(1)})-\widehat{\bm f}(\vec{\bm z}^{(2)})\|_k\leq \frac{L_{c'}(1+\beta L_k)}{N\lambda_N}\|\vec{\bm z}^{(1)}-\vec{\bm z}^{(2)}\|, \forall 
\vec{\bm z}^{(1)}, \vec{\bm z}^{(2)}\in {B}(\vec{\bm z}(\theta),\delta),
\edeq
where ${B}(\vec{\bm z}_j,\delta)$ denotes a closed ball centered at $\vec{\bm z}_j$ with radius $\delta$ under the distance $\|\cdot\|$ over $Z^{\otimes N}$.
By the finite covering theorem, we can construct a $\delta$-net $\vec{\bm z}_1,\cdots, \vec{\bm z}_J$ 
over $Z_{[\vec{\bm z}', \vec{\bm z}'']}$ 
such that (i) $Z_{[\vec{\bm z}', \vec{\bm z}'']}\subset \cup_{j=1}^J {B}(\vec{\bm z}_j,\delta)$, (ii) a unique implicit function is well-defined over each of the balls and (iii) 
these implicit functions are connected to form a continuous function over $\cup_{j=1}^J {B}(\vec{\bm z}_j,\delta)$.
Following a similar argument to that in Step 1, we can show that 
\bgeqn
\label{eq:local_Lip}
\|\widehat{\bm f}(\vec{\bm z}^{(1)}_j)-\widehat{\bm f}(\vec{\bm z}^{(2)}_j, \lambda_N)\|\leq \frac{L_{c'}(1+\beta L_k)}{N\lambda_N}\|\vec{\bm z}^{(1)}_j-\vec{\bm z}^{(2)}_j\|, \forall \vec{\bm z}^{(1)}_j,\vec{\bm z}^{(2)}_j\in {B}(\vec{\bm z}_j,\delta),~~j=1,\cdots,J.
\edeqn

Let $\hat{\vec{\bm z}}_1\in \{\vec{\bm z}_j\}_{j=1}^J$ 
be such that $\vec{\bm z}'\in {B}(\hat{\vec{\bm z}}_1,\delta)$ 
and $\theta_1$ be the smallest value in $(0,1)$ such that $\vec{\bm z}(\theta_1)$ 
lies at the boundary of $B(\hat{\vec{\bm z}}_1, \delta)$ and
in the next ball ${B}(\hat{\vec{\bm z}}_2,\delta)$, 
where $\hat{\vec{\bm z}}_2 \in \{\vec{\bm z}_j\}_{j=1}^J$. 
Next, we let $\theta_2$ be the smallest value in 
$[\theta_1,1)$ such that $\vec{\bm z}(\theta_2)$
lies at the boundary of ${B}(\hat{\vec{\bm z}}_2, \delta)$ and in the next ball labelled 
 ${B}(\hat{\vec{\bm z}}_3,\delta)$,
where $\hat{\vec{\bm z}}_3\in \{\vec{\bm z}_j\}_{j=1}^J$. 
Continuing the process, 
we let $\theta_{J-1}$ be the smallest value in
$[\theta_{J-2},1)$ such that $\vec{\bm z}(\theta_{J-1})$ lies at the boundary of ${B}(\hat{\vec{\bm z}}_{J-1},\delta)$
and in the next ball ${B}(\hat{\vec{\bm z}}_J,\delta)$,
where 
$\hat{\vec{\bm z}}_{J-1} \in \{\vec{\bm z}_j\}_{j=1}^J$,
$\hat{\vec{\bm z}}_{J}\in  \{\vec{\bm z}_j\}_{j=1}^J$,
and $\vec{\bm z}''\in {B}(\hat{\vec{\bm z}}_{J}, \delta)$.

We  can set $\delta$ to be sufficiently small such that
inequality (\ref{eq:local_Lip}) holds in each of the ball.
Consequently, 
we have
\bgeq
\|\widehat{\bm f}(\vec{\bm z}')
-\widehat{\bm f}(\vec{\bm z}'')\|_k
&\leq& \|\widehat{\bm f}(\vec{\bm z}')
-\widehat{\bm f}(\vec{\bm z}(\theta_1))\|_k+\sum_{j=1}^{J-2}\|\widehat{\bm f}(\vec{\bm z}(\theta_j)) 
-\widehat{\bm f}(\vec{\bm z}(\theta_{j+1}))\|_k \nonumber  \\
&+&\|\widehat{\bm f}(\vec{\bm z}(\theta_{J-1})) -{\bm f}(\vec{\bm z}'')\|_k  \nonumber  \\
&\leq& \frac{L_{c'}(1+\beta L_k)}{N\lambda_N} \|\vec{\bm z}'-\vec{\bm z}(\theta_1)\| 
+\sum_{j=1}^{J-2}\frac{L_{c'}(1+\beta L_k)}{N\lambda_N} \|\vec{\bm z}(\theta_j)-\vec{\bm z}(\theta_{j+1})\|\\
&& \frac{L_{c'}(1+\beta L_k)}{N\lambda_N}
\|\vec{\bm z}(\theta_{J-1})-\vec{\bm z}''\| \nonumber\\
&\leq& 
\frac{L_{c'}(1+\beta L_k)}{N\lambda_N}
\big[(1-\theta_1) \|\vec{\bm z}'-\vec{\bm z}''\|
+\sum_{j=1}^{J-2}(\theta_{j+1}-\theta_{j})
 \|\vec{\bm z}'-\vec{\bm z}''\|\\
&&+ 
 \theta_{J-1}  \|\vec{\bm z}'-\vec{\bm z}''\|
\big]\nonumber\\
&=&\frac{L_{c'}(1+\beta L_k)}{N\lambda_N}
\|\vec{\bm z}'-\vec{\bm z}''\|.
\edeq
The rest is analogous to the proof of Theorem~\ref{t:QSR_SGE}, we omit the details.
\hfill
\Box

Theorem~\ref{thm:Qunt-SR-New} differs
from Theorem~\ref{t:QSR_SGE} on several aspects.
First, inequality (\ref{eq:Quant-SR-A}) holds for 
all $N$ sufficiently large as opposed to  all $N$
as in (\ref{eq:Quant-SR}). This is because only when
$N$ is sufficiently large, we will be able to ensure  
${\bm f}_{Q_N,\lambda_N}$ to lie in the interior of ${\cal F}$ under condition (d). Second, 
condition (d) is not very demanding because 
when $N$ is sufficiently large, we know by the law of large numbers that 
$\bbe_{Q_N}[\psi({\bm z})]$ is close to $\bbe_Q[\psi({\bm z})]$ with probability $1$.
So condition (d) is satisfied with high probability 
when $Q\in {\cal M}_{Z,M_3/2}$.
Third, the local Lipschitz continuity of $\widehat{\bm f}(\cdot,\lambda_N)$
is derived in a completely different manner from
that of Theorem \ref{thm:lip_solution-lambda}~(ii). This is because we use the new implicit function theorem, Theorem~\ref{thm:implicit_thm}, as opposed to Lemma~\ref{lem:Lipschitz-perturbed-phiz}. One of the main advantages is that we no longer require
strong regularity condition and replace it with 
strong monotone of $\psi(\vec{\bm z}^N,{\bm f})$. A sufficient condition for this is the cost function $c({\bm z},{\bm f}({\bm x}))$ is convex in ${\bm f}$ {\color{black}in ${\rm (C1)}$}.
Moreover, we no longer require second order continuous differentiability of the cost function.
Fourth, when $\psi({\bm z})=\|{\bm z}\|^p$, we can compare
set $
{\cal M}_{Z,M_2}^{p\gamma}$ with 
set ${\cal M}_{Z,M_3}$. 
Let $M_2=M_3$.
Then we can see that 
${\cal M}_{Z,M_2}^{p\gamma}$ is only a subset of 
${\cal M}_{Z,M_3}$ which means the quantitative statistical robustness result is applicable to 
a large set of probability distributions (and hence 
a larger class of data sets).

\section{Single data perturbation}
\label{Sec:Single-data-IF}
{\color{black}
We now move to discuss  
 the impact of single data 
 perturbation on the optimal solution of the 
 regularized problem 
$
\min_{{\bm f}\in {{\cal F}}} \bbe_{P}[c({\bm z},{\bm f}({\bm x}))]
+ \lambda \|{\bm f}\|_k^2.
$
}
To this end, we consider the perturbation of the true probability distribution $P$ by a Dirac
distribution $\delta_{\tilde{{\bm z}}}$ at a perturbed data point $\tilde{{\bm z}}$.
For $\tilde{{\bm z}}\in \R^n$ and $t\in (0,1)$, let $(1-t)P+t\delta_{\tilde{{\bm z}}}$ denote the mixture distribution between $P$ and $\delta_{\tilde{{\bm z}}}$.
Consider
\bgeqn
\label{eq:single-perturbation-KKT}
0\in \bbe_{(1-t)P+t\delta_{\tilde{{\bm z}}}}[c'_2({\bm z},{\bm f}({\bm x})){\color{black}k_{\bm x}}]+2\lambda {\bm f}+{\cal N}_{\cal F}({\bm f}).
\edeqn
In this section, we assume ${\rm (C1)}$ hold and thus the solution of 
$0\in \bbe_{P}[c'_2(({\bm z},{\bm f}({\bm x})){\color{black}k_{\bm x}}]+2\lambda {\bm f}
+{\cal N}_{\cal F}({\bm f})$ and 
(\ref{eq:single-perturbation-KKT}) is a singleton.
Let 
${\bm f}_{(1-t)P+t\delta_{\tilde{{\bm z}}},\lambda}$
denote the 
solution of (\ref{eq:single-perturbation-KKT}).
We will investigate the impact 
of probability perturbation on the solution by the so-called influence function of
${\bm f}_{(\cdot,\lambda)}$.
\begin{definition}
(Influence function, \cite[Definition 10.4]{StC08})
\label{def:IF}
Let  $\tilde{{\bm z}} \in Z$.
The influence function ${\rm IF}:Z\rightarrow {\cal F}$ of ${\color{black}{\bm f}_{(\cdot,\lambda)}}:\mathcal{P}(Z) \rightarrow {\cal F}$ at a point $\tilde{{\bm z}}$ for a distribution $P\in \mathcal{P}(Z)$ is given by 
\bgeq
{\rm IF}\left(\tilde{{\bm z}};{\bm f}_{(\cdot,\lambda)},P\right):=\lim_{t\downarrow 0} \frac{{{\bm f}}_{(1-t)P+t\delta_{\tilde{{\bm z}}},\lambda}-{{\bm f}}_{P,\lambda}}{t}
\edeq
provided that the limit exists.
\end{definition}

{\color{black}The influence function is closely related to the directional derivative at $P$ along
direction $\delta_{\tilde{{\bm z}}}-P$.
It quantifies the sensitivity (instead of stability) of 
{\color{black}
${\bm f}_{P,\lambda}$ 
}
when the true probability $P$ 
is perturbed along the direction.
From a data perspective, 
the mixture distribution means
that there is a 
probability $t$, the sample data could be an outlier $\tilde{\bm z}$
and the influence function allows one to quantify the sensitivity of the optimal solution w.r.t. the change. 
}
In general, 
a closed form of ${\bm f}$ 
is not obtainable.
However, since ${\bm f}$  satisfies the first order optimality condition
(\ref{eq:ML-Rf-min-opti-KKT}),
we may 
compute 
{\color{black}
$
{\rm IF}\left(\tilde{{\bm z}};{{\bm f}}_{(\cdot,\lambda)},P\right)$
}
by applying the
implicit function theorem 
to (\ref{eq:ML-Rf-min-opti-KKT}).
Indeed, this is what
Steinwart and Christmann \cite{StC08} do
for the case when $\F$ is the whole space.
The main technical challenge here is that
$\F$ is not necessarily the whole space.
In this case, the implicit function theorem
would have to involve differentiation 
of the normal cone 
${\cal N}_{\cal F}({\bm f})$.
To this end, we use the proto-derivative of set-valued mapping 
in the Hilbert space, see e.g.~\cite{
AdR21}.

Let $H$ be a Hilbert space.
Recall that a set-valued mapping $\Gamma: H \rightrightarrows H$ is said to be proto-differentiable at a point ${\bm f}$ and for a particular element ${\bm h}\in \Gamma({\bm f})$ if the set-valued mappings 
$
\Delta_{{\bm f},{\bm h},\epsilon}: d \rightarrow \frac{\Gamma({\bm f}+\epsilon {\bm d})-{\bm h}}{\epsilon},
$
regarded as a family indexed by $\epsilon>0$, 
graph-converge as $\epsilon \downarrow 0$.
The limit,
if exits, 
is denoted by $\Gamma'_{{\bm f},{\bm h}}$ and called the {\em proto-derivative} of $\Gamma$ at ${\bm f}$ for ${\bm h}$.
The graph-convergence of $\Delta_{{\bm f},{\bm h},\epsilon}$ means that 
$$
\limsup_{\epsilon \downarrow 0}\gph(\Delta_{{\bm f},{\bm h},\epsilon})=\liminf_{\epsilon \downarrow 0} \gph(\Delta_{{\bm f},{\bm h},\epsilon})=\gph(\Gamma'_{{\bm f},{\bm h}}),
$$
where {\color{black}
$ \limsup_{t\to t_0} A(t):= \{a \in H\times H: \exists t_k\to t_0, \exists  a_k \to a \inmat{ with }a_k \in  A(t_k)\}$,
$ \liminf_{t\to t_0} A(t):= \{a \in H\times H: \forall t_k\to t_0, \exists {\cal N}\in {\cal N}_{\infty}, a_k \underset{{\cal N}}{\to} a \inmat{ with }a_k \in  A(t_k)\}$,
${\cal N}_{\infty}:=\{{\cal N}\subset \mathbb{N}:\mathbb{N}\backslash{{\cal N}} \inmat{ finite}\}$,
see \cite[Chapter 5]{RoW98},
and
}
\bgeq
\gph(\Delta_{{\bm f},{\bm h},\epsilon})
= \left\{(d,g)\in H\times H| g\in \frac{\Gamma({\color{black}{\bm f}}+\epsilon d)-h}{\epsilon} \right\}
= \frac{\gph{\Gamma}-({\bm f},{\bm h})}{\epsilon},
\edeq
see e.g.~\cite[Definition 5.32]{RoW98}
and \cite{Chi92}.
In the case that the proto-derivative exits, by \cite[Proposition 2.3]{Roc89}, we have
\bgeq
 \gph(\Gamma'_{{\bm f},{\bm h}})
=\limsup_{\epsilon \downarrow 0} \gph(\Delta_{{\bm f},{\bm h},\epsilon})
=\left\{({\bm d},{\bm g})\in H\times  H:{\bm g}\in \limsup_{\epsilon \downarrow 0, {\bm d}'\rightarrow {\bm d}} \frac{\Gamma({\bm f}+\epsilon {\bm d}')-{\bm h}}{\epsilon} \right\},
\edeq
which means the proto-derivative of 
$\Gamma$ equals to
\bgeqn
\label{eq:calculate-Gamma'}
\Gamma'_{{\bm f},{\bm h}}({\bm d})=\limsup_{\epsilon \downarrow 0, {\bm d}'\rightarrow {\bm d}} \frac{\Gamma({\bm f}+\epsilon {\bm d}')-{\bm h}}{\epsilon}.
\edeqn

A closely related concept is {\em graphical derivative},
see e.g.~\cite[Definition 8.33]{RoW98},
which is defined by the outer limit of $(\gph(\Delta_{{\bm f},{\bm h},\epsilon}))_{\epsilon>0}$ (the graphical outer limit of the set-valued mappings $\{\Delta_{{\bm f},{\bm h},\epsilon}: \epsilon >0\}$,
denoted by $D\Gamma({\bm f}|{\bm h})$,
that is, 
\bgeq
\gph(D\Gamma({\bm f}|{\bm h}))
:=\limsup_{\epsilon \downarrow 0}\gph(\Delta_{{\bm f},{\bm h},\epsilon})
=\limsup_{\epsilon \downarrow 0, {\bm d}'\rightarrow {\bm d}} \frac{\Gamma({\bm f}+\epsilon {\bm d}')-{\bm h}}{\epsilon}.
\edeq
The graphical derivative always exists and its graph is the {\em contingent (Bouligand) cone} of
$\gph \Gamma$ at point $({\bm f},{\bm h})$, see the definition in \cite[page xv]{BoS00}.
By the definitions of the proto-derivative and graphical derivative,
they are equal when the proto-derivative exists.
When $\Gamma$ is single-valued, the
proto-derivative $\Gamma'_{{\bm f},{\bm h}}({\bm d})$ with ${\bm h}=\Gamma({\bm f})$ 
reduces to a directional derivative in which case 
we denote it 
by $\Gamma'({\bm f};{\bm d})$ 
to ease the exposition.

Let 
$$
\Phi_{\tilde{\bdz},P,{\color{black}\lambda}}({\bm f},t):= \bbe_{P}[c'_2 ({\bm z},{\bm f}({\bm x})){\color{black}k_{\bm x}}]+t( c'_2 (\tilde{{\bm z}},{\bm f}(\tilde{{\bm x}}))k_{\tilde{\bm x}}-\bbe_{P}[ c'_2 ({\bm z},{\bm f}({\bm x})){\color{black}k_{\bm x}}]){\color{black}+2\lambda {\bm f}}.
$$
Then 
(\ref{eq:single-perturbation-KKT}) 
can be succinctly written as 
\bgeqn
0\in \Phi_{\tilde{{\bm z}},P,{\color{black}\lambda}}({\bm f},t)+{\cal N}_{\cal F}({\bm f}).
\label{eq:sngl-Phi-t}
\edeqn
{\color{black}
Differing from Theorem~\ref{thm:lip_solution-lambda},
here we investigate the derivative of the implicit function ${\bm f}$
defined in (\ref{eq:sngl-Phi-t})
when the true probability $P$ 
is perturbed along direction $\delta_{\tilde{{\bm z}}}-P$.
}
For fixed $\tilde{\bm z}$, $P$ {\color{black}and ${\lambda}$},
we assume that (\ref{eq:sngl-Phi-t})
has a unique solution 
${\bm f}_{(1-t)P+t\delta_{\tilde{{\bm z}}},{\color{black}\lambda }}$ parameterized by $t$, 
denoted by $\tilde{\bm f}_t$,
i.e., 
\bgeq
\tilde{\bm f}_t
 :={\bm f}_{(1-t)P+t\delta_{\tilde{{\bm z}}},{\color{black}\lambda }}.
\edeq
By the definition of proto-derivative,
the proto-derivative of 
 $\tilde{\bm f}_{t}
:t \mapsto {\cal H}_k$ at $t=0$ 
(single-valued)
for {\color{black} ${\bm f}_{P,\lambda}$}
reduces to 
the directional derivative of $\tilde{\bm f}_{(\cdot)}$ at $0$ with direction $s$,
i.e.,
$$
\tilde{\bm f}'_{0;s}=\limsup_{\epsilon \downarrow 0, s'\rightarrow s}\;({\bm f}_{(1- \epsilon s')P+\epsilon s' \delta_{\tilde{\bm z}},{\color{black}\lambda}}-{\bm f}_{P,{\color{black}\lambda }})/\epsilon,
$$
which is equal to the directional derivative of 
{\color{black}${\bm f}_{(\cdot,\lambda)}$}
at $P$ for direction $\delta_{\tilde{\bm z}}-P$,
that is, equals to the influence function 
${\rm IF}(\tilde{{\bm z}};{\bm f}_{(\cdot,{\color{black}\lambda})},P)$,
that is,
$$
{\rm IF}(\tilde{{\bm z}};{\bm f}_{(\cdot,{\color{black}\lambda})},P) ={\bm f}_{(\cdot,{\color{black}\lambda})}'(P;\delta_{\tilde{\bm z}}-P)=\limsup_{\epsilon \downarrow 0, s'\rightarrow s}\;({\bm f}_{(1- \epsilon s')P+\epsilon s' \delta_{\tilde{\bm z}},{\color{black}\lambda}}-{\bm f}_{P,{\color{black}\lambda}})/\epsilon.
$$
In the case that $c({\bm z},{\bm f}({\bm x}))$ is continuously differentiable in ${\bm f}$ for almost every ${\bm z}\in Z$ and ${\cal N}_{\cal F}$ is proto-differentiable,
we may represent ${\rm IF}(\tilde{{\bm z}};{\bm f}_{(\cdot,{\color{black}\lambda })},P)$
in terms of the derivatives of $c$ and 
the proto-derivative of ${\cal N}_{\cal F}$. 
The next proposition states this.

\begin{proposition}[Expression of IF of ${\bm f}_{\color{black}P,\lambda}$]\label{prop:influence-solution}
Assume:
(a) 
${\rm (C1)}$, 
${\rm (C2)}$, 
${\rm (C3)'}$, ${\rm(C4)}$ and ${\rm(C5)}$
hold;
(b)  $P\in \hat{\cal P}\cap \widetilde{\cal P}_{\bm h}$,
for all ${\bm h}\in {\cal H}_k$,
where $\hat{\cal P}$ and $\widetilde{\cal P}_{\bm h}$ are defined as in (\ref{eq:set_Q}) and (\ref{eq:tilde_P}) respectively;
(c)

normal cone mapping ${\cal N}_{\cal F}$ is proto-differentiable at ${\bm f}_{\color{black}P,\lambda}$ for ${\bm u}^*$,
where ${\bm u}^*=-\bbe_{P}[ c'_2({\bm z},{\bm f}_{P}({\bm x})){\color{black}k_{\bm x}}]$.
 If influence function ${\rm IF}(\tilde{{\bm z}};{\bm f}_{(\cdot,{\color{black}\lambda})},P)$ 
 is well-defined, then it has the following expression:
 \bgeq
{\rm IF}(\tilde{{\bm z}}; {\bm f}_{(\cdot,{\color{black}\lambda})},P)&=&
\left\{{\bm h}\in {\cal H}_k:
\begin{array}{ll}
&- c'_2 (\tilde{{\bm z}},{\bm f}_{P,{\color{black}\lambda}}({\tilde{\bm x}}))k_{\tilde{\bm x}}+\bbe_{P}[ c'_2 ({\bm z},{\bm f}_{P,{\color{black}\lambda}}({\bm x})){\color{black}k_{\bm x}}] \nonumber \\
&  - \bbe_{P,{\color{black}\lambda}}[ c''_2 ({\bm z},{\bm f}_{P,{\color{black}\lambda}}({\bm x}))T_{\bm x}{\bm h}] {\color{black}-2\lambda {\bm h}}\in ({\cal N}_{\cal F})'_{{\bm f}_{P,{\color{black}\lambda}},{\bm u}^*}({\bm h})
\end{array}
\right\}.
\edeq
Moreover,
if ${\cal F}={\cal H}_k$
and $\bbe_P[ c''_2 ({\bm z},{\bm f}({\bm x})){T_{\bm x}}]: {\cal H}_k \rightarrow {\cal H}_k$ is one-to-one and onto, then
\bgeqn
\label{eq:influence-equation}
 {\rm IF}(\tilde{{\bm z}};{\bm f}_{(\cdot,{\color{black}\lambda})},P)
 =-(\bbe_{P}[c''_2 ({\bm z},{\bm f}_{P,{\color{black}\lambda}}({\bm x})){T_{\bm x}}])^{-1} (c'_2(\tilde{{\bm z}},{\bm f}_{P,{\color{black}\lambda }}(\tilde{{\bm x}}))k_{\tilde{{\bm x}}}
  -\bbe_{P}[c'_2 ({\bm z},{\bm f}_{P,{\color{black}\lambda}}({\bm x})){\color{black}k_{\bm x}}]).
\edeqn
\end{proposition}

\noindent{\bf Proof}.
This is a special case of \cite[Theorem~4.1]{LeR94} with ${\bm z}=0$, $p=t$.
Under ${\rm (C4)}$,
$c'_2({\bm z}, w)$ is continuous in $w$.
Combining with the continuity of ${\bm f}({\bm x})$ w.r.t ${\bm f}$ for fixed ${\bm x}$,
we obtain that
$c'_2({\bm z},{\bm f}({\bm x})){\color{black}k_{\bm x}}$ is continuous in ${\bm f}$.
Then we have the following direction derivative
\bgeq
&& \Phi'_{\tilde{{\bm z}},P,{\color{black}\lambda }}({\bm f}_{\color{black}P,\lambda},0;{\bm h},1)\\
&=&\lim_{\epsilon \downarrow 0, (\tilde{{\bm h}},\tilde{t})\rightarrow ({\bm h},1)} [\Phi_{\tilde{{\bm z}},P,{\color{black}\lambda}}({\bm f}_{P,{\color{black}\lambda}}+\epsilon \tilde{{\bm h}}, 0+\epsilon\tilde{t}) -\Phi_{\tilde{{\bm z}},P,{\color{black}\lambda}}({\bm f}_{P,{\color{black}\lambda}},0)]/\epsilon\\
&=& \lim_{\epsilon \downarrow 0, (\tilde{{\bm h}},\tilde{t})\rightarrow ({\bm h},1)}
\epsilon \tilde{t}(c'_2(\tilde{{\bm z}},({\bm f}_{P,{\color{black}\lambda}}+\epsilon \tilde{{\bm h}})(\tilde{\bm x}))k_{\tilde{\bm x}}-\bbe_{P}[c'_2({\bm z},({\bm f}_{P,{\color{black}\lambda}}+\epsilon \tilde{{\bm h}})({\bm x})){\color{black}k_{\bm x}}])/\epsilon\\
&& \qquad \qquad \quad  +
\{\bbe_{P}[c'_2({\bm z},({\bm f}_{P,{\color{black}\lambda}}+\epsilon \tilde{{\bm h}})({\bm x})){\color{black}k_{\bm x}}]-\bbe_P[c'_2 ({\bm z},{\bm f}_{P,{\color{black}\lambda}}({\bm x})){\color{black}k_{\bm x}}]\}/\epsilon\\
&& \qquad \qquad \quad 
{\color{black} +\{2\lambda ({\bm f}_{P,{\color{black}\lambda}}+\epsilon \tilde{{\bm h}})-2\lambda {\bm f}_{P,{\color{black}\lambda}}\}/\epsilon
}\\
&=& c'_2 (\tilde{{\bm z}},{\bm f}_{P,{\color{black}\lambda}}(\tilde{{\bm x}}))k_{\tilde{{\bm x}}}-\bbe_{P}[c'_2 ({\bm z},{\bm f}_{P,{\color{black}\lambda}}({\bm x})){\color{black}k_{\bm x}}] + \bbe_{P}[c''_2 ({\bm z},{\bm f}_{P,{\color{black}\lambda}}({\bm x})){\color{black}k_{\bm x}}{\bm h}({\bm x})] 
{\color{black}+ 2\lambda {\bm h}}.
\edeq
{\color{black}By (\ref{eq:tx1}),
$T_{\bm x}{\bm h}={\color{black}k_{\bm x}}{\bm h}({\bm x})$.
}
Then it follows from \cite[Theorem~4.1]{LeR94} that 
\bgeq
\tilde{\bm f}'_{0;1}
&=&\left\{{\bm h}\in {\cal H}_k: 
\begin{array}{ll}
&-c'_2 (\tilde{{\bm z}},{\bm f}_{P,{\color{black}\lambda}}(\tilde{{\bm x}}))k_{\tilde{{\bm x}}}+\bbe_{P}[c'_2 ({\bm z},{\bm f}_{P,{\color{black}\lambda}}({\bm x})){\color{black}k_{\bm x}}] \\
& -\bbe_{P}[c''_2 ({\bm z},{\bm f}_{P,{\color{black}\lambda}}({\bm x}))T_{\bm x}{\bm h}] {\color{black}-2\lambda {\bm h} } \in ({\cal N}_{\cal F})'_{{\bm f}_{P,{\color{black}\lambda}},u^*}({\bm h})
\end{array}
\right\}.
\edeq
In the case that  ${\cal F}={\cal H}_k$,
${\cal N}_{\cal F}({\bm f})=\{0\}$.
Since $\bbe_P[c''_2 ({\bm z},{\bm f}({\bm x})){T_{\bm x}}]: {\cal H}_k \rightarrow {\cal H}_k$ is one-to-one and onto,
then
(\ref{eq:influence-equation}) holds 
by the classical implicit function theorem.  
\hfill
\Box

The next theorem states a sufficient condition for the boundedness of influence function ${\rm IF}(\tilde{{\bm z}};{\bm f}_{(\cdot,{\color{black}\lambda})},P)$.
Let $d({\bm f},S):=\inf_{{\bm f}'\in S} \|{\bm f}-{\bm f}'\|_k$ denote the distance from a function ${\bm f}\in {\cal H}_k$ to a set $S\subset {\cal H}_k$.
\begin{theorem}
[Boundedness of ${\rm IF}$]
\label{thm:bound_IF}
Assume the setting and conditions of Proposition~\ref{prop:influence-solution}.
 Suppose that the generalized equation 
 $$
 0\in \bbe_{P}[ c''_2 ({\bm z},{\bm f}_{P,{\color{black}\lambda}}({\bm x})){\color{black}k_{\bm x}} {\bm h}({\bm x})] {\color{black}+2\lambda {\cal I}}
 +({\cal N}_{\cal F})'_{{\bm f}_{P,{\color{black}\lambda}},{\bm u}^*}({\bm h})
 $$
 has a unique solution $0$,
where ${\cal I}$ is the identity operator from  ${\cal H}_k$ to $\mathcal{H}_k$.
 Then 
the following assertions hold.
\begin{itemize}
    \item[(i)]  
    ${\rm IF}(\tilde{{\bm z}};{\bm f}_{(\cdot,{\color{black}\lambda})},P)$ is bounded for every $\tilde{{\bm z}}\in Z$ provided that it is well-defined.
    \item[(ii)] 
    If, in addition,
$\Psi({\bm f}):=
\bbe_{P}[c'_2 ({\bm z},{\bm f}({\bm x})){\color{black}k_{\bm x}}]{\color{black}+2\lambda {\bm f}} +{\cal N}_{\cal F}({\bm f}) 
$
    is strongly metrically subregular at ${\bm f}_{P,{\color{black}\lambda}}$ for $0$ with regular modulus $\kappa'$, i.e.,
    {\color{black}
    there exist a constant $\kappa'>0$,
     neighborhoods $U_{{\bm f}_{P,{\color{black}\lambda}}}\subset {\cal H}_k$ of ${\bm f}_{P,{\color{black}\lambda}}$ and $U_{0}\subset {\cal H}_k$
     of $0$ such that
     $\|{\bm f}-{\bm f}_{\color{black}P,\lambda}\|_k\leq \kappa' d(0,\Psi({\bm f})\cap U_{0})$ for all ${\bm f}\in U_{{\bm f}_{\color{black}P,\lambda}}$,
    }
    then  {\color{black}there exists a constant $t_0>0$ such that}
\bgeqn
\label{eq:boundedness-IF}
    \sup_{\tilde{{\bm z}}\in Z} \|{\rm IF}(\tilde{{\bm z}};{\bm f}_{(\cdot,{\color{black}\lambda})},P)\|_k \leq \kappa' \sup_{\tilde{{\bm z}}\in Z} 
    \Upsilon(\tilde{{\bm z}}),
\edeqn
where
\bgeq 
 \Upsilon(\tilde{{\bm z}}) :=
\sup_{
t\in[0,t_0]
}\|\bbe_{P}[c'_2({\bm z},{\bm f}_{(1-t)P+t\delta_{\tilde{\bm z}},{\color{black}\lambda}}({\bm x})){\color{black}k_{\bm x}}]-c'_2(\tilde{{\bm z}},{\bm f}_{(1-t)P+t\delta_{\tilde{\bm z}},{\color{black}\lambda}}(\tilde{{\bm x}}))k_{\tilde{{\bm x}}}\|_k.
\edeq
\end{itemize}
\end{theorem}

\noindent{\bf Proof}.
Part (i).
 By a similar argument to that of Proposition~\ref{prop:influence-solution},
we have 
 $\tilde{\bm f}'_{0;0}
  =\{{\bm h}\in {\cal H}_k : 0\in \bbe_P[ c''_2 ({\bm z},{\bm f}_{\color{black}P,\lambda}({\bm x})){\color{black}k_{\bm x}}{\bm h}({\bm x})] {\color{black}+2\lambda {\cal I}} +({\cal N}_{\cal F})'_{{\bm f}_{P,{\color{black}\lambda}},{\bm u}^*}({\bm h})\}=
  \{0\}, 
 $
 which implies 
$D\tilde{\bm f}_t
(0|{\bm f}_{P,{\color{black}\lambda}}))(0)=\{0\}$ for each $\tilde{{\bm z}}\in \Z$. 
Then for each fixed $\tilde{{\bm z}}\in Z$, by \cite[Proposition~2.1]{KiR92},
 there exist constants $\tau>0$, $\kappa>0$ and $\beta\in (0,1)$ depending on $\tilde{{\bm z}}$ such that
 $
 \{\tilde{\bm f}_t\}\cap (\tilde{\bm f}_0+\tau {\cal B}) \in \tilde{\bm f}_0+\kappa t {\cal B}
$
for all $t\in (0,\beta)$,
{\color{black}where ${\cal B}$ is a closed unit ball in ${\cal H}_k$.
}
Recall that $\tilde{\bm f}_t={\bm f}_{(1-t)P+t\delta_{\tilde{\bm z}},{\color{black}\lambda}}$
and $\tilde{\bm f}_0={\bm f}_{\color{black}P,\lambda}$.
Moreover,
$
\{\frac{{\bm f}_{(1-\epsilon s)P+\epsilon s \delta_{\tilde{\bm z}},{\color{black}\lambda}}-{\bm f}_{\color{black}P,\lambda}}{\epsilon}\}\cap \frac{\tau}{\epsilon} {\cal B} \subset \kappa s {\cal B}
$
for $\epsilon s\in (0,\beta)$.
Thus
$
{\rm IF}(\tilde{{\bm z}};{\bm f}_{(\cdot,{\color{black}\lambda})},P)=
\limsup_{\epsilon \downarrow 0, s\rightarrow 1}
\frac{{\bm f}_{(1-\epsilon s)P+\epsilon s \delta_{\tilde{\bm z}},{\color{black}\lambda}}-{\bm f}_{P,{\color{black}\lambda}}}{\epsilon}
\in  \kappa {\cal B},
$
which implies 
$
\|{\rm IF}(\tilde{{\bm z}}; {\bm f}_{(\cdot,{\color{black}\lambda})}, P)\|_k\leq \kappa.
$

Part (ii).
Under the subregularity condition, we have
\bgeq
\|{\bm f}_{Q,{\color{black}\lambda}}-{\bm f}_{P,{\color{black}\lambda}}\|_k\leq \kappa' d(0,\Psi({\bm f}_{Q,{\color{black}\lambda}})\cap U_{0}), \;\forall\; {\bm f}_{Q,{\color{black}\lambda}}\in U_{{\bm f}_{P,{\color{black}\lambda}}} \inmat{ with }  Q\in \mathscr{P}(Z).
\edeq
On the other hand,
since ${\bm f}_Q$ satisfies that $0\in \bbe_{Q}[c'_2 ({\bm z},{\bm f}_{Q,{\color{black}\lambda}}({\bm x})){\color{black}k_{\bm x}}]{\color{black}+2\lambda {\bm f}} +{\cal N}_{\cal F}({\bm f}_{Q,{\color{black}\lambda}}) $,
then
\bgeq
\Delta_Q:= \bbe_{P}[c'_2 ({\bm z},{\bm f}_{Q,{\color{black}\lambda}}({\bm x})){\color{black}k_{\bm x}}]-\bbe_{Q}[c'_2 ({\bm z},{\bm f}_{Q,{\color{black}\lambda}}({\bm x})){\color{black}k_{\bm x}}]\in 
\Psi({\bm f}_{Q,{\color{black}\lambda}}).
\edeq
Let $Q\in \mathscr{P}(Z)$ be such that 
$\Delta_Q\in U_0$. 
Then $\Delta_Q \in \Psi({\bm f}_{Q,{\color{black}\lambda}})\cap U_{0}$ and hence
\bgeqn
\|{\bm f}_{Q,{\color{black}\lambda}}-{\bm f}_{P,{\color{black}\lambda}}\|_k\leq 
\kappa' d(0,\Delta_Q)
=
\kappa' \|\bbe_{P}[c'_2 ({\bm z},{\bm f}_{Q,{\color{black}\lambda}}({\bm x})){\color{black}k_{\bm x}}]-\bbe_{Q}[ c'_2 ({\bm z},{\bm f}_{Q,{\color{black}\lambda}}({\bm x})){\color{black}k_{\bm x}}]\|_k. \nonumber\\
 \label{eq:metric-subregularity}
\edeqn
Let $Q=(1-t)P+t \delta_{\tilde{{\bm z}}}$
and $t_0>0$ and $\tau>0$ be 
positive constants such that ${\bm f}_{(1-t)P+t\delta_{\tilde{\bm z}},{\color{black}\lambda}}
\in U_{{\bm f}_P}$ for all $t\in [0,t_0]$ and
 $U_0\subset \tau {\cal B}$.
Then by (\ref{eq:metric-subregularity}),
\bgeq
&& {\bm f}_{(1-t)P+t\delta_{\tilde{\bm z}},{\color{black}\lambda}}\\
&\in & 
{\bm f}_{P,{\color{black}\lambda}}+\kappa' \min\left\{t \sup_{
t\in[0,t_0]
}\|\bbe_{P}[c'_2({\bm z},{\bm f}_{(1-t)P+t\delta_{\tilde{\bm z}},{\color{black}\lambda}}({\bm x})){\color{black}k_{\bm x}}]-c'_2(\tilde{{\bm z}},{\bm f}_{(1-t)P+t\delta_{\tilde{\bm z}},{\color{black}\lambda}}(\tilde{{\bm x}}))k_{\tilde{{\bm x}}}\|_k,\tau\right\} {\cal B}
\edeq
for all $t\in [0,t_0]$, which implies
\bgeq
\frac{{\bm f}_{(1-\epsilon s)P+\epsilon s\delta_{\tilde{\bm z}},{\color{black}\lambda}}-{\bm f}_{P,{\color{black}\lambda}}}{\epsilon} 
\in \kappa' \min\left\{s 
\Upsilon(\tilde{{\bm z}}),
\frac{\tau}{\epsilon}\right\} {\cal B}
\edeq
for $\epsilon s\in [0,t_0]$. Thus
\bgeq
{\rm IF}(\tilde{{\bm z}}; {\bm f}_{(\cdot,{\color{black}\lambda})},P)
&=& \limsup_{\epsilon \downarrow 0, s\rightarrow 1} \; ({\bm f}_{(1-\epsilon s)P+\epsilon s\delta_{\tilde{\bm z}},{\color{black}\lambda}}-{\bm f}_{P,{\color{black}\lambda}})/\epsilon
\in \kappa' 
\Upsilon(\tilde{{\bm z}})
{\cal B}, 
\edeq
which implies
$
\|{\rm IF} (\tilde{{\bm z}};{\bm f}_{(\cdot,{\color{black}\lambda})},P)\|_k \leq \kappa' 
\Upsilon(\tilde{{\bm z}}).
$
By taking the supremum w.r.t. $\tilde{{\bm z}}\in Z$, we obtain (\ref{eq:boundedness-IF}).
The proof is complete. 
\hfill 
\Box

We make 
some comments about how the established theoretical results in Proposition~\ref{prop:influence-solution} and Theorem~\ref{thm:bound_IF}
may be applied to data-driven problems where
the true probability distribution $P$ is unknown and the sample contains an outlier $\tilde{{\bm z}}$.
Let $Z_N:=\{{\bm z}^1,\cdots,{\bm z}^{N-1},\tilde{{\bm z}}\}$ where the first $N-1$ samples are iid~and 
are not 
{\color{black} perturbed}
whereas $\tilde{{\bm z}}$ 
is an outlier.
Let 
$$
P_{N-1}:=\frac{1}{N-1}\sum_{i=1}^{N-1}\delta_{{\bm z}^i}\;
\inmat{and} \;
Q_N:=\frac{1}{N} \left( \sum_{i=1}^{N-1} \delta_{{\bm z}^i}+\delta_{\tilde{{\bm z}}}\right)=\left(1-\frac{1}{N}\right)
P_{N-1}
+\frac{1}{N} \delta_{\tilde{{\bm z}}}.
$$
First, we can 
compute the influence function ${\rm IF}(\tilde{{\bm z}};{\bm f}_{(\cdot)},P_{N-1})$ by Proposition~\ref{prop:influence-solution}.
Next, 
by Definition~\ref{def:IF},
we can 
estimate 
the difference between the 
solutions 
based on $P_{N-1}$
and $Q_N$,
\bgeqn
\label{eq:estimate_f_QN}
{\bm f}_{Q_N}-{\bm f}_{P_{N-1}} \approx  \frac{1}{N} {\rm IF}(\tilde{{\bm z}};{\bm f}_{(\cdot)},P_{N-1}),
\edeqn
when 
${\rm IF}(\tilde{{\bm z}};{\bm f}_{(\cdot,{\color{black}\lambda})},P_{N-1})$ is singleton.
In the case that the outlier is known,
(\ref{eq:estimate_f_QN}) 
may give us guidance on the difference between the confidence regions constructed with the two estimators. 
In the case that the outlier is unknown,
we cannot use (\ref{eq:estimate_f_QN}) to compute the influence function, rather we have to treat every data as a possible outlier and 
compute 
an upper bound of the true 
unknown influence function
by 
\bgeq
    \|{\rm IF}(\tilde{{\bm z}};{\bm f}_{(\cdot,{\color{black}\lambda})},P_{N-1})\|_k  \leq \kappa' \sup_{{{\bm z}}\in Z_N} 
    \sup_{t\in [0,1]}
    && \|\frac{1}{N-1}\sum_{{\bm z}^i\in Z_N\backslash {\bm z}}c'_2({\bm z}^i,{\bm f}_{(1-t)
     P_{Z_N\backslash {\bm z}}+t\delta_{\bm z},{\color{black}\lambda}}({\bm x}^i))k_{{\bm x}^i}
  \\
  && - c'_2 ({{\bm z}},{\bm f}_{(1-t)
  P_{Z_N\backslash {\bm z}}
  +t\delta_{\bm z},{\color{black}\lambda}}({{\bm x}}))k_{{\bm x}}\|_k,
\edeq 
where $ P_{Z_N\backslash {\bm z}}:=
\frac{1}{N-1}\sum_{{\bm z}^i\in Z_N\backslash \{{\bm z}\}}\delta_{{\bm z}^i}$.
We can 
then use the bound and 
 ${\bm f}_{Q_N,{\color{black}\lambda}}$
to obtain a conservative 
confidence region of
 ${\bm f}_{P_{N-1},{\color{black}\lambda}}$.
In practice, it will be difficult to compute the 
right-hand side
of the inequality because it is difficult to obtain ${\bm f}_{(1-t)
  P_{Z_N\backslash \{{\bm z}\}}
  +t\delta_{\bm z},{\color{black}\lambda}}$. It might be interesting to derive a bound for the right-hand side  of 
  the inequality which is independent of ${\bm f}$.
  We leave 
  all these for future research.
  
  \section{Numerical tests}
   \label{sec:numerical tests}
 
 We have undertaken some numerical tests 
on the established theoretical results in the previous section with two small academic examples.

\subsection{All data perturbation}
\label{sec:all_data_example}
\begin{figure}[!ht]
\begin{minipage}[t]{0.32\linewidth}
\centering
\includegraphics[width=1.7in]{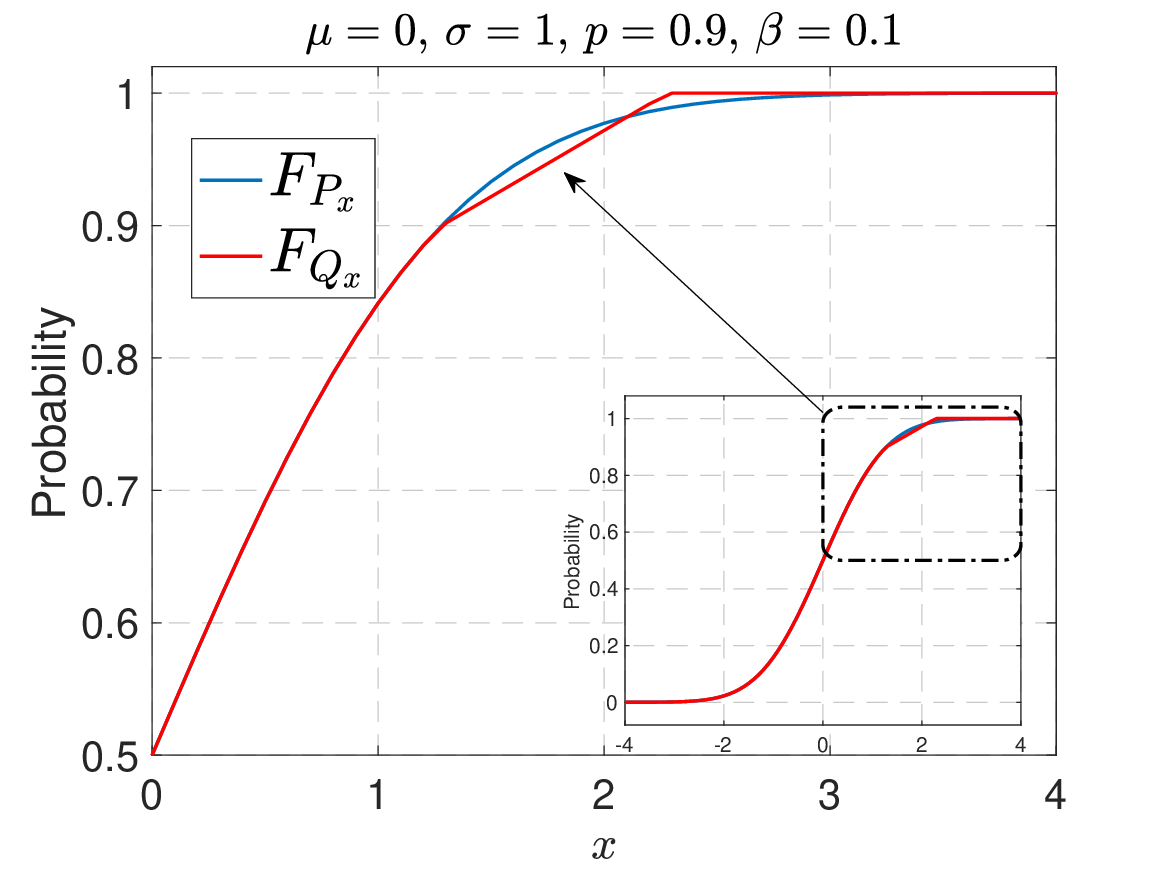}
\text{\tiny{(a) cdfs correspond to $P_{\bm x}$ and $Q_{\bm x}$.}}
\end{minipage}
\hspace{-1em}
\begin{minipage}[t]{0.3\linewidth}
\centering
\includegraphics[width=1.7in]{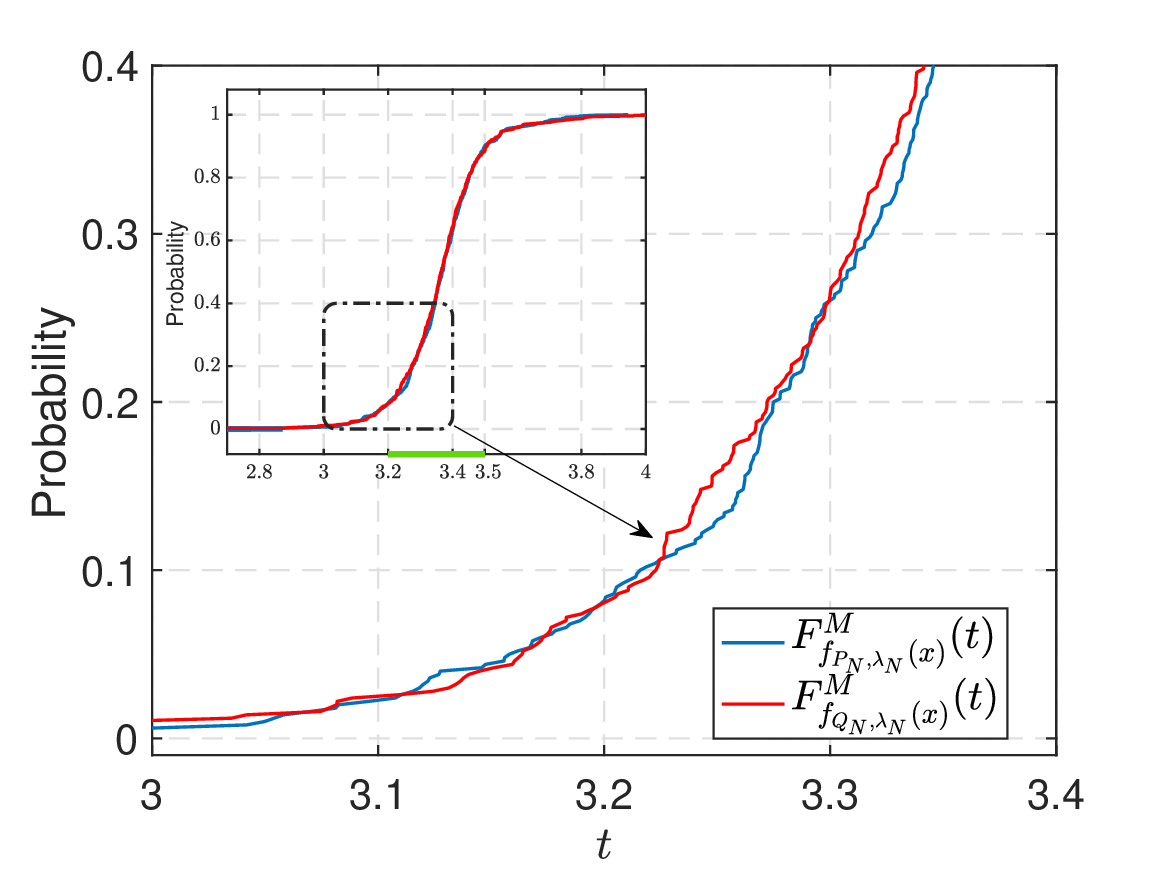}
\text{\tiny{(b) cdfs at $x=-1.9$.}}
\end{minipage}
\hspace{-0.5em}
 \begin{minipage}[t]{0.36\linewidth}
 \centering
\includegraphics[width=1.7in]{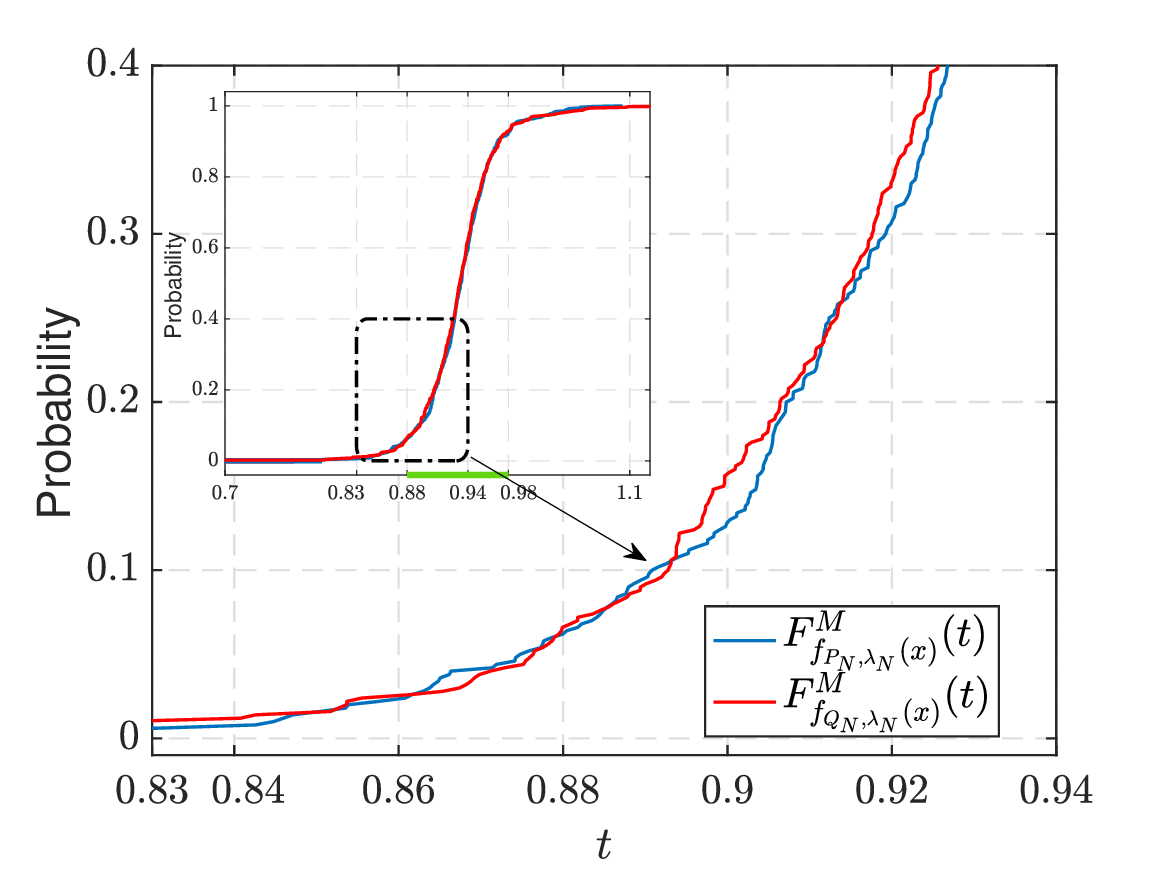}
\text{\tiny{(c) cdfs at $x=-1$.}}
\end{minipage}
 \captionsetup{font=footnotesize}
\caption{
(a) $F_{P_{\bm x}}$ and $F_{Q_{\bm x}}$.
(b)-(c): 
Empirical distributions of ${\bm f}_{P_N,\lambda_N}(x)$
and ${\bm f}_{Q_N,\lambda_N}(x)$
at points $x=-1.9$ and $x=-1$ with $M=500$. The embedded smaller graphs are plotted over larger ranges whereas the 
bigger graphs are enlargements of 
the curves over a specific range
where the gap between blue and red curves is significant.
The red curve
represents 
empirical distributions of ${\bm f}_{Q_N,\lambda_N}(x)$,
the 
blue curve
represents empirical distributions of ${\bm f}_{P_N,\lambda_N}(x)$.}
\label{fig:CDFs}
\end{figure}

Consider problem (\ref{eq:ML-Rf-min}), where
\bgeqn
\label{eq:cost_example}
c({\bm z},f(x)):=\frac{1}{2}(f(x)-y)^2,\;x\in \R,\;y\in \R,
\edeqn
 and kernel function
$
k(x_1,x_2):=\langle x_1,x_2\rangle^2.
$
Assume that
$x$
follows a normal distribution with mean value $\mu$ and standard deviation $\sigma$.
Let
$
F_{P_{\bm x}} 
$
be the cumulative distribution function (cdf) of $x$.
In this case, the regularized 
SAA
problem 
(\ref{eq:ML-saa-Q_N}) 
can be written as
$$
\min_{{\bm \alpha}\in {\color{black}[{\bm a},{\bm b}]}} \sum_{i=1}^N \frac{1}{N} \big(\tilde{y}^i- \sum_{j=1}^N \alpha_j k(\tilde{x}^j,\tilde{x}^i)\big)^2 +\lambda_N {\bm \alpha}^T \tilde{\bm K} {\bm \alpha},
$$
 where
{\color{black}
$[{\bm a},{\bm b}]:=[a_1,b_1]\times \cdots [a_N,b_N]$ with $a_i=-10$ and $b_i=10$, for $i=1,\cdots,N$,
}
${\bm \alpha}=(\alpha_1,\cdots,\alpha_N)^T$,
$\tilde{\bm K}=(\tilde{K}_{i,j})\in \R^{N\times N}$
with
$
\tilde{\bm K}_{i,j}= k(\tilde{x}^i,\tilde{x}^j), i,j=1,\cdots,N.
$
Here we write the feasible solution 
of (\ref{eq:ML-saa-Q_N})
as
$
\sum_{j=1}^N \alpha_j k(\tilde{x}^j,x)$.
The sample data $\tilde{x}^j$ and $\tilde{y}^j$, $j=1,\cdots,N$
are generated as follows. First, we consider 
a perturbation $Q_{\bm x}$ of $P_{\bm x}$
with  cdf:
\bgeq
F_{Q_{\bm x}}(x):=\left\{
\begin{array}{ll}
F_{P_{\bm x}}(x), & \inmat{for } x\leq x_0,\\
p+\beta(x-x_0),& \inmat{for } x_0\leq x\leq x_1,\\
1,& \inmat{for } x>x_1,
\end{array}
\right.
\edeq
where $x_0=F_{P_{\bm x}}^{-1}(p)$,
$x_1=x_0+\frac{1}{\beta}(1-p)$, $\beta$ and $p\in (0,1)$ are fixed positive constants.
Next, we use $Q_{\bm x}$ to generate samples $\tilde{x}^j, j=1,\cdots,N$ and 
$\tilde{y}^j=(\tilde{x}^j)^2$.
Let $\tilde{\bm z}^j=(\tilde{x}^j,\tilde{y}^j)$,
and $Q_N=\frac{1}{N}\sum_{j=1}^N \delta_{\tilde{\bm z}^j}$, where  $\delta_{\tilde{{\bm z}}}$ denotes 
the Dirac
distribution at a perturbed data point $\tilde{{\bm z}}$.
Let $P_N$ be defined in a similar way except 
that $x^j$, $j=1,\cdots,N$ are generated by 
$P_{\bm x}$.
The optimal solution
of  (\ref{eq:ML-saa-Q_N})
 can be written as
${\bm f}_{Q_N,\lambda_N}(x)=
\sum_{j=1}^N \alpha_j k(\tilde{x}^j,x)$.
The Kantorovich distance between 
$P_{\bm x}^{\otimes N}\circ ({\bm f}_{P_N,\lambda_N})^{-1}$
and $Q_{\bm x}^{\otimes N}\circ ({\bm f}_{Q_N,\lambda_N})^{-1}$
is
$\dd_{K}\left(P_{\bm x}^{\otimes N}\circ ({\bm f}_{P_N,\lambda_N})^{-1}, Q_{\bm x}^{\otimes N}\circ ({\bm f}_{Q_N,\lambda_N})^{-1}\right)
$.
For a fixed number $x\in \R$, let
$G_1$ and $G_2$ be 
the cdfs of ${\bm f}_{P_N,\lambda_N}(x)$ and ${\bm f}_{Q_N,\lambda_N}(x)$ respectively.
In this case,
\bgeq
\Delta_1(x):=\dd_{K}\left(P_{\bm x}^{\otimes N}\circ ({\bm f}_{P_N,\lambda_N}(x))^{-1}, Q_{\bm x}^{\otimes N}\circ ({\bm f}_{Q_N,\lambda_N}(x))^{-1}\right)
= \int_{-\infty}^{\infty} |G_1(t)-G_2(t)|dt.
\edeq
We use $Q_{\bm x}$ of $P_{\bm x}$
to generate 
$M$
 groups of samples each of which with size
$N=100$ and calculate
${\bm f}_{P_N,\lambda_N}(x)=
\sum_{j=1}^N \alpha_j^{P_N} k(\tilde{x}^j,x)$
and 
${\bm f}_{Q_N,\lambda_N}(x)=
\sum_{j=1}^N \alpha_j^{Q_N} k(\tilde{x}^j,x)$
for each group of samples.
For fixed $x$, we can then 
obtain the $M$ 
data points
and use to construct empirical
cdfs.
Figure~\ref{fig:CDFs} (a) depicts 
the difference between the cdfs of $P_{\bm x}$ and $Q_{\bm x}$.
Figure~\ref{fig:CDFs} (b) 
depicts
the cdfs (more precisely the cumulative frequency)  of ${\bm f}_{P_N,\lambda_N}(x)$ and ${\bm f}_{Q_N,\lambda_N}(x)$ at 
point $x=-1.9$.
We can see 
that most data are located within the range $[3.2,3.5]$
(the cdfs are very steep, see smaller embedded graph) 
and 
the 
red curve
approximates the blue
ones very well.
Figure~\ref{fig:CDFs} (c) 
displays similar phenomena 
when $x=-1$. In this case,
most data points fall within the range $[0.88, 0.98]$.

\begin{figure}[!ht]
\begin{minipage}[t]{0.32\linewidth}
\centering
\includegraphics[width=1.75in]{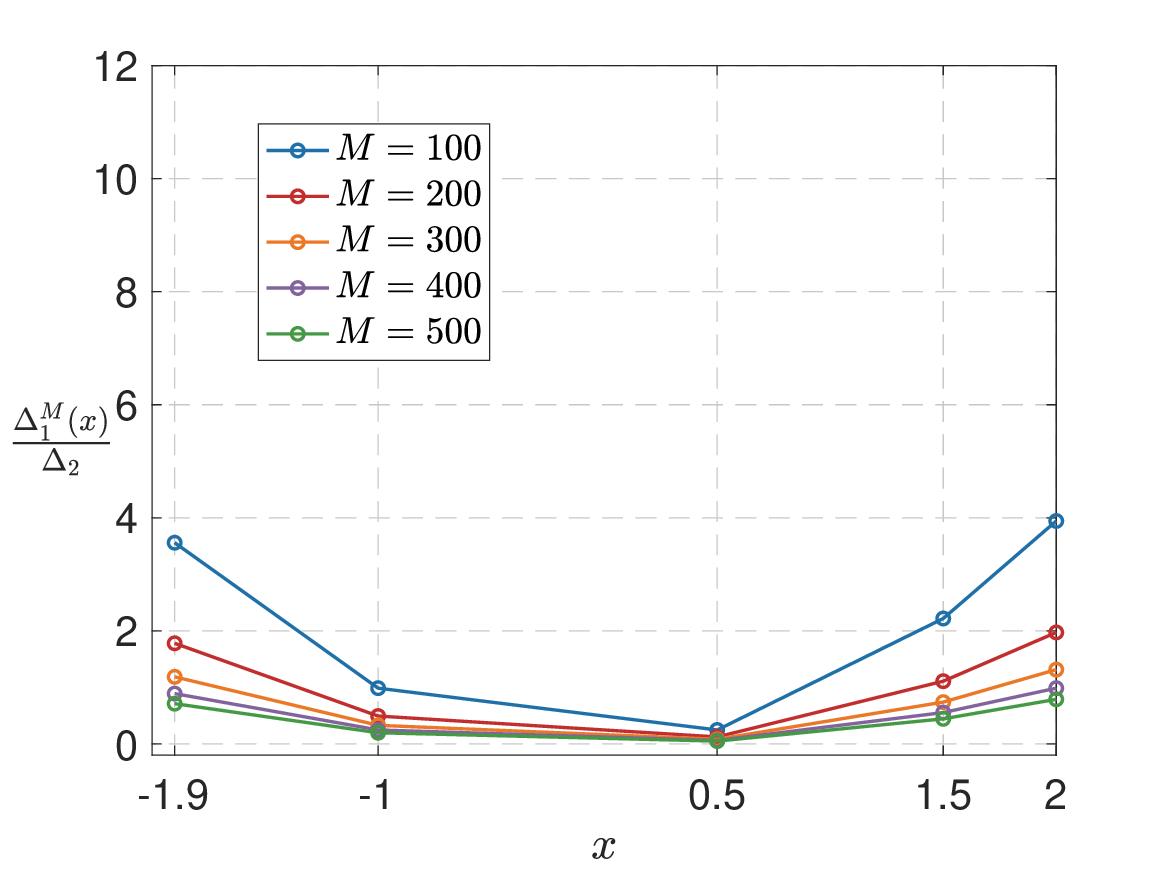}
\text{\tiny{(a) $\frac{\Delta_1^M(x)}{\Delta_2}$ at points $1$-$5$.}}
\end{minipage}
 \hspace{-1em}
\begin{minipage}[t]{0.32\linewidth}
\centering
\includegraphics[width=1.75in]{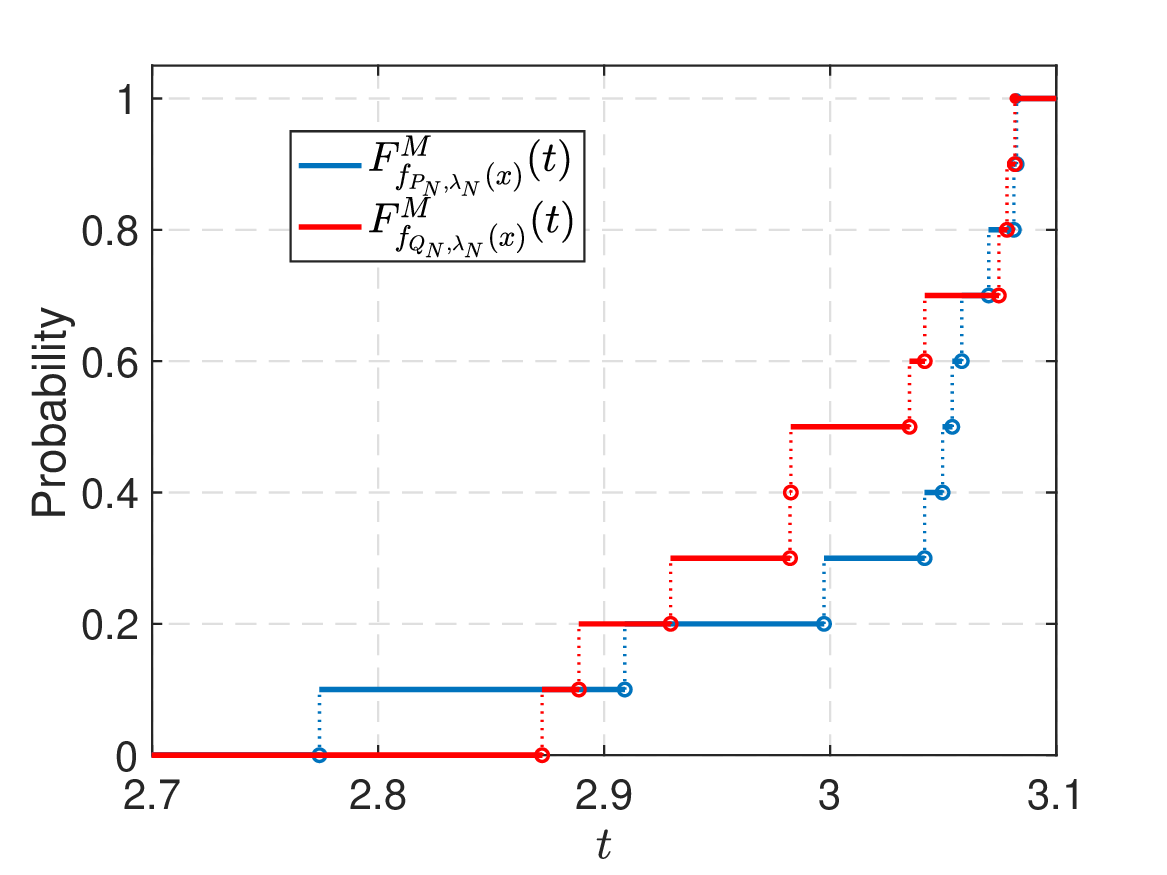}
\text{\tiny{(b) Simulation: $M=10$.}}
\end{minipage}
 \hspace{-1em}
\begin{minipage}[t]{0.36\linewidth}
\centering
\includegraphics[width=1.75in]{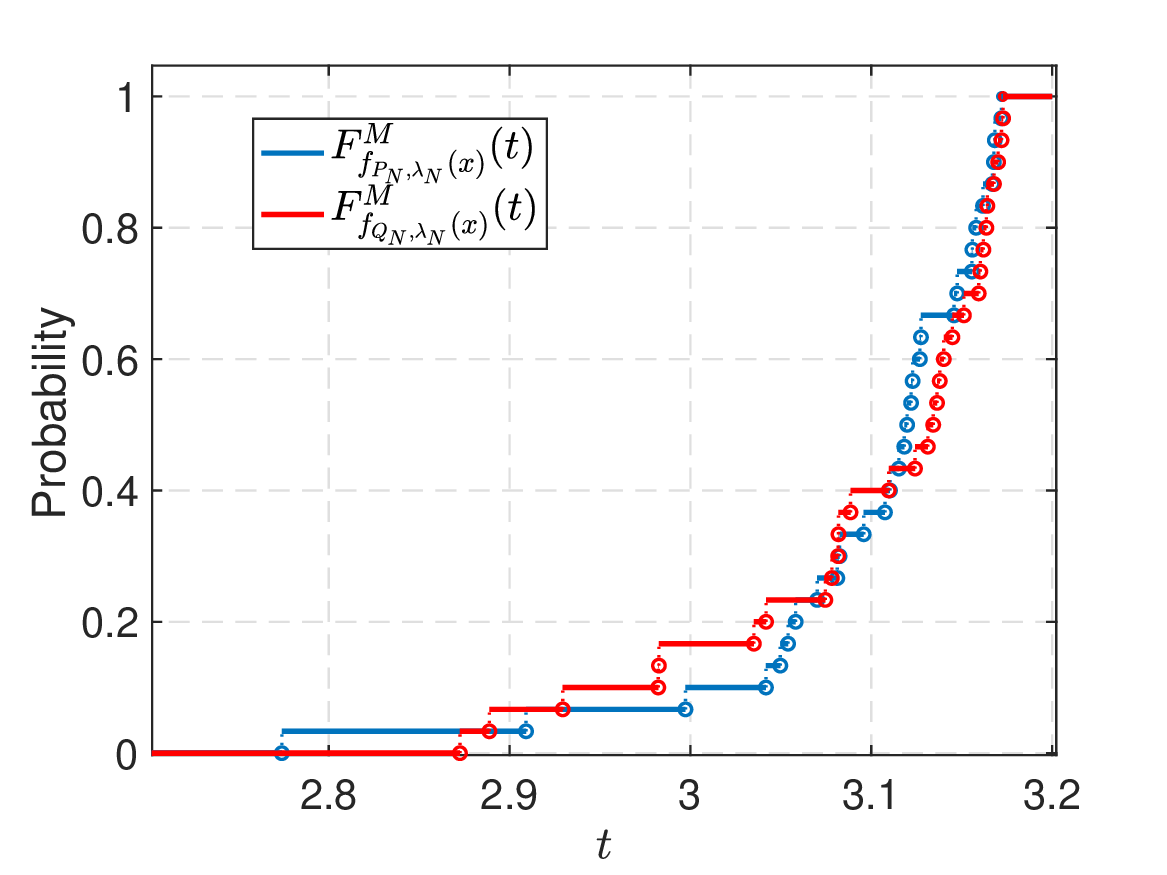}
\text{\tiny{(c) Simulation: $M=30$.}}
\end{minipage}
\hspace{-1em}
 \captionsetup{font=footnotesize}
\caption{
(a) $\frac{\Delta_1^M(x)}{\Delta_2}$ by 
$M$ simulations at five points when $M$ varies from $100$ to $500$.
(b)-(c): 
Empirical distributions of ${\bm f}_{P_N,\lambda_N}(x)$ and ${\bm f}_{Q_N,\lambda_N}(x)$
at point $x=-1.9$ 
with $M=10$ and $M=30$
simulations.
}
\label{fig:step_like}
\end{figure}

To examine the approximation of the optimal solutions more 
closely, we consider the Kantorovich distance 
between the cdfs of 
${\bm f}_{P_N,\lambda_N}(x)$
and ${\bm f}_{Q_N,\lambda_N}(x)$
relative to the 
 Kantorovich distance 
between the cdfs of the input data
generated by $P_{\bm x}$ and $Q_{\bm x}$. 
Let
\bgeq
\Delta_2&:=& \dd_{K}(P_{\bm x},Q_{\bm x})
\; \inmat{and} \;
\Delta_1^M(x):=\dd_{K}\left(F_{{\bm f}_{P_N,\lambda_N}(x)}^M,F_{{\bm f}_{Q_N,\lambda_N}(x)}^M\right),
\edeq
 where 
 {\color{black}
 $F_{{\bm f}_{P_N,\lambda_N}(x)}^M$ is the empirical distribution of random variable ${\bm f}_{P_N,\lambda_N}(x)$ using $M$ samples,
 that is,
 we first generate $N$ iid samples
 to compute ${\bm f}_{P_N,\lambda_N}(x)$, and then
 simulate $M$ times to obtain $M$ samples $\{{\bm f}_{P_N^m,\lambda_N}(x)\}_{m=1}^M$, 
 consequently we obtain $F_{{\bm f}_{P_N,\lambda_N}(x)}^M(y):=\frac{1}{M} \sum_{m=1}^M \mathds{1}_{y\geq {\bm f}_{P_N^m,\lambda_N}(x)}(y)$.
 }
 Figure~\ref{fig:step_like} (a)
 depicts the ratios 
 $\frac{\Delta_1^M(x)}{\Delta_2}$
 at five different points,
 we can see that as $M$ increases the ratios converge. 
   The ratio $\frac{\Delta_1^M(x_l)}{\Delta_2}$, $l=1,\cdots, 5$,
 is somehow related to 
 the constant $\widetilde{C}_{\kappa}$ in (\ref{eq:Quant-SR}).

\subsection{Single data perturbation} 
{\color{black}

We consider the same example as in the previous subsection with the cost function
(\ref{eq:cost_example}).
The true distribution $P$
of random vector ${\bm z}=(x,y)$ is defined  
as follows:
component $x$ follows a normal distribution with the mean value 
$\mu$ and standard deviation $\sigma$
whereas component $y$ is equal to
$x^2+\epsilon$,
where 
$\epsilon$ follows a
normal distribution with mean $0$ and variance $0.01$. 
We follow the standard procedures (see e.g.~\cite{Pen16})
to generate $N$ samples ${\bm z}^j_P=(x^j,y^j)$, $j=1,\cdots,N$ of ${\bm z}=(x,y)$ 
with $P$.
To generate $N$ samples 
 ${\bm z}^j_{(1-t)P+t\delta_{\tilde{\bm z}}}=(x^j_{(1-t)P+t\delta_{\tilde{\bm z}}},y^j_{(1-t)P+t\delta_{\tilde{\bm z}}})$, $j=1,\cdots,N$ of ${\bm z}=(x,y)$ with the mixture distribution
 $(1-t)P+t\delta_{\tilde{\bm z}}$,
where $\tilde{\bm z}=(\tilde{x}, \tilde{y})$ with $\tilde{y}=\tilde{x}^3$ is an outlier, we
begin by generating $N$ samples with distribution $(1-t)P+t\delta_{\tilde{\bm z}}$ 
using the Dirac distribution
$\delta_{\tilde{x}}$ at $\tilde{x}$ and
a switching variable $X$,
which takes the value $X=1$ with probability $1-t$ and $X=2$ with probability $t$.
Here $X$ is
independent of $P$ and $\delta_{\tilde{z}}$.
Specifically, we generate 
$$
{z}^j_{(1-t)P+t\delta_{\tilde{\bm z}}}
:=\left\{
\begin{array}{ll}
{\bm z}_P^j:=(x^j_P, y^j_P) &  \inmat{if }X^j=1,\\
\tilde{\bm z} &  \inmat{if }X^j=2,
\end{array}
\right. \; \inmat{for} \; j=1,\cdots,N,
$$
where 
$x_P^j$,
$j=1,\cdots,N$, 
are the iid samples generated by $P_x$.
The empirical distributions can be respectively written as
$$
P_{N}:=\frac{1}{N}\sum_{j=1}^N \delta_{{\bm z}^j_{P}}
\quad
\inmat{and}
\quad
((1-t)P+t\delta_{\tilde{\bm z}})_N:=\frac{1-t}{N-1} \sum_{j=1}^{N-1}\delta_{{\bm z}^j_P}+t\delta_{\tilde{\bm z}}.
$$
In this case, 
the optimal solutions of 
the regularized SAA problem (\ref{eq:ML-saa-Q_N}) 
with probability distributions $P$ and $(1-t)P+t\delta_{\tilde{\bm z}})$
take the form of 
$$
{\bm f}_{P_N,{\color{black}\lambda}}(x)=\sum_{j=1}^N \alpha_j k(x^j_{P},x)
\quad
\inmat{and}
\quad 
{\bm f}_{((1-t)P+t\delta_{\tilde{\bm z}})_N,{\color{black}\lambda}}(x)=
\sum_{j=1}^N \alpha_j k\left(x^j_{(1-t)P+t\delta_{\tilde{\bm z}}},x\right).
$$
The regularized SAA problem (\ref{eq:ML-saa-Q_N}) and its single data perturbation problem  can be solved via the following two problems
\bgeq
\min_{{\bm \alpha}\in \R^N} \sum_{i=1}^N \frac{1}{N} \Big({y}^i- \sum_{j=1}^N \alpha_j k({x}^j,{x}^i)\Big)^2 +{\color{black}\lambda} {\bm \alpha}^T {\bm K}_P {\bm \alpha}
\edeq
and 
\bgeq
\min_{{\bm \alpha}\in \R^{N}} \sum_{i=1}^{N-1} \frac{1-t}{N-1} \Big({y}^i- \sum_{j=1}^{N} \alpha_j k({x}^j,{x}^i)\Big)^2
+t\Big(y^{N}-\sum_{j=1}^{N} \alpha_j k({x}^j,{x}^{N})\Big)^2+{\color{black}\lambda} {\bm \alpha}^T {\bm K}_{(1-t)P+t\delta_{\tilde{\bm z}}} {\bm \alpha}, 
\edeq
where ${\bm K}_P$
and ${\bm K}_{(1-t)P+t\delta_{\tilde{\bm z}}} $ are the Gramer matrices with
$(K_P)_{ij}:=k(x^i_P,x^j_P)$,
and $(K_{(1-t)P+t\delta_{\tilde{\bm z}}})_{ij}:=k(x^i_{(1-t)P+t\delta_{\tilde{\bm z}}},x^j_{(1-t)P+t\delta_{\tilde{\bm z}}})$, for $i,j=1,\cdots,N$.

Figure~\ref{fig:IF} depicts the change of the influence function as the outlier varies. We can see that as $\tilde{x}$ increases (hence $\tilde{y}$ also changes accordingly), the influence function increases, which means that 
when the outlier shifts away from the normal sample data, it has a greater effect on the sensitivity of the kernel learning estimator.
Moreover, we can see that as $\lambda_N$ increases, the influence function value 
decreases, which is consistent with the fact that a larger $\lambda_N$ enhances the stability of the kernel learning estimator and hence makes it less sensitive to the outlier.

\begin{figure}[!ht]
\begin{minipage}[t]{\linewidth}
\centering
\includegraphics[width=2.5in]{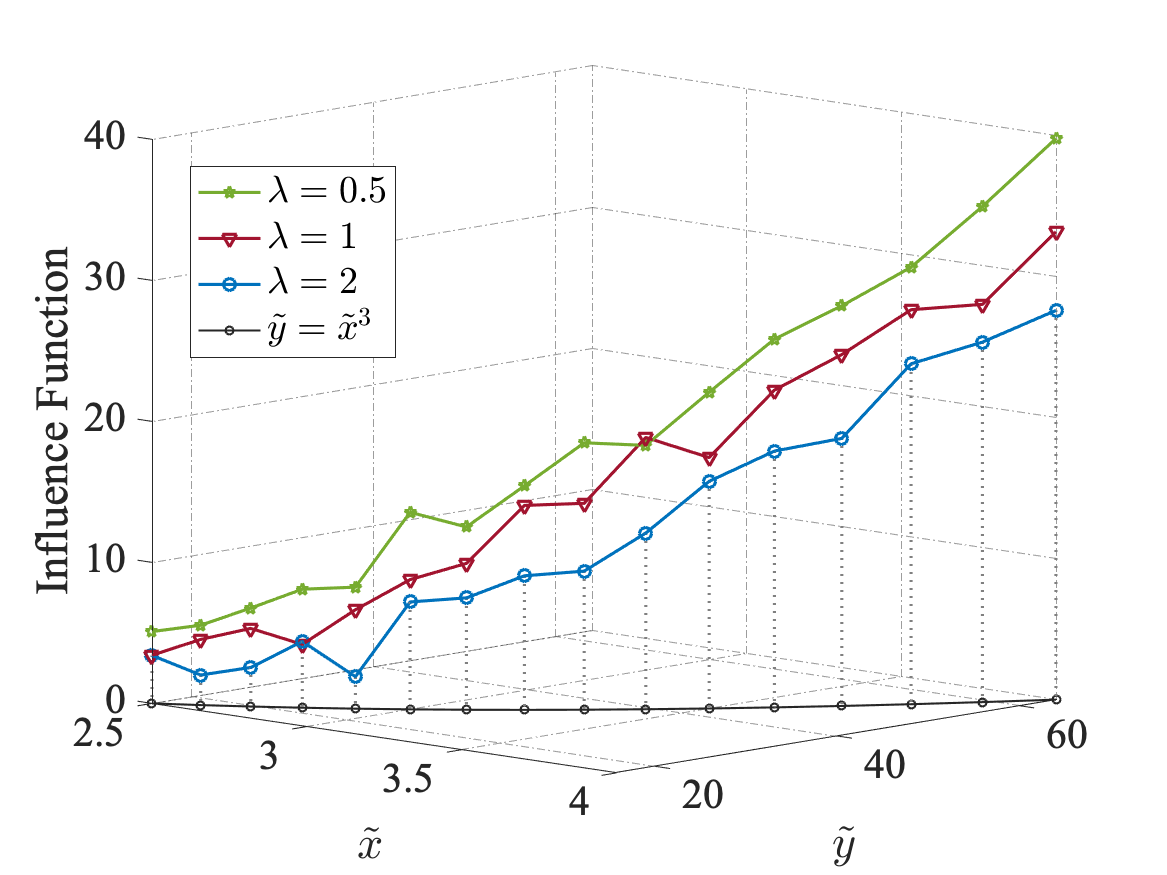}
\end{minipage}
 \caption{Performance of the influence function ${\rm IF}(\tilde{\bm z};{\bm f}_{(\cdot,{\color{black}\lambda})}(x_0),P)$ at point $x_0=0.5$ when
 outlier $\tilde{\bm z}$ satisfies $\tilde{y}=(\tilde{x})^3$, and
 $\mu=0$, $\sigma=1$.
 }
 \label{fig:IF}
 \end{figure}

  \section{Concluding remarks}
  \label{sec:conclusion}
 In this paper, we investigate 
 quality of 
 the kernel learning estimator 
 (the optimal solution to the empirical risk minimization problem) obtained with perturbed
 data. In single data perturbation case, 
 we use the notion of influence function
 to measure 
 the sensitivity of the kernel learning estimator w.r.t. 
 the perturbation of the data; in the case when all data are potentially 
 {\color{black} perturbed}, we derive conditions under which the learning 
 kernel learning estimator  based on 
 the {\color{black} perturbed}
 data is close to the one that is based on real data (without 
 {\color{black} perturbation}
 ) in the sense that the distributions of the two learning estimators are linearly bounded by the difference of the true probability distributions generating the training data under the Kantorovich metric.
This is a step forward from the traditional stability analysis which measures the difference of two statistical estimators based on each set of sample data.
 While the theoretical results are useful in some practical applications as outlined in Sections~\ref{sec:Quali-SR-KLE} and \ref{sec:Quant-SR-KLE}
  and the numerical tests confirm the theoretical results in the small academic examples, it would be more interesting 
 to apply the theoretical results in substantial practical cases. For instances, 
 Kern et al.~\cite{KSZ20} apply the influence function approach to a Markov decision-making process when the transition probability deviates from the true one. Besbes et al.~\cite{besbes2022beyond} propose distributionally robust optimization models 
in data-driven decision-making
 such as newsvendor, pricing, and ski rental under heterogeneous environments.
 It 
 might be interesting to apply our theoretical results to these circumstances.
 {\color{black} 
 It will also be interesting to consider
 an approach similar to the In-CVaR of Liu and Pang \cite{LiP23} to reduce the effect of data perturbation when all data are potentially contaminated.}
 We leave all these for future research.

\section*{Acknowledgments.}
{\color{black}
This project is supported by RGC grant 
14204821.
The third author is supported by National Natural Science Foundation of China (Grant Nos. 12122108).}

\begin{appendix}
\section{Berge's maximum theorem}
\label{Sec:Berge-max}
\begin{theorem}[Berge's maximum theorem,  Page 116 \cite{berge1877topological}]
\label{T-berge-max}
Let $\phi$ be a continuous numerical function in a topological space $X\times Y$
and $\Gamma$ be a continuous mapping of  $X$ into $Y$ such that
for each $x$, $\Gamma(x) \neq \emptyset$. Then the numerical function $M$ defined by $M(x):=\max\{\phi(y): y\in \Gamma(x)\}$
is continuous in $X$ and the mapping $\Phi$ defined by
$\Phi(x):=\{y: y\in \Gamma(x),\phi(y)=M(x)\}$ is a upper semicontinuous mapping
of $X$ into $Y$.
\end{theorem}

\section{Differentiability of a function defined over \texorpdfstring{${\cal H}_k$}{}}
\label{Def:Gateaux-diffentibility}

Recall that a functional
$g:{\cal H}_k\rightarrow \R$ 
is said to be 
directionally differentiable at 
${\bm f}\in {\cal H}_k$ 
if the limit
$$
g'({\bm f};{\bm h}):=\lim_{t\downarrow 0} \frac{g({\bm f}+t{\bm h})-g({\bm f})}{t}
$$
exists for all direction ${\bm h}\in {\cal H}_k$.
If, in addition, the
directional derivative 
$g'({\bm f};{\bm h})$ is linear and continuous in ${\bm h}$, then
$g$ is said to be 
G\^{a}teaux differentiable at ${\bm f}$ with derivative $Dg({\bm f})$ satisfying $$
g'({\bm f};{\bm h})=Dg({\bm f}){\bm h}.
$$
If the derivative 
$Dg(\cdot
):{\cal H}_k\to {\cal H}_k$ is continuous on an open set $S\subset {\cal H}_k$ in terms of the norm $\|\cdot\|_{k}$, 
then $g$ is said to be continuously differentiable on $S$.

For an operator $
{\color{black} {\cal A}}
:{\cal H}_k\to {\cal H}_k$,
it is said to be directionally 
differentiable at 
${\bm f}\in {\cal H}_k$ 
if the limit
$$
{\color{black} {\cal A}}'({\bm f};{\bm h}):=\lim_{t\downarrow 0} \frac{{\color{black} {\cal A}}({\bm f}+t{\bm h})-{\color{black} {\cal A}}({\bm f})}{t}
$$
exists for 
any fixed
direction ${\bm h}\in {\cal H}_k$.
If, in addition, the
directional derivative 
${\color{black} {\cal A}}'({\bm f};{\bm h})$ is linear and continuous in ${\bm h}$, then
${\color{black} {\cal A}}$ is said to be 
G\^{a}teaux differentiable at ${\bm f}$ with derivative $D {\color{black} {\cal A}}({\bm f})$ satisfying 
$$ {\color{black} {\cal A}}'({\bm f};{\bm h})=D {\color{black} {\cal A}}({\bm f}){\bm h}.
$$
Let ${\cal L}({\cal H}_k)$ be the space of linear continuous operators ${\color{black} {\cal A}}:{\cal H}_k\to {\cal H}_k$
equipped with the operator norm  
\bgeqn
\label{eq:defi-Lnorm}
\|{\color{black} {\cal A}}\|_{{\cal L}}:=\sup_{{\bm f}\in {\cal B}}\|{\color{black} {\cal A}}({\bm f})\|_k,
\edeqn
where ${\cal B}$
is the 
closed unit ball in ${\cal H}_k$.
If the derivative 
$D {\color{black} {\cal A}}(\cdot):{\cal H}_k\to {\cal L}({\cal H}_k)$ is continuous on an open set $S\subset {\cal H}_k$ in terms of the norm $\|\cdot\|_{{\cal L}}$,
then
{\color{black}
${\color{black} {\cal A}}$
}
is  said to be continuously differentiable on $S$, see \cite[pages 35-36]{BoS00}.
In what follows, 
we 
consider the case that 
$g({\bm f}) := c({\bm z},{\bm f}({\bm x}))$ for some fixed ${\bm z}$,
where $c$ is differentiable in the second argument,
and $ {\color{black} {\cal A}}({\bm f}) :=Dg({\bm f})$.
We will use $D_{\bm f}({\color{black} {\cal A}}({\bm f}))$ to denote the derivative of functional ${\color{black} {\cal A}}$ w.r.t.~${\bm f}$.

\section{Proof of Proposition~\ref{prop:metricize}}

The proof is inspired by that of 
\cite[Lemma~9]{KSZ20}\footnote{The second author would like to thank Henryk Zh\"ale for referring him to the lemma.}.
By \cite[Corollary A.48]{FoS16},
the $\psi$-weak topology is metrizable.
Thus the notions of continuity and sequential continuity coincide and it suffices
to prove that for any $\{\mu_N\}_{N=1}^{\infty}\subset {\cal M}_Z^{\psi}$ and $\mu\in {\cal M}_Z^{\psi}$,
$$
{\mu_N} \xrightarrow[]{\tau_\psi} \mu
\Longleftrightarrow \dd_{\rm FM}(\mu_N,\mu)\to 0.
$$
Note that 
${\mu_N} \xrightarrow[]{\tau_\psi} \mu$
if and only if 
${\mu_N} \xrightarrow[]{w} \mu$
and $\int_{Z} \psi({\bm z}) d \mu_N({\bm z})\to \int_Z \psi({\bm z}) d \mu({\bm z})$.
We proceed the rest of the proof in two steps.

\underline{Step 1.}
Suppose $\dd_{\rm FM}(\mu_N,\mu)\to 0$.
First, we show that 
${\mu_N} \xrightarrow[]{w} \mu$.
By the Portmanteau Theorem (\cite[Theorem~13.16]{Kle13}), 
it suffices to show that 
for any bounded {\em Lipschitz continuous} function $h:Z\to \R$,
\bgeqn 
\lim_{N\to\infty} \int_Z h({\bm z}) d\mu_N({\bm z}) = \int_Z h({\bm z}) d\mu({\bm z}).
\label{eq:Weak-convg-h} 
\edeqn 
Let $L_h$ be the Lipschitz modulus of $h$ and
$\tilde{h} :=h/L_h$. Then $\tilde{h}\in {\cal F}_{p}$. Thus $\dd_{\rm FM}(\mu_N,\mu)\to 0$ implies
$\lim_{N\to\infty} \int_Z \tilde{h}({\bm z}) d\mu_N({\bm z}) = \int_Z \tilde{h}({\bm z}) d\mu({\bm z}),$
and hence (\ref{eq:Weak-convg-h}) as desired. 
Next, we 
show 
$\int_Z \psi({\bm z}) d\mu_N({\bm z})\to \int_Z \psi({\bm z}) d \mu({\bm z})$.
It 
is enough to verify that there exists $L>0$ such that 
$\psi/L\in {\cal F}_{p}$. Observe that
\bgeq
|\psi({\bm z}_1)-\psi({\bm z}_2)|
&=& |1+\max\{\|{\bm z}_1\|,\|{\bm z}'\|\}^{p-1}\|{\bm z}_1-{\bm z}'\|-(1+\max\{\|{\bm z}_2\|,\|{\bm z}'\|\}^{p-1}\|{\bm z}_2-{\bm z}'\|)|\\
& \leq & \max\{1,\|{\bm z}'\|\}\max\{1,\|{\bm z}_1\|,\|{\bm z}_2\|\}^{p-1} \|{\bm z}_1-{\bm z}_2\|.
\edeq
Let $L:=\max\{1,\|{\bm z}'\|\}$.
Then $\psi/L\in {\cal F}_{p}$.
By the condition that $\dd_{\rm FM}(\mu_N,\mu)\to 0$, 
$\int_Z \psi({\bm z}) d\mu_N({\bm z})\to \int_Z \psi({\bm z}) d \mu({\bm z})$.

\underline{Step 2.} Conversely,
suppose ${\mu_N} \xrightarrow[]{\tau_\psi} \mu$. 
We need to show that $\dd_{\rm FM}(\mu_N,\mu)\to 0$,
that is,
for any $\epsilon>0$,
there exists some $N_0\in \mathbb{N}$ such that 
\bgeqn
\label{eq:d_FM0}
\sup_{h\in {\cal F}_{p}}\left|\int_Z h({\bm z}) d\mu_N({\bm z})-\int_Z h({\bm z}) d\mu({\bm z})\right|\leq \epsilon,\quad \forall N\geq N_0.
\edeqn
Let $K$ be a sufficiently large positive constant and
$
h_K:=h \mathds{1}_{|h|\leq K}+K\mathds{1}_{h>K}-K\mathds{1}_{h<-K},
$
and $h^K:=h-h_K$.
Then 
\bgeq
 \sup_{h\in {\cal F}_{p}} \left|\int_Z h({\bm z}) d\mu_N({\bm z})-\int_Z h({\bm z}) d\mu({\bm z})\right|
\leq R_1+R_2,
\edeq
where 
\bgeq
R_1 &:=&  
\sup_{h\in {\cal F}_{p}} \left|
\int_Z h_K({\bm z})d\mu_N({\bm z})-\int_Z h_K({\bm z}) d\mu({\bm z})
\right|,\\
R_2&:=&  \sup_{h\in {\cal F}_{p}} \left|\int_Z h^K({\bm z}) d\mu_N({\bm z})-\int_Z h^K({\bm z}) d\mu({\bm z})\right|.
\edeq
It suffices to show that 
$R_1,R_2\to 0$ as $N\to\infty$.
Since the set ${\cal F}_{p}$ has a one-to-one 
correspondence relationship with
\bgeq
\tilde{\cal F}_{p}:=\{h:Z\to \R: h({\bm z}')=0,
|h({\bm z}_1)-h({\bm z}_2)|\leq L_{p}({\bm z}
_1,{\bm z}_2)\|{\bm z}_1-{\bm z}_2\|,\forall {\bm z}_1,{\bm z}_2\in Z
\}, \nonumber 
\edeq
then we may consider $\tilde{\cal F}_{p}$. 
Consequently, we have 
\bgeqn
\label{eq:h-psi}
|h({\bm z})|=|h({\bm z})-h({\bm z}')|
&\leq & \max\{1,\|{\bm z}\|,\|{\bm z}'\|\}^{p-1} \|{\bm z}-{\bm z}'\|
<\psi({\bm z}).
\edeqn
In this case,
\bgeq
|h^K({\bm z})|
&=&|h({\bm z})-(h({\bm z})\mathds{1}_{|h({\bm z})|\leq K}+K\mathds{1}_{h({\bm z})>K}-K\mathds{1}_{h({\bm z})<-K})|\\
&=& |h({\bm z})\mathds{1}_{h({\bm z})>K} +h({\bm z})\mathds{1}_{h({\bm z})<-K}-(K\mathds{1}_{h({\bm z})>K}-K\mathds{1}_{h({\bm z})<-K})|\\
&=& |(h({\bm z})-K)\mathds{1}_{h({\bm z})>K} +(h({\bm z})+K)\mathds{1}_{h({\bm z})<-K}|\\
&\leq& |h({\bm z})| \mathds{1}_{|h({\bm z})|>K}
\leq \psi({\bm z})\mathds{1}_{|\psi({\bm z})| > K},
\edeq
where the last inequality  holds because (\ref{eq:h-psi}) implies $\mathds{1}_{|h({\bm z})|>K} \leq \mathds{1}_{|\psi({\bm z})|>K}$.
Thus
\bgeq
R_2 &\leq &\sup_{h\in {\cal F}_{p}} \left|\int_Z h^K({\bm z}) d\mu_N({\bm z})\right|+
\sup_{h\in {\cal F}_{p}} \left|\int_Z h^K({\bm z}) d\mu({\bm z})\right|\\
&\leq& \int_Z \psi({\bm z})\mathds{1}_{|\psi({\bm z})|> K}d \mu_N({\bm z})+\int_Z \psi({\bm z}) \mathds{1}_{|\psi({\bm z})|> K}d\mu({\bm z})=:R_{21}+R_{22}.
\edeq
Since 
 ${\mu_N} \xrightarrow[]{\tau_\psi} \mu$, 
then for any $\epsilon$, 
we can set $K$ sufficiently large such that
$R_{22}=\int_Z \psi({\bm z})\mathds{1}_{|\psi({\bm z})|> K} d\mu({\bm z})\leq \frac{\epsilon}{5}$.
Moreover,
\begin{eqnarray}
\label{eq:bounded_R21}
 R_{21}&\leq& \left|\int_Z \psi({\bm z}) \mathds{1}_{|\psi({\bm z})|> K}d \mu_N({\bm z})-\int_Z \psi({\bm z})\mathds{1}_{|\psi({\bm z})|> K} d\mu({\bm z})\right|
  + \int_Z \psi({\bm z}) \mathds{1}_{|\psi({\bm z})|> K} d\mu({\bm z}) 
 \nonumber \\
&\leq &  \left|\int_Z \psi({\bm z}) \mathds{1}_{|\psi({\bm z})|> K} d \mu_N({\bm z})-\int_Z \psi({\bm z})\mathds{1}_{|\psi({\bm z})|> K} d\mu({\bm z})\right| 
 +\frac{\epsilon}{5}\nonumber \\
&\leq &\left|\int_Z \psi({\bm z})d \mu_N({\bm z})-\int_Z \psi({\bm z}) d\mu({\bm z})\right| \nonumber \\
&& +\left|\int_Z \psi({\bm z})\mathds{1}_{|\psi({\bm z})|\leq  K}  d \mu_N({\bm z})-\int_Z \psi({\bm z}z) \mathds{1}_{|\psi({\bm z})|\leq  K} d\mu({\bm z})\right| +\frac{\epsilon}{5}.
\end{eqnarray}
Since 
${\mu_N} \xrightarrow[]{\tau_\psi} \mu$,
the first term in the last inequality of (\ref{eq:bounded_R21}) converges to $0$ as $N\to \infty$.
Thus we can find $N_0\in \mathbb{N}$ such that it is bounded above by $\frac{\epsilon}{5}$ for $N\geq N_0$.
Since $\mu \circ \psi^{-1}$ as a probability measure on the real line has at most countably many atoms,
we choose $K>0$ such that $\mu({\bm z}\in Z:\psi({\bm z})=K)=0$.
Since $\mu_N\to \mu$ ($\psi$-weakly and thus) weakly,
it follows by the Portmanteau theorem that the second term in the last inequality of (\ref{eq:bounded_R21}) converges to $0$ as $N\to \infty$.
By possibly increasing $N_0$ we obtain that the second term in the last inequality of (\ref{eq:bounded_R21}) is at most $\epsilon/5$ for all $N\geq N_0$.
So far we have shown that 
$R_2$ is bounded above by $\frac{4\epsilon}{5}$ for $N\geq N_0$.
Since the functions in $\{h_K:h\in {\cal F}_{p}\}$ are uniformly bounded and equicontinuous,
\cite[Corollary 11.3.4]{Dud02} ensures that one can increase $N_0$ further such that $R_1\leq \frac{\epsilon}{5}$ for all $N\geq N_0$.
Therefore,
we have (\ref{eq:d_FM0}) holds.
\hfill
\Box

\section{Proof of Proposition~\ref{prop:secon-order-d}}
\label{sec:proof_2d}

By Assumption~\ref{assu-solution-P}{\color{black}${\rm (C5)}$ and ${\rm (K1)}'$},
for every
${\bm z}\in Z$, 
\bgeq
&&|c''_2 ({\bm z},{\bm f}({\bm x}))|\\
&\leq & \Big(|c''_2 ({\bm z}_0,{\bm f}({\bm x}_0))|+C_2L_{p_2}({\bm z},{\bm z}_0)(\|{\bm z}-{\bm z}_0\|+|{\bm f}({\bm x})-{\bm f}({\bm x}_0)|)\Big)\\
&\overset{(\ref{eq:z_f-difference})}{\leq }&
 \Big(|c''_2 ({\bm z}_0,{\bm f}({\bm x}_0))|+ 2C_2(1+C_3\|{\bm f}\|_k)\max\{1,\|{\bm z}\|\}\max\{1,\|{\bm z}_0\|\} )L_{p_2}({\bm z},{\bm z}_0)\\
 &\overset{(\ref{eq:Lp})}{\leq} &  
 \Big(|c''_2 ({\bm z}_0,{\bm f}({\bm x}_0))|+ 2C_2(1+C_3\|{\bm f}\|_k)\Big)
 \max\{1,\|{\bm z}_0\|\}^{p_2} \max\{1,\|{\bm z}\|\}^{p_2}\\
 &=:& \widetilde{C}_{\bm f}\max\{1,\|{\bm z}\|\}^{p_2},
  \edeq
 where $\widetilde{C}_{\bm f}:=\Big(|c''_2 ({\bm z}_0,{\bm f}({\bm x}_0))|+ 2C_2(1+C_3\|{\bm f}\|_k)\Big)
 \max\{1,\|{\bm z}_0\|\}^{p_2}$,
 and thus 
 \bgeq
|c''_2 ({\bm z},{\bm f}({\bm x}))|\|T_{\bm x}{\bm h}\|_k= |c''_2 ({\bm z},{\bm f}({\bm x}))\langle {\bm h}, k_{\bm x}\rangle|\| k_{\bm x}\|_k
 \leq |c''_2 ({\bm z},{\bm f}({\bm x}))| \| k_{\bm x}\|_k^2 \overset{(\ref{eq:kx_norm})}{\leq}  \widehat{M}_1^2 \widetilde{C}_{\bm f} \max\{1,\|{\bm z}\|\}^{p_2+2}.
\edeq
Then $\bbe_Q[|c''_2 ({\bm z},{\bm f}({\bm x}))\langle {\bm h}, k_{\bm x}\rangle|\|k_{\bm x}\|_k]<\infty$ for all $Q\in {\cal M}_Z^{p_2+2}$.
 Let 
 $$
\widetilde{\mathcal{P}}_{\bm h}=\{Q\in \mathscr{P}(Z):\bbe_{Q}[|c'_2({\bm z},{\bm f}({\bm x}))|\|k_{\bm x}\|_k]<\infty, \bbe_{Q}[|c''_2 ({\bm z},{\bm f}({\bm x}))\langle {\bm h}, k_{\bm x}\rangle|\|k_{\bm x}\|_k]<\infty\}.
$$
Combining Remark~\ref{rem:well-defined}(iv)-(v) that,
under
${\rm (C4)}'$ and ${\rm (K1)'}$,
$\bbe_Q[|c'_2 ({\bm z},{\bm f}({\bm x}))|\|k_{\bm x}\|_k]$ is well-defined for $Q\in {\cal M}_Z^{{p}_1+1}$.
Then 
$\widetilde{\mathcal{P}}\subset {\cal M}_Z^{\max\{p_1+1,p_2+2\}}$.

Let 
$
\widetilde{\mathfrak{R}}_t({\bm z};{\bm h}):=t^{-1} [c'_2 ({\bm z},({\bm f}+t{\bm h})({\bm x})){\color{black}k_{\bm x}}-c'_2 ({\bm z},{\bm f}({\bm x})){\color{black}k_{\bm x}}].
$
Then
$$
\|\widetilde{\mathfrak{R}}_t({\bm z};{\bm h})\|_k
\leq \Psi({\bm z}) \|{\color{black}k_{\bm x}}\|_k |{\bm h}({\bm x})|
\leq C_2L_p({\bm z},{\bm z})\|{\bm h}\|_k\|{\color{black}k_{\bm x}}\|_k^2.
$$
Analogous to 
Proposition~\ref{pro:derivative-objective}, we can obtain  by the Lebesgue Dominated Convergence Theorem
\bgeq
&&D_{\bm f}(\bbe_P[c'_2 ({\bm z},{\bm f}({\bm x})){\color{black}k_{\bm x}}])({\bm h})\\
&=&\bbe_P\left[\lim_{t\downarrow 0}\frac{c'_2 ({\bm z},({\bm f}+t{\bm h})({\bm x})){\color{black}k_{\bm x}}- c'_2 ({\bm z},{\bm f}({\bm x})){\color{black}k_{\bm x}}}{t}\right] \\
&=& \bbe_P\left[\lim_{t\downarrow 0}\frac{c'_2 ({\bm z},\langle({\bm f}+t{\bm h}),{\color{black}k_{\bm x}}\rangle )- c'_2 ({\bm z},\langle{\bm f},{\color{black}k_{\bm x}}\rangle)}{t}{\color{black}k_{\bm x}}\right]\\
&=& \bbe_{P}[c''_2 ({\bm z},{\bm f}({\bm x})) \langle {\bm h}, {\color{black}k_{\bm x}}\rangle {\color{black}k_{\bm x}}],
\edeq
for all $P\in \widetilde{\mathcal{P}}$.
This shows
$D_{{\bm f}}(\bbe_{P}[c'_2 ({\bm z},{\bm f}({\bm x})){\color{black}k_{\bm x}}])=\bbe_{P}[c''_2 ({\bm z},{\bm f}({\bm x}))T_{{\bm x}}]$, 
where $T_{{\bm x}}$ is defined as in
(\ref{eq:tx1}).
\hfill
\Box

\section{Proof of Proposition~\ref{prop:Lip-c'-c''}}
\label{sec:proof_Lip2d}

Part (i). 
Under ${\rm (K1)'}$,
\bgeqn
\label{eq:kx1}
\|{k}_{{\bm x}_1}\|_k\leq \|k_{{\bm x}_0}\|_k+C_3\|{\bm x}_1-{\bm x}_0\|
\leq \widehat{M}_1\max\{1,\|{\bm z}_1\|\}.
\edeqn
Under ${\rm (C3)}'$, ${\rm (C4)}'$ and ${\rm (K1)'}$,
we have 
$
| c'_2 ({\bm z}_2,{\bm f}({\bm x}_2))|
\leq \widehat{M}_2 \max\{1,\|{\bm z}_2\|\}^{p_1}.
$
By ${\rm (K1)'}$,
${\rm (C4)}'$, 
and inequality (\ref{eq:f-LLip}),
we obtain that 
\bgeq
&& |\langle {\bm d}, c'_2 ({\bm z}_1,{\bm f}({\bm x}_1))k_{{\bm x}_1} \rangle-\langle {\bm d}, c'_2 ({\bm z}_2,{\bm f}({\bm x}_2))k_{{\bm x}_2} \rangle| \\ 
&\leq & \|{\bm d}\|_k \| c'_2 ({\bm z}_1,{\bm f}({\bm x}_1))k_{{\bm x}_1}
-c'_2 ({\bm z}_2,{\bm f}({\bm x}_2))k_{{\bm x}_2}\|_k \nonumber \\
&\leq & \|{\bm d}\|_k \Big( \| c'_2 ({\bm z}_1,{\bm f}({\bm x}_1))k_{{\bm x}_1}
-c'_2 ({\bm z}_2,{\bm f}({\bm x}_2))k_{{\bm x}_1}\|_k
+ \| c'_2 ({\bm z}_2,{\bm f}({\bm x}_2))k_{{\bm x}_1}
-c'_2 ({\bm z}_2,{\bm f}({\bm x}_2))k_{{\bm x}_2}\|_k \Big)\nonumber \\
&\leq &\|{\bm d}\|_k \Big( C_1 L_{p_1}({\bm z}_1,{\bm z}_2) \|k_{{\bm x}_1}\|_k\|{\bm z}_1-{\bm z}_2\|
+ | c'_2 ({\bm z}_2,{\bm f}({\bm x}_2))| C_3\| {\bm x}_1-{\bm x}_2\| \Big)\nonumber \\
&=&\|{\bm d}\|_k\Big( C_1 \widehat{M}_1\max\{1,\|{\bm z}_1\|\} L_{p_1}({\bm z}_1,{\bm z}_2)\|{\bm z}_1-{\bm z}_2\| 
+ C_3\widehat{M}_2\max\{1,\|{\bm z}\|_2\}^{p_1}\|{\bm z}_1-{\bm z}_2\|\Big)  \nonumber \\
&\leq & (C_1\widehat{M}_1+C_3\widehat{M}_2)
 L_{p_1+1}({\bm z}_1,{\bm z}_2)\|{\bm z}_1-{\bm z}_2\| \nonumber \\
&=: &C_{{\bm f}_0} L_{{p}_1+1}({\bm z}_1,{\bm z}_2)\|{\bm z}_1-{\bm z}_2\|, 
\nonumber 
\edeq
where
$C_{{\bm f}_0}:=(C_1\widehat{M}_1+C_3\widehat{M}_2)$,
 $\widehat{M}_1$ and $\widehat{M}_2$ are defined as in Remark~\ref{rem:well-defined}~(ii) and (iii),
$C_1$ and $C_3$ are defined as in Assumption~\ref{assu-solution-P}${\rm (C4)}'$ and  Assumption~\ref{ass:assumption-Z-k}${\rm (K1)}'$.
Moreover, by Assumption~\ref{assu-solution-P}${\rm (C4)'}$,
\bgeqn
|c'_2 ({\bm z},{\bm f}({\bm x}))|
&\leq &
|c'_2 ({\bm z}_0,{\bm f}_0({\bm x}_0))|+
C_1 L_{p_1}({\bm z}_0,{\bm z}_0) |{\bm f}({\bm x}_0)-{\bm f}_0({\bm x}_0)|\nonumber\\
&\leq & 
|c'_2 ({\bm z}_0,{\bm f}_0({\bm x}_0))|+C_1 L_{p_1}({\bm z}_0,{\bm z}_0) |\langle{\bm f}-{\bm f}_0, k_{{\bm x}_0}\rangle |\nonumber\\
&\leq & | c'_2 ({\bm z}_0,{\bm f}_0({\bm x}_0))|+C_1 L_{p_1}({\bm z}_0,{\bm z}_0)\|k_{{\bm x}_0}\|_k\epsilon_{\cal V}, 
\forall {\bm f}\in {\cal V}_{{\bm f}_0}.
\label{eq:cyf}
\edeqn

Part (ii).
Let ${\bm d}_1,{\bm d}_2\in {\cal H}_k$ with $\|{\bm d}_1\|_k\leq 1$, $\|{\bm d}_2\|_k\leq 1$.
Note that $|{\bm d}_i({\bm x})|=|\langle {\bm d}_i, {\color{black}k_{\bm x}}\rangle|\leq \|{\bm d}_i\|_k\|{\color{black}k_{\bm x}}\|_k$.
By \eqref{eq:tx1} and \eqref{eq:tx2}, 
 we have 
\bgeq
&& |\langle \big(c''_2 ({\bm z}_1,{\bm f}({\bm x}_1))T_{{\bm x}_1}- c''_2 ({\bm z}_2,{\bm f}({\bm x}_2))T_{{\bm x}_2}\big){\bm d}_1,{\bm d}_2\rangle|  \\
&\leq & |c''_2 ({\bm z}_1,{\bm f}({\bm x}_1)){\bm d}_1({\bm x}_1){\bm d}_2({\bm x}_1)
-c''_2 ({\bm z}_1,{\bm f}({\bm x}_1)){\bm d}_1({\bm x}_2){\bm d}_2({\bm x}_2)|\\
&& +|c''_2 ({\bm z}_1,{\bm f}({\bm x}_1)){\bm d}_1({\bm x}_2){\bm d}({\bm x}_2)
-c''_2 ({\bm z}_2,{\bm f}({\bm x}_2)){\bm d}_1({\bm x}_2){\bm d}_2({\bm x}_2)|\\
&\leq & |c''_2 ({\bm z}_1,{\bm f}({\bm x}_1))|
\|{\bm d}_1\|_k\|{\bm d}_2\|_k
\left(\|k_{{\bm x}_1}\|_k+\|k_{{\bm x}_2}\|_k \right)\|k_{{\bm x}_1}-k_{{\bm x}_2}\|_k\\
&&
+\|k_{{\bm x}_2}\|^2_k C_2L_{p_2}({\bm z}_1,{\bm z}_2) (\|{\bm z}_1-{\bm z}_2\|
+\|{\bm f}\|_k\|k_{{\bm x}_1}-k_{{\bm x}_2}\|_k)\\
&\leq &\sup_{{\bm f}\in {\cal V}_{{\bm f}_0}}  \Big(\widetilde{C}_{\bm f}  \max\{1,\|{\bm z}_1\|\}^{p_2}\Big)\Big(2 \widehat{M}_1 \max\{1,\|{\bm x}_1\|,\|{\bm x}_2\|\} \Big) C_3  \|{\bm x}_1-{\bm x}_2\|\\
&& +\widehat{M}_1^2\max\{1,\|{\bm x}_2\|\}^{2}
C_2 L_{p_2}({\bm z}_1,{\bm z}_2) (\|{\bm z}_1-{\bm z}_2\|+\|{\bm f}\|_k C_3\|{\bm x}_1-{\bm x}_2\|)\\
&\leq &  \widehat{C}_{{\bm f}_0} L_{p_2+2}({\bm z}_1,{\bm z}_2)\|{\bm z}_1-{\bm z}_2\|,
\edeq
where 
$
\widehat{C}_{{\bm f}_0}:= \sup_{{\bm f}\in {\cal V}_{{\bm f}_0}} \widetilde{C}_{\bm f} 2\widehat{M}_1 C_3 +\widehat{M}_1^2 C_2(1+\|{\bm f}\|_kC_3).
$ 
Analogous to the derivation of \eqref{eq:cyf}, 
we can show that $|c''_2 ({\bm z}_0,{\bm f}({\bm x}_0))|$ and 
hence $\widetilde{C}_{{\bm f}}$ are uniformly bounded over ${\cal V}_{{\bm f}_0}$. 
By 
 the definition of 
$\widehat{C}_{{\bm f}_0}$,  $\widehat{C}_{{\bm f}_0}$is bounded.
\hfill 
\Box

\section{Implicit function theorem}
\label{app:Implicit}

\begin{theorem}[Implicit Function Theorem for Strictly Monotone Functions]
\label{thm:implicit_thm}
Let $U$ be a
Banach space
equipped with norm $\|\cdot\|$ and $V$ be a Hilbert space. 
Consider a function $\psi: U \times V \to V$ and 
a point $(\bar{u},\bar{v}) \in \inmat{int dom} \psi$ satisfying
$\psi( \bar{u}, \bar{v}) = 0$. 
Assume: (a) there are neighborhoods ${\cal B}_U$ of $\bar{u}$ and ${\cal B}_V$ of $\bar{v}$ such that $\psi$ is
continuous on ${\cal B}_U \times {\cal B}_V$, 
(b) $\psi(u, \cdot)$
is strongly monotone on ${\cal B}_V$ uniformly in $u \in  {\cal B}_U$ and
(c)
$\psi(\cdot,v)$ is Lipschitz continuous on
$U$ uniformly in $v \in {\cal B}_V$.
Then following assertions hold.
\begin{itemize}
\item[(i)] 
The solution mapping
\bgeq 
S(u) =\{v\in V: \psi(u,v)=0\}
\edeq
for $u \in U$
has a single-valued localization around $\bar{u}$ for $v$.

\item[(ii)]
$S(\cdot)$ is Lipschitz continuous 
around $\bar{u}$ for $\bar{v}$.

\end{itemize}

\end{theorem}

\noindent
\textbf{Proof.} The proof is similar to \cite[Theorem 1H.3]{DoR09} which is in finite dimensional space.
Here we give a sketch for completeness particularly 
because we will use Shauder's fixed point theorem as opposed to Brower's fixed point theorem in the proof of \cite[Theorem 1H.3]{DoR09}.

Part (i). Observe that if $u\in \inmat{dom} S\cap {\cal B}_U$, then $S(u)\cap {\cal B}_V$ consists of one element, if any. 
In fact, if there existed two elements $v_1,v_2\in 
  S(u)\cap {\cal B}_V$ with $v_1\neq v_2$,
  then from the
strong monotonicity, we would have 
\bgeq 
0=\langle \psi(u,v_1)-\psi(u,v_2),v_1-v_2\rangle >0,
\edeq
which leads to a contradiction.
Thus, all we need is to establish that 
$\inmat{dom} S$ 
contains a neighborhood 
of $\bar{u}$.
Without loss of generality, 
let $\bar{v}= 0$ and choose
$\delta>0$ such that 
$\delta \mathbb{B} \subset {\cal B}_V$. 
For $\epsilon \in (0,\delta]$  define
\bgeqn
d(\epsilon) := \inf_{\epsilon\leq \|v\|\leq \delta}
\frac{\langle v, \psi(\bar{u},v)\rangle}{\|v\|}.
\edeqn
Choose $v\in V$ with $\|v\|\in [\epsilon,\delta]$.
Since $\psi(\bar{u},0)=0$, 
from strict monotonicity of
$\psi(\bar{u},\cdot)$, we obtain that
$$
\langle \psi(\bar{u},v),v\rangle =\langle \psi(\bar{u},v)-\psi(\bar{u},0),v-0\rangle
>0,
$$                      
hence $d(\epsilon)\geq 0$ for all $\epsilon\in (0,\delta]$. We claim that $d(\epsilon)>0$. Assume for the sake of a contradiction that there exists $\epsilon_0>0$ such $d(\epsilon_0)=0$. Then there exists a sequence $\{v_k\}$ with $\|v_k\|\in [\epsilon_0,\delta]$ such that 
$$
\frac{\langle v_k, \psi(\bar{u},v_k)\rangle}{\|v_k\|}\to 0.
$$
On the other hand, it follows by the strong monotonicity of $\psi(\bar{u},\cdot)$
$$
\frac{\langle v_k, \psi(\bar{u},v_k)\rangle}{\|v_k\|}
=\frac{\langle v_k-0, \psi(\bar{u},v_k)-\psi(\bar{u},0)\rangle}{\|v_k\|}\geq \tau\|v_k\|\geq \tau\epsilon.
$$

Part (ii). We proceed in three steps.

\noindent
\underline{Step 1}.
Let $\mu\in (0,d(\delta))$ and $\nu>0$ be such that ${B}(\bar{u},\nu) \subset {\cal B}_U$.
From the continuity of $\psi$ in $u\in {B}(\bar{u},\nu)$ uniformly with respect to $v\in \{v: \|v\|=\delta\}$ there exists $\alpha_\mu\in (0,\nu)$ such that for any $u\in {B}(\bar{u},{\alpha_\mu})$ and $v\in \{v: \|v\|=\delta\}$,
$\|\psi(u,v)-\psi(\bar{u},v)\|\leq \mu$.
Consequently
\bgeqn
\frac{\langle v, \psi(u,v)\rangle}{\|v\|}
=
\frac{\langle v, \psi(\bar{u},v)\rangle}{\|v\|}
+
\frac{\langle v, \psi(u,v)-\psi(\bar{u},v)\rangle}{\|v\|}
\geq \frac{\langle v, \psi(\bar{u},v)\rangle}{\|v\|}-\mu \geq d(\delta)-\mu>0.
\label{eq:IMT-FNT-THM-PROOF}
\edeqn

\noindent
\underline{Step 2}.
For fixed $u \in {B}(\bar{u},\nu)$, consider the function
\bgeqn 
\Psi(v) := \Pi_{\delta\mathbb{B}}(v-\psi(u,v)), \inmat{for}\; v\in \delta\mathbb{B} 
\edeqn
where $\Pi_{\delta\mathbb{B}}$ is a orthogonal 
project in $V$. We use Shauder's fixed point theorem 
to show that $\Phi(\cdot)$ has a fixed point in the interior of $\delta\mathbb{B}$. To this end, we need to show that $\Psi: \delta\mathbb{B} \to \delta\mathbb{B}$ is a compact operator.
By Arzela-Ascoli theorem, it suffices to what that 
(a) $\Pi_{\delta\mathbb{B}}(v-\psi(u,v))$ is uniformly bounded, i.e.,
$$
\sup_{v\in \delta \mathbb{B}}
\|\Pi_{\delta\mathbb{B}}(v-\psi(u,v))\|\leq \delta,
$$
which is obvious and
(b)
$\Pi_{\delta\mathbb{B}}(v-\psi(u,v))$ is equicontinuous, i.e., for any $\epsilon>0$, there exists a $\beta>0$ such that
\bgeqn
\sup_{ v\in\delta\mathbb{B}}
\|\Pi_{\delta\mathbb{B}}(v-\psi(u',v))
- \Pi_{\delta\mathbb{B}}(v-\psi(u'',v))\| \leq \epsilon,
\quad \forall u',u'' \in {B}(\bar{u},\nu), \|u'-u''\|\leq\beta.
\label{eq:equi-cont-Imp-fnt-thm}
\edeqn
Inequality (\ref{eq:equi-cont-Imp-fnt-thm})
is guaranteed by the uniform continuity of $\psi(u,v)$ in $u$ w.r.t. $v\in {\cal B}_V$ because
$$
\|\Pi_{\delta\mathbb{B}}(v-\psi(u',v))
- \Pi_{\delta\mathbb{B}}(v-\psi(u'',v))\|
\leq 
\|\psi(u',v))
- \psi(u'',v))\|.
$$
By Shauder's fixed point theorem, there exists $v^*\in \delta\mathbb{B}$ such that
$$
v^*=\Psi(v^*)= \Pi_{\delta\mathbb{B}}(v^*-\psi(u,v^*)).
$$
Next, we show that $\|v^*\|<\delta$.
Assume for the sake of contradiction 
$\|v^*\|=\delta$. By the definition of 
orthogonal projection
$
\langle 0-v^*, v^*-\psi(u,v^*)-v^*\rangle \leq 0, 
$
which implies
$
\frac{\langle v^*, \psi(u,v^*)\rangle}{\|v^*\|} \leq 0. 
$
On the other hand, it follows by (\ref{eq:IMT-FNT-THM-PROOF})
\bgeqn
\frac{\langle v^*, \psi(u,v^*)\rangle}{\|v^*\|}
\geq d(\delta)-\mu>0,
\label{eq:IMT-FNT-THM-PROOF-0}
\edeqn
a contradiction. Thus
$$
v^* = \Pi_{\delta\mathbb{B}}(v^*-\psi(u,v^*))
= v^*-\psi(u,v^*)
$$
and hence $\psi(u,v^*)=0$.

\noindent 
\underline{Step 3}.
To prove continuity of $S(\cdot)$ at $\bar{u}$, we note that by the strong monotonicity of $\psi(u,\cdot)$, 
\bgeq
\|S(u')-S(u'')\|\| \psi(u', S(u'))-
\psi(u', S(u''))\| &\geq& \langle S(u')-S(u''), \psi(u', S(u'))-
\psi(u', S(u''))\rangle \\
&\geq& \tau \|S(u')-S(u'')\|^2.
\edeq 
Since $\psi(u', S(u')) = \psi(u'', S(u''))=0$, then
by the Lipschitz continuity of $\psi(\cdot,v)$,
\bgeq
\| \psi(u', S(u'))-
\psi(u', S(u''))\|  &=&
\| \psi(u'', S(u''))-
\psi(u', S(u''))\| \\
&\leq& 
L \|u'-u''\|.
\edeq
Combining the two inequalities above, we obtain
\bgeq
\|S(u')-S(u'')\| \leq \frac{L}{\tau}\|u'-u''\|. 
\edeq
The proof is complete. \hfill $\Box$

\end{appendix}

\end{document}